\documentclass[10pt]{amsart} 

\usepackage{a4wide}
\usepackage[T1]{fontenc}
\usepackage[utf8]{inputenc}
\usepackage[english]{babel}
\usepackage[pdfauthor={Lorenzo Mantovani}, %
                pdftitle={Localizations and completions of stable $\infty$-categories}%
            dvips,colorlinks=true]{hyperref}
\hypersetup{
    bookmarksnumbered=true,
    linkcolor=red,
    citecolor=blue,
}
\usepackage{amsmath,amsthm,stmaryrd}
\usepackage{mathrsfs,amssymb,amsfonts}
\usepackage[mathscr]{eucal}
\usepackage{mathcomp}
\usepackage[all,cmtip]{xy}
\usepackage{enumitem}
\usepackage{xcolor}
\usepackage{verbatim}
\usepackage[style=alphabetic, doi=false, url=false]{biblatex}
\usepackage[switch, displaymath, mathlines,right,left,pagewise]{lineno}

\newcommand{\op}{\operatorname}

\newcommand{\spt}{\mathbf{Spt}}

\newcommand{\sm}[1]{\mathbf{Sm}_{#1}}

\newcommand{\bb}{\mathbb}
\renewcommand{\bf}{\mathbf}
\newcommand{\scr}{\mathscr}
\newcommand{\cal}{\mathcal}

\renewcommand{\rm}{\mathrm}

\newcommand{\M}{\scr{M}}
\DeclareMathOperator\Hom{Hom}

\newcommand{\Map}{\mathrm{Map}}

\DeclareMathOperator\spec{Spec}

\newcommand{\coker}{\text{coker}}
\newcommand{\ra}{\longrightarrow}
\newcommand{\bc}[1]{\big <#1 \big >}

\DeclareMathOperator\cofib{cofib}
\DeclareMathOperator\fib{fib}
\DeclareMathOperator\colim{colim}
\renewcommand{\lim}{\rm{lim}}

\DeclareRobustCommand{\longlongrightarrow}{\relbar\joinrel \relbar\joinrel\relbar\joinrel\rightarrow}
\newcommand{\rra}{\longlongrightarrow}

\DeclareRobustCommand{\longlonglongrightarrow}{\relbar\joinrel \relbar\joinrel \relbar\joinrel\relbar\joinrel\rightarrow}
\newcommand{\rrra}{\longlonglongrightarrow}

\newcommand{\medslant}[2]{{\raisebox{.15em}{$#1$}\left/\raisebox{-.15em}{$#2$}\right.}}

\theoremstyle{plain}
\newtheorem{thm}[subsubsection]{Theorem}
\newtheorem{prop}[subsubsection]{Proposition}
\newtheorem{lemma}[subsubsection]{Lemma}
\newtheorem{cor}[subsubsection]{Corollary}

\theoremstyle{definition}
\newtheorem{defin}[subsubsection]{Definition}
\newtheorem{ass}[subsubsection]{Assumption}
\newtheorem{notat}[subsubsection]{Notation}

\theoremstyle{remark}
\newtheorem{rmk}[subsubsection]{Remark}
\newtheorem{ex}[subsubsection]{Example}

\DeclareRobustCommand{\gobblefive}[5]{}

\title{Localizations and completions of stable $\infty$-categories}
\author{Lorenzo Mantovani}
\address{Department of Mathematics, Universit\"at Z\"urich, Z\"urich, Switzerland}
\email{lorenzo.mantovani@math.uzh.ch}
\begin{document}

\begin{abstract}
  We extend some classical results of Bousfield on homology localizations and nilpotent completions to a presentably symmetric monoidal stable $\infty$-category $\scr{M}$ admitting a multiplicative left-complete $t$-structure. If $E$ is a homotopy commutative algebra in $\scr{M}$ we show that $E$-nilpotent completion, $E$-localization, and a suitable formal completion agree on bounded below objects when $E$ satisfies some reasonable conditions.  
\end{abstract}

\maketitle

\setcounter{tocdepth}{1}
\tableofcontents

\section{Introduction} 
  \label{sec:introduction}

  \subsubsection*{Motivation}
  Let $\scr M$ be a presentably symmetric monoidal $\infty$-category, with monoidal product $\wedge$ and unit $\bf 1$, and let $E$ be an object of $\scr M$. One can thus construct a homology localisation of $\scr M$ by inverting all the maps $\phi$ in $\scr M$ such that $\phi\wedge E$ is an equivalence. This construction was first introduced for the topological stable category $\scr{SH}$ in \cite{MR551009}. The associated localization functor $X\mapsto X_E$ is called $E$-homology localization, and has a particularly simple universal property: $X_E$ is the initial $E$-local object with a map $\lambda(X): X\ra X_E$. Unfortunately very little can be said about $X_E$ for a general $E$, even in the case $X=\bf 1$.

  Assume now that $E$ is a commutative algebra object in $\scr M$. In this situation, given any object $ X \in \scr{M}$, we can perform a second construction $X^\wedge_E$, by setting
  \[X^\wedge_E:=\lim_{\Delta} X\wedge E^{\wedge\bullet}.\]
  By construction there is a natural map $\alpha(X): X \ra X^\wedge_E$, which factors through $\lambda$. The object $X^\wedge_E$ is called \emph{nilpotent completion of $X$ at $E$}.

  Instances of this construction have appeared in a wide range of contexts. When $\scr M$ is (the nerve of) the category of modules over a ring $A$, and $E$ is a commutative $A$-alegbra, the map $\alpha(X)$ can be interpreted as the obstruction to recover $X$ from its associated descent datum along $A \ra E$. In algebraic topology nilpotent completions where used by Adams and many others in relation to computing homotopy groups of spectra. Indeed starting with the simplicial object 
  \[X\wedge E^{\wedge\bullet}\]
  one can construct the tower of its partial totalizations and the inverse limit of such a tower recovers $X^\wedge_E$. In addition the tower of partial totalizations gives rise to a Bousfield-Kan spectral sequence conditionally converging to $X^\wedge_E$. In some particularly favourable situation the page of this spectral sequence is amenable to computations, and sometimes a good deal of information on $X^\wedge_E$ can be understood via this spectral sequence.

  Here is a crucial fact that follows from combining several parts of \cite{MR551009}.
  \begin{thm}[Bousfield]
    \label{thm:}
    When $E$ is $(-1)$-connected commutative algebra in $\scr{SH}$, with $\pi_0(E) \simeq \bb Z/n\bb Z$, then for every $k$-connected spectrum $X$ the natural map $X^\wedge_E \ra X_E$ is an equivalence. Furthermore, under the above assumptions, $X_E$ is naturally identified with the derived completion $X^\wedge_n=\lim_k X/n^k$.
  \end{thm}
  
  \subsubsection*{Actual Content} In this paper we axiomatize some of the techniques used by Bousfield and adapt them to work in a presentably symmetric monoidal stable $\infty$-category $\scr{M}$. Our aim is to reach a formal analogue for $\scr{M}$ of the above theorem of Bousfield. The main assumption we need on $\scr{M}$ is that it comes endowed with an accessible $t$-structure which has the following properties:
  \begin{itemize}
    \item $\bf 1 \in \scr{M}_{\geq0}$;
    \item $\scr {M}_{\geq p}\wedge\scr{M}_{\geq q}\subseteq\scr{M}_{\geq p+q}$;
    \item the $t$-structure is left-complete; i.e. for every $X \in \scr{M}$, $X\simeq \lim_n P^n(X)$, where $P^n(-)$ denotes the $\leq n$ truncation functor.
  \end{itemize}
  The main application we have in mind being motivic homotopy theory, we have decided to dedicate Section \ref{sec:a_reminder_on_motivic_categories} to recollecting some well known facts about the motivic stable category $\scr{SH}(S)$ and about the categories of modules $\scr{Mod}_A(S)$ over a commutative algebra $A \in \scr{SH}(S)_{\geq0}$. In particular we review how the $t$-structure $\scr{Mod}_A(S)$ inherits from $\scr{SH}(S)$ has the above three properties when $S$ is a Noetherian scheme of finite Krull dimension.

  In oder to work with an abstract symmetric monoidal $\infty$-category $\scr{M}$ we have chosen to axiomatize the elements of $\pi_0(\bb S)\simeq \bb Z$ in terms of maps $L\ra \bf 1$ where $L$ is a $\wedge$-invertible object of $\scr{M}$ such that $L\wedge-$ respects both $\scr{M}_{\geq 0}$ and $\scr{M}_{\leq 0}$; objects satisfying this properties are called tif objects (tif stands for "tensor invertible and flat").

  Another choice we have made is to work with localizations at homotopy commutative algebras of $\scr{M}$, i.e. with a commutative algebra of the homotopy category $h\scr{M}$. We refer to \ref{subsub:hom_algebras} for a clearer definition.
  
  In this framework the main assumptions on $E$ is essentially the following. $E$ is a homotopy commutative algebra of $\scr{M}_{\geq 0}$. Furthermore there exists a finite set of tif objects $\{L_i\}_{i=1}^r$ and maps $f_i: L_i \ra \bf 1$ such that the unit $\bf 1 \ra E$ induces an isomorphism $\tau_0(\bf 1)/(f_1,\dots,f_r)\simeq \tau_0(E)$. For a more precise statement of this technical assumption we direct the reader to \ref{sub:assumption_A1}.

  The main results we obtain in this general framework are then condensed in the following.
  \begin{thm}[\ref{thm:red_to_moore_spt} and \ref{prop:localization_at_mod_x_1-x_n_moore_spectrum}]
    \label{thm:1}
    Let $E$ be a homotopy commutative algebra in $\scr{M}$ satisfying assumption \ref{sub:assumption_A1} in the special case of $J=\emptyset$. Then for every $k$-connected object $X$ in $\scr{M}$ there is a canonical isomorphism $X_E\simeq X^{\wedge}_{f_1} {\cdots}^\wedge _{f_r}$ compatible with the localization map $\lambda_E(X): X\ra X_E$ and the formal completion map $\chi_{\underline f}(X): X \ra X^{\wedge}_{f_1} {\cdots}^\wedge _{f_r}$.
  \end{thm}

  The proof combines two main steps, where one compares formal completions and homology localization with some other intermediate object. This involves a sort of axiomatization of Moore spectra which is performed in Section \ref{sec:moore_spectra}. In the context of \ref{thm:1} the relevant Moore object is $M=C(f_1)\wedge\dots\wedge C(f_r)$, where $C(f_i)$ denotes the cofiber of $f_i$. The construction of $M$ is what dictates the rather strong assumptions we have on $\tau_0(E)$. With this notation the first main step consists in proving the well known statement that for every $X \in \scr{M}$, $X_M\simeq X^{\wedge}_{f_1} {\cdots}^\wedge _{f_r}$. The second step, performed in Section \ref{sec:localizations}, consists in showing that when $X$ is $k$-connected $X_M\simeq X_E$. This is done with a careful use of Bousfield classes, and uses crucially the left-completeness of the $t$-structure.

  Section \ref{sec:examples} contains a list of relevant examples and applications in the motivic setting. We mention here that, as an application, we partially recover a conservativity result for motives of Bachmann (c.f. \cite{Bacmann_Conservativity}). For this we work with $\scr{M}=\scr{SH}(K)$ where $K$ is a field and $p$ is the exponential characteristic of $K$. We have thus functors $M: \scr{SH}(K)\ra\scr{DM}(K)$ and $\widetilde M: \scr{SH}(K)\ra \widetilde{\scr{DM}}(K)$ associating with every spectrum its motive $M(X)$ and its Chow-Witt motive $\widetilde M(X)$.
  \begin{cor}[\ref{cor:conserv}]
    \label{thm:cons}
    Let $K$ be a field of exponential characteristic $p\not =2$, and let $X$ be a $k$-connected spectrum for some $k$. 
    \begin{enumerate}
      \item Assume that $-1$ is a sum of squares in $K$. If $M(X)=0$ then $X[\frac{1}{2}]=0$.
      \item Assume that $K$ has finite \'etale $2$-cohomological dimension and that $X$ is a dualizable $\bb S[\frac{1}{p}]$-module. If $M(X)=0$ then $X=0$.
      \item Assume that $X$ is a $\bb S[\frac{1}{p}]$-module. If $\widetilde M(X)=0$ then $X=0$.
    \end{enumerate}
  \end{cor}
  Actually analogous results hold for categories of motives that arise as categories of modules over a homotopy commutative algebra $E$ in $\scr{SH}(K)$ whose $\tau_0(E)$ is Milnor (or Milnor-Witt) $K$-theory. We direct the reader to \ref{rmk:our_cons_res_vs_Bachmann} for a precise description of the relation with the work of Bachmann.

  Section \ref{sec:the_E_based_MANSS} contains a construction of $E$-nilpotent completions using an axiomatized version of the Adams tower, and a construction of the associated spectral sequence. The Adams tower is an alternative to the tower of partial totalizations mentioned above, and allows to avoid the assumption that $E$ be a commutative algebra of $\scr{M}$. 

  Section \ref{sec:nilpotent_resolutions} contains the axiomatization of nilpotent resolutions and a general proof of some of their properties. Using these we can provide a universal property for the Adams tower as a pro-object. It is by comparing the pro-objects associated with the Adams tower and with the formal completion that we obtain the following result.

  \begin{thm}[\ref{thm:convergence_mod_stuff} and \ref{thm:convergence_inv_stuff}]
     \label{thm:2}
     Let $E$ be an homotopy commutative algebra in $\scr{M}$ satisfying assumption \ref{sub:assumption_A1} in the special case where either $J=\emptyset$ or $I=\emptyset$. Then for every $k$-connected object $X$ in $\scr{M}$ the natural map $X^{\wedge}_E \ra X_E$ is an equivalence in $\scr{M}$.
   \end{thm} 

  This result has already appeared in \cite{MR2811716}, although in a very special case, and with a different proof. It was the reading of this specific work that stimulated our interest in the topic. Our approach to the problem is in fact very different in spirit from that of \cite{MR2811716}.

  In their recent work \cite{topologmodels} the authors have generalized and streamlined the arguments of a previous version of our results which appeared in \cite{locandcomp}. Although our results of \cite{locandcomp} are phrased for the motivic stable homotopy category $h\scr{SH}(K)=\cal{SH}(K)$ of a perfect field $K$, the structure and the arguments of the present paper are essentially the same as those of \cite{locandcomp}. As a consequence the present paper and the second section of \cite{topologmodels} present similar results with similar techniques, although \cite{topologmodels} has a more direct and a simpler approach. The present paper was written, independently of \cite{topologmodels}, during the spring of 2020. We wish to thank T. Bachmann and P. A. Østvær for having allowed us to publish our work despite the overlap with theirs.

  \subsubsection*{Acknowledgments} I wish to thank M. N. Levine, J. Ayoub and A. Kresch for their encouragement in editing the arguments of \cite{locandcomp}. A special thank goes to T. Bachmann for the interest in a previous version of this paper. Another special thank goes to J. I. Kylling. His enthusiasm and interest on this project have helped me focusing on this work in several occasions. Furthermore he found a mistake in a previous version of \cite{locandcomp}, and pointed out a number of typos. I wish to thank the anonymous referee for her/his careful reading, for suggesting further references, and for tangibly improving the overall quality of the exposition. I also wish to thank M. Porta who helped me navigating the higher categorical language. Finally I wish to thank the algebraic geometry group of the University of Zürich for providing a comfortable working environment. This work was supported by the \emph{Swiss National Science Foundation} (SNF), project 200020\_178729.


\section{Preliminaries} 
  \label{sec:a_reminder_on_motivic_categories}

  	We begin this section by introducing the categorical framework that we will be working with in order to fix some ideas and notation. This will happen in \ref{sub:cate_framework}. In \ref{sub:motivic_stable_cateogory} we review some well-known facts about motivic stable categories, and in \ref{sub:homotopy_t_structure} we review Morel's homotopy $t$-structure, which is by far the most important tool we need. In \ref{sub:homology_localizations} we introduce $E$-homology localizations associated to an object $E$, and recall some formal properties of these constructions.  We conclude with a review of the formalism of Bousfield classes in \ref{sub:bousfield_classes}, which turns out to be very useful for keeping track of the mutual relations between various localization functors appearing at the same time.

	\subsection{Categorical Framework}
    \label{sub:cate_framework}
    
    \subsubsection{}
	    \label{subsub:ass_on_t-str}
	    Along this paper we will be working with a presentably symmetric monoidal stable $\infty$-category $\scr{M}$. We will denote the monoidal product by $-\wedge-$ and by $\bf 1$ the unit. We also introduce the symbol $\Map(-,-)$ to denote the mapping space between two objects. Finally we denote by $h\scr{M}$ the homotopy category of $\scr{M}$, and by $[-,-]$ the Hom groups of the homotopy category, so that for every non-negative integer $k$ and every pair of objects $M$ and $N$ in $\scr{M}$ we have
      \[ \pi_k \Map(M,N)\simeq \pi_0\Omega^k \Map(M,N) \simeq \pi_0 \Map(\Sigma^k M,N)\simeq [\Sigma^k M,N].\]

      We will assume that $\scr{M}$ is endowed with an $t$-structure: the objects of $\scr{M}_{\geq n}$ will be referred to as $(n-1)$-connected. We will denote by $P^{n}(-)$ (resp. $P_n(-)$, resp. $P^n_n(-)$) the $\leq n$ (resp. $\geq n$, resp. $=n$) truncation functors. The heart
  	    \[\scr{M}^{\heartsuit}:=\scr{M}_{\geq 0}\cap \scr{M}_{\leq 0}\] 
  	  is thus the nerve of an abelian category, and for every object $X$ of $\scr{M}$, $\tau_k(X)=\Sigma^{-k}P^k_k(X)$ will be referred to as the $k$-th homotopy object of $X$. Most of this paper works under the following list of assumptions on the $t$-structure on $\scr{M}$:
	    \begin{enumerate}
	      \item the $t$-structure is accessible (c.f. Definition 1.1.4.2 of \cite{HA});
	      \item the $t$-structure is left-complete, i.e. $X \simeq \lim_k P^k(X)$ for all $X \in \scr{M}$;
	      \item $\bf 1\in \scr{M}_{\geq 0}$;
	      \item $\scr{M}_{\geq p}\wedge\scr{M}_{\geq q}\subseteq \scr{M}_{\geq p+q}$;
	    \end{enumerate}
      Some sections actually work with less assumptions: for this we direct the reader to the introduction of each section and to the specific statements. The only additional assumption that appears on the $t$-structure of $\scr{M}$ is the following:
      \begin{enumerate}[resume]
        \item the $t$-structure is compatible with filtered colimits, i.e. the truncation functors $P^k$ commute with filtered colimits.
      \end{enumerate}
      This last assumptions is needed, in our opinion, to ensure that inverting homotopy elements is a $t$-exact functor (c.f. \ref{cor:inv_elts_is_t_exact}). This assumption is used only in Proposition \ref{cor:HK^MW_mod_stuff_vs_H_KMW_mod_stuff} when $J\not= \emptyset$. In any case this assumption is not needed for the main theorems of the paper.

    \begin{lemma}
    \label{lemma:tens_prod_right_exact}
      If the $t$-structure on $\scr{M}$ is multiplicative, the monoidal product $-\otimes^\heartsuit-$ induced on $\scr{M}^{\heartsuit}$ is right exact.
    \end{lemma}
    
    \begin{proof}
      Let  
      \[F \overset{f}{\ra} G \overset{g}{\ra} H \ra 0\] be an exact sequence of objects of $\scr{M}^{\heartsuit}$. Let $C:=\cofib(f)$, so that $H\simeq P^0(\cofib(f))$. Thus given any $D \in \scr{M}^{\heartsuit}$ we have an induced fiber sequence
      \[D\wedge F \ra D\wedge G \ra D \wedge C,\]
      and we only need to check that the natural map
      \[P^0(D \wedge C)\ra P^0(D\wedge P^0(C))\]
      is an equivalence. This follows from the fact that $D\wedge_A P_1(C) \in \scr{M}_{\geq 1}$.
    \end{proof}

    \subsubsection{} 
	    \label{ssub:mult_prop_t_str}

	    In this situation $\scr{M}_{\geq 0}$ has a natural structure of presentably symmetric monoidal $\infty$-category induced via restriction along $\scr{M}_{\geq 0}\subseteq \scr{M}$. The inclusion $\scr{M}_{\geq 0}\subseteq \scr{M}$ is a symmetric monoidal functor, while its left adjoint $P_0: \scr{M} \ra \scr{M}_{\geq 0}$ is a map of $\infty$-operads, and hence preserves algebra and module categories. All the claims follow at once form Proposition 2.2.1.1 of \cite{HA}.

	    For every non-negative integer $n$, the functors $P^n(-)$ restrict to symmetric monoidal endofunctors on $\scr{M}_{\geq 0}$; the restrictions of the projections maps $\pi_n$ are monoidal natural transformations. In particular for every integer $k\leq n$, the natural projection map $\pi^n_k: P^n(\bf 1) \ra P^k(\bf 1)$ is an algebra map. Moreover if $X \in \scr{M}_{\geq 0}$ then the $P^k(X)$ has a canonical structure of $P^n(\bf 1)$-module for every $n\geq k$, and this structure is compatible with the projections $\pi_n$ for varying $n$.

	    The inclusion
	    \[\scr{M}^{\heartsuit} \subseteq \scr{M}_{\geq 0}\]
	    induces on the heart of the $t$-structure, the structure of a symmetric monoidal $\infty$-category. In addition the above inclusion is a map of $\infty$-operads and thus respects algebras and modules, while the truncation $P^0: \scr{M}_{\geq 0} \ra \scr{M}^{\heartsuit}$ is actually symmetric monoidal. These assertions follow at once from \cite[Propositions 2.2.1.8 and 2.2.1.9]{HA}. In particular we deduce the usual formulas
	    \[-\otimes^{\heartsuit}-\simeq \tau_0(-\wedge -)\]
	    for the induced symmetric monoidal product, and 
	    \[\bf 1^{\heartsuit
	    }\simeq P^0(\bf 1)\simeq\tau_0(\bf 1)\]
	    for the monoidal unit in $\scr{M}^{\heartsuit}$.

    \subsubsection{}
      \label{subsub:hom_algebras}
      We will also use the notion of \emph{homotopy commutative algebra} in $\scr{M}$. With this expression we mean an object $E$ of $\scr{M}$ together with maps $e:\bf 1 \ra E$ and $\mu: E\wedge E\ra E$ and suitable homotopies, making $(E,e,\mu)$ a commutative monoid in the homotopy category $h\scr{M}$. We will similarly use the notion of homotopy $E$-module, defined in an analogous way. For instance, given a choice of composition $\mu_E\circ (\op{id}_E\wedge_A e_E)$ is only homotopic to the identity of $E$, i.e. there is a $2$-simplex of $\scr{M}$ filling
      \[
      \xymatrix{
      E\ar[rr]\ar[dr]_{\op{id}_E} & & E\ar[dl]^{\mu\circ (\op{id}_E\wedge e)}\\
      & E. & \\
      }\]
    \subsubsection{}
      \label{subsub:trunc_of_hom_mod}
    	On a similar note, if $E$ is a homotopy commutative algebra and $X$ is a homotopy $E$-module, the homotopy objects $\tau_k(X)$ inherit a natural structure of $\tau_0(E)$-modules in $\scr{M}^{\heartsuit}$.      

  \subsection{Motivic Stable Categories} 
    \label{sub:motivic_stable_cateogory}
    
    \subsubsection{} 
    \label{ssub:sh}
    
    A \emph{base scheme} is a Noetherian scheme of finite Krull dimension. Given a base scheme $S$, we denote by $\scr{SH}(S)$ the Morel-Voevodsky stable $\infty$-category. The associated homotopy category $\cal{SH}(K)$ was first constructed in \cite{morel:stabA1ht,MR2175638}, and more generally, over any base scheme $S$, by Jardine \cite{jard:motsym} via model model categories. We recall that $\scr{SH}(S)$ is a presentably symmetric monoidal stable $\infty$-category; its monoidal product is denoted by $\wedge$ and the unit, the motivic sphere spectrum, is denoted by $\bb S$. We redirect the reader to \ref{sec:cat_rec} for a reference to the previous claim and a quick review of the higher categorical terminology.

    The $\infty$ category $\scr {SH}(S)$ is actually compactly generated by the set 
      \[ \{\Sigma ^{p+q\alpha}\Sigma^\infty X_+ : \; X \in \sm S, \; p,q \in \bb Z \}, \]
      where $\Sigma^{p+q\alpha}$ is defined as $\Sigma_{S^1}^{p}\Sigma_{\bb G_m}^{q}$.

    \subsubsection{} 
    \label{ssub:modules_over_ring_spectra}
    
    To any commutative algebra $A$ of $\scr{SH}(S)$ (cf \ref{definB:com_alg}) we associate a category $\scr {Mod}_A(S)$ whose object are called $A$-module spectra, or simply $A$-modules. $\scr{Mod}_A(S)$ inherits from $\scr{SH}(S)$ the property of being a presentably symmetric monoidal stable $\infty$-category. Once again references and definitions are postponed to \ref{ssubB:infty_mods}. The monoidal product is denoted by $-\wedge_A-$, or simply by $-\wedge -$ when no confusion arises; the monoidal unit is denoted by $\bf 1_A$.

    In addition we have a free-forget adjunction 
    \[ F_A: \scr{SH}(S) \rightleftarrows \scr {Mod}_A(S): U_A.\]
    The forgetful functor $U_A$ is right adjoint of $F_A$: it commutes with all small limits and colimits, and it is conservative. The functor $F_A$ is symmetric monoidal and commutes with all small colimits. Since $U_A$ is conservative and commutes with colimits, the category $\scr{Mod}_A(S)$ is compactly generated by the set 
    \[ \{F_A (\Sigma^{p+q\alpha} \Sigma^\infty X_+) : \; X \in \sm S, \; p,q \in \bb Z \}. \] 
    We conclude by observing that the composition $U_A\circ F_A\simeq A\wedge-$, and thus the
    monoidal unit $\bf 1_A \in \scr{Mod}_A(S)$ is mapped to $A$ in $\scr{SH}(S)$. We will abuse the language and confuse $\bf 1_A$ and $A$. We give a reference for all these facts in \ref{ssubB:infty_mods}.

  \subsection{Homotopy t-structure} 
    \label{sub:homotopy_t_structure}

    \subsubsection{} 
    \label{ssub:acc_t_str}
    
	    We define $\scr {Mod}_A(S)_{\geq n}$ as the smallest full sub-$\infty$-category of $\scr{M}_A(S)$ closed under small colimits and extensions that contains the collection 
	    \[\{ F_A(\Sigma ^{p+q\alpha} \Sigma ^\infty X_+) \; : X \in \sm S, \; p \geq n,\; q\in \bb Z \}.\]
	    Furthermore we denote by $\scr{Mod}_A(S)_{\leq n}$ the full sub-$\infty$-category spanned by those $A$-modules $Z$ such that $\Map_A(X,Z)\simeq \ast$ for every $X \in \scr{Mod}_A(S)_{\geq n+1}$. 
	    
	    By Proposition 1.4.4.11 of \cite{HA} the pair of subcategories $\big (\scr{Mod}_A(S)_{\geq 0},\scr{Mod}_A(S)_{\leq -1} \big )$ defines an accessible $t$-structure on $\scr{Mod}_A(S)$ (c.f. Definition 1.4.4.12 of \cite{HA}), called \emph{homotopy $t$-structure}. It follows that for every $n \in \bb Z$ we have a cofiber sequence in $\scr {Mod}_A(S)$
	    \begin{equation}
	    \label{eqn:postfibseq}
	    P_n(X) \overset{\delta_n}{\ra} X \overset{\pi_{n-1}}{\ra} P^{n-1}(X)
	    \end{equation}
	    which is functorial in the $A$-module $X$, with $P_n(X)\in \scr{Mod}_A(S)_{\geq n}$ and $P^{n-1}(X)\in \scr{Mod}_A(S)_{\leq n-1}$. It is clear from the choice of the generators that when $A$ is $(-1)$-connected object for the homotopy $t$-structure on $\scr{SH}(S)$, the homotopy $t$-structure on $\scr{Mod}_A(S)$ is multiplicative.

    \begin{thm}
    \label{thm:morelTstr}
    Let $S$ be a Noetherian scheme of finite Krull dimension and $A$ be a commutative algebra in $\scr{SH}_{\geq 0}$. The homotopy $t$-structure on $\scr{M}_A(S)$ is left-complete and compatible with filtered colimits.
    \end{thm}

    \begin{proof}
      We need to show that $\lim_n P_n(X)\simeq 0$. Since the forgetful functor commutes with limits and we have natural equivalences $U_A(P^n(X))\simeq P^n(U_A(X))$, we are reduced to the case of $A\simeq S$, which is treated in Corollary 3.8 of \cite{MR3781430}. On the other hand the subcategory $\scr{M}_A(S)_{\leq n}$ is closed under filtered colimits since the generators of $\scr{M}_A(S)_{\geq n+1}$ are compact. In particular the truncation functor $P^n$ commutes with filtered colimits.
    \end{proof}
    
    \begin{cor}
      Let $S$ be a Noetherian scheme of finite Krull dimension, and $A$ be a commutative algebra in $\scr{SH}_{\geq 0}$. Then the $\infty$-category $\scr{Mod}_A(S)$ is a presentably symmetric monoidal stable $\infty$-category and it is compactly generated. The homotopy $t$-structure in accessible, left-complete, multiplicative and compatible with filtered colimits.
    \end{cor}
    
    \subsubsection{} 
    \label{ssub:pi_0S}
    
    When $K$ is a perfect field of characteristic not $2$ we Morel constructs a map (see \cite[Section 6.1]{MR2061856})
    \[ \sigma:K^{MW}_\ast(K) \ra [\bb S, \bb G_m^{\wedge \ast}]_{\bb S} \]
    by defining it on the generators and checking that it passes to the quotient through the defining relations of $K^{MW}_\ast(K)$ (c.f. Definition 3.1 pg. 49 of \cite{MR2934577} for precise formulas).  

    In Theorem 6.2.1 of \cite{MR2061856} the author shows that $\sigma$ is actually an isomorphism, and in Corollary 6.4.1 he  shows that $\sigma$ extends uniquely to an isomorphism of homotopy modules $\sigma: \cal K^{MW}_\ast \overset{\simeq}{\ra} \underline \pi_0 \bb S$.

  \subsection{Homology Localizations} 
  \label{sub:homology_localizations}
    
    \begin{defin}
      \label{defin:E-local_stuff}
      Let $\scr{M}$ be a presentably symmetric monoidal stable $\infty$-category and let $E$ be an object of $\scr{M}$. We say that an arrow $f:X\ra Y$ in $\scr{M}$ is an $E$-homology equivalence (or, shortly, an $E$-equivalence) if the induced map $f\wedge \mathrm{id} :X\wedge E \ra Y\wedge E$ is an equivalence in $\scr{M}$. We say that anobject $C$ is \emph{$E$-acyclic} if $E\wedge C\simeq 0$ in $\scr{M}$. Finally we say that an object $X\in \scr{M}$ is \emph{$E$-local} if for every $E$-homology equivalence $X\ra Y$, the induced map on mapping spaces $\Map(Y,Z)\ra\Map(X,Z)$ is an equivalence of spaces. 
    \end{defin}
    \begin{rmk}
      \label{rmk:prop_of_acyc_obj}
      One sees immediately that the full sub-$\infty$-category $\rm{Ac}(E)$ spanned by $E$-acyclic objects is closed under arbitrary (small) colimits and retracts, and that $E$-acyclic objects have the 2-out-of-3 property in fiber sequences in $\scr{M}$. 
      More precisely if $K$ is any space and $p: K\ra \rm{Ac}(E)$ is a diagram whose composition with the inclusion $\rm{Ac}(E) \subseteq \scr{M}$ extends to a colimit diagram $\bar p : K^{\triangleright} \ra \scr{M}$, then there exists a unique lift $\bar p: K^{\triangleright} \ra \rm{Ac}(E)$ and such lift is again a colimit diagram. Furthermore since the inclusion functor $\rm{Ac}(E)\subseteq \scr{M}$ has a right adjoint (c.f. \ref{prop:infty_loc}), then whenever $\bar p: K^{\triangleright}\ra \rm{Ac}(E)$ is a colimit diagram, its composition $K^{\triangleright}\ra\scr{M}$ is again a colimit diagram; the converse to this statement is instead implied by the previous observation.
      Note that $\rm{Ac}(E)$ is also closed under smashing with an arbitrary object.
    \end{rmk}
    \begin{rmk}
      \label{rmk:prop_of_loc_obj}
      Similarly it is immediate to see that full sub-$\infty$-category $\rm{Loc}(E)$ spanned by $E$-local objects is closed under arbitrary (small) limits and retracts, and that $E$-local objects have the 2-out-of-3 property in fiber sequences. More precisely if $K$ is any space and $p: K\ra \rm{Loc}(E)$ is a diagram whose composition with the inclusion $\rm{Loc}(E) \subseteq \scr{M}$ extends to a limit diagram $\bar p : K^{\triangleleft} \ra \scr{M}$, then there exists a unique lift of $\bar p: K^{\triangleleft} \ra \rm{Loc}(E)$ and such lift is again a limit diagram.
      We finally note that an object $Z$ is $E$-local if and only if for every $E$-acyclic object $C$ the space $\Map(C,Z)$ is contractible, and this happens if and only if for every $E$-acyclic object $C$, the group $[C,Z]=0$.
    \end{rmk}

    \begin{rmk}
      \label{rmk:E-hom_eq_vs_E-loc_eq}
      Note that we could define, a priori, a more general notion of equivalence: we call a map $M\ra N$ an \emph{$E$-local equivalence} if, for every $E$-local object $X$ the natural map $\Map(N,X)\ra \Map(M,X)$ is an equivalence of spaces. Clearly all $E$-equivalences are $E$-local equivalences. The reverse holds if and only if the class $\rm{Ac}(E)$ of $E$-acyclic objects coincides with its double orthogonal 
      \[
      {}^{\perp} (\rm{Ac}(E)^{\perp}):=\{X\in \scr{M} \; s.t. \; \forall \; C\in \rm{Loc}(E), \Map(X,C)\simeq\ast\}.
      \]
      This follows immmediatly from the existence of a left adjoint to the inclusion $\rm{Loc}(E)\subseteq \scr{M}$
     \end{rmk}

    \begin{prop}
        \label{prop:infty_loc}
        Let $\scr{M}$ be a presentable symmetric monoidal stable $\infty$-category and let $E$ be an object of $\scr{M}$. Then 
        \begin{enumerate}[label=l.\arabic*]
          \item \label{l1} the inclusion $i_{\rm{Loc}}: \rm{Loc}(E) \subseteq \scr{M}$ has a left adjoint $\op L_E: \scr{M} \ra \rm{Loc}(E)$;
          \item \label{l2} for every $X$ in $\scr{M}$ there is an $E$-equivalence $X\ra X'$ with target $X'\in \rm{Loc}(E)$;
          \item \label{l3} the $\infty$-category $\rm{Loc}(E)$ is presentable;
          \item \label{l4} a map $f:X\ra Y$ in $\scr{M}$ is an $E$-equivalence if and only if $L_E(f)$ is an equivalence. 
        \end{enumerate}
        Moreover
        \begin{enumerate}[label=a.\arabic*]
          \item \label{a1} the inclusion $i_{\rm{Ac}}: \rm{Ac}(E)\subseteq \scr{M}$ has a right adjoint $\op A_E: \scr{M} \ra \rm{Ac}(E)$;
          \item \label{a2} for every $X$ in $\scr{M}$ there is a co-local equivalence $X''\ra X$ with source $X'' \in  \rm{Ac}(E)$;
          \item \label{a3} the $\infty$-category $\rm{Ac}(E)$ is presentable;
        \end{enumerate}
    \end{prop}
    
    \begin{proof}
        Let us concentrate on $E$-local objects first. The class of $E$-equivalences in $\scr{M}$ is a strongly saturated class of morphisms according to Definition 5.5.4.5 of \cite{HTT}. Moreover the collection of $E$-equivalences is a strongly saturated class of small generation: this can be seen combining Proposition 5.5.4.16 and of Remark 5.5.4.7 in \cite{HTT}), given that the functor $-\wedge E$ is presentable.
        
        Let us denote by $\scr T_E$ a choice of a small set generating $E$-equivalences as a strongly saturated class. As a consequence we can apply Proposition 5.5.4.15 of \cite{HTT} and conclude immediately for \eqref{l1},\eqref{l2} and \eqref{l3}. We also deduce part of \eqref{l4}, i.e. that $f:X\ra Y$ is an $E$-equivalence if and only if $f$ is a $\scr T_E$-local equivalence, if and only if $L_E(f)$ is an equivalence. 

        Let us now tun to $E$-acyclic objects. We immediately note that applying Proposition 5.2.7.8 to $\scr{M}^{op}$ we get that \eqref{a2} is equivalent to \eqref{a1}. In order to prove \eqref{a1} we take, for every $A$-module $X$ an $E$-equivalence to an $E$-local object $\lambda: X\ra X'$ and define $X'' \ra X$ as the fiber of $\lambda$. Then clearly $X''$ is $E$-acyclic and for every $E$-acyclic object $C$ composing with $\lambda$ induces an equivalence $\Map(C,X'') \ra \Map(C, X)$, given that $\Map(C,X')$ is contractible by assumption. In order to prove \eqref{a3} we use Proposition 1.4.4.13 of \cite{HA}        
    \end{proof}

    \begin{defin}
    \label{defin:E-localization}
      We call \emph{$E$-localization functor} a choice of composition $(-)_E:= L_E\circ i_{\rm{Loc}}$, while a unit transformation $\lambda_E: \op{id} \ra (-)_E$ is called \emph{$E$-localization map}. Similarly we call \emph{$E$-acyclicization functor} a choice of composion ${}_E(-): i_{\rm{Ac}}\circ \op A_E$ and a co-unit transformation ${}_E(-)\ra \op {id}$ is called \emph{$E$-acyclicization map}.
    \end{defin}
    
    \begin{prop}(c.f. 2.2.1.9 of \cite{HA})
    \label{prop:prop_of_loc}
      The functor $L_E$ is symmetric monoidal and the natural transformation $\op{id}\ra \iota_E \circ L_E$ is monoidal.
    \end{prop}

  \subsection{Bousfield classes} 
    \label{sub:bousfield_classes}

    \subsubsection{}
    Let $\scr{M}$ be a presentably symmetric monoidal stable $\infty$-category. We introduce an equivalence relation on the class of equivalence classes of objects of $\scr{M}$ following \cite{MR543337}. We set $E\sim_B F$ if, for every $X\in \scr{M}$, we have that $E\wedge X=0$ if and only if $F\wedge X=0$. By Theorem \ref{prop:infty_loc} localization functors at $E$ and $F$ exist and they are (canonically) isomorphic exactly when $E\sim_BF$. We denote by $\cal A(\scr{M})$ the class of Bousfield classes in $\scr{M}$ and by $\bc E$ the element in $\cal A(\scr{M})$ represented by an object $E$. On $\cal A(\scr{M})$ we introduce a partial ordering by setting $\bc E \leq \bc F$ if every $F$-acyclic object is $E$-acyclic. 
      
    \subsubsection{}
    Given a possibly infinite collection of Bousfield classes $\bc {E_i}_{i \in I}$ we have a \emph{join} operation which is defined as $\vee_{i \in I}\bc {E_i}:=\bc{\vee_{i \in I}E_i}$. We note that the join is always the minimal upper bound of its summands. 

    Similarly, given a finite collection of Bousfield classes $\bc {E_i}_{i \in I}$ we have a \emph{meet} operation which is defined as $\wedge_{i \in I}\bc {E_i}:=\bc{{\wedge}_{i \in I}E_i}$. We note that the meet operation is a lower bound for its factors, but in general does not need to be the maximal lower bound.

    \subsubsection{}
    Following \cite{MR543337} we denote by $\cal{DL}(\scr{M})$ the subclass of $\cal A(\scr{M})$ of those Bousfield classes $\bc E$ such that $\bc E \wedge \bc E =\bc E$. The operations of meet and join restrict to $\cal{DL}(\scr{M})$ and it is elementary to check that $\cal{DL}(\scr{M})$ satisfies the axioms of a distributive lattice. The partial ordering $\leq$ on $\cal A(\scr{M})$ restricts to a partial ordering on $\cal{DL}(\scr{M})$. We wish to observe that for given $\bc E,\bc F \in \cal{DL}(\scr{M})$ their meet $\bc{E\wedge F}$ is actually their greatest upper bound. Most of the objects we will consider later actually belong to this subclass: for instance, every homotopy algebra $E$ belongs to $\cal {DL}(\scr{M})$, since $E$ is a retract of $E\wedge E$.

    \subsubsection{}
    We say that a Bousfield class $\bc{E} \in \cal A(\scr{M})$ has a complement if there is another Bousfield class $\bc{F}$ such that $\bc{E}\wedge\bc{F}=\bc{0}$ and $\bc{E}\vee \bc{F}=\bc{\bf 1}$. If $\bc E$ has a complement then such complement is unique, and we denote it by $\bc{E}^c$. Moreover when $\bc E$ has a complement, $\bc{E} \in \cal{DL}(\scr{M})$. We denote by $\cal{BA}(\scr{M})$ the sub-lattice of $\cal{DL}(\scr{M})$ of those Bousfield classes admitting a complement. Lemma 2.7 of \cite{MR543337} shows that the inclusion $\cal{BA}(\scr{M})\subseteq \cal{DL}(\scr{M})$ is in general strict. Assume that both $\bc E, \bc F$ have complements, then the following equalities are satisfied:
    \begin{equation}
      \label{eqn:bc_complements_formulas}
      \bc E^{cc}=\bc E, \;\;\; \bc{E\vee F}^c=\bc{E}^c\wedge \bc{F}^c, \;\;\; \bc{E\wedge F}^c=\bc{E}^c\vee \bc{F}^c.
    \end{equation}

\section{Moore objects} 
  \label{sec:moore_spectra}

  In the topological stable category $\scr{SH}$, one can associate to every abelian group $A$ a Moore spectrum $\bb S A$. Up to isomorphism $\bb S A$ is characterized by the following properties: $\bb S A\in \scr{SH}_{\geq 0}$, $\pi_0(\bb S A)\simeq A$, and $\op H \bb Z\wedge \bb S A\simeq \op HA$. Moore spectra are fundamental in the study of Bousfield classes for two reasons. One is that $\bc{\bb S A}$ depends only on torsion and divisibility properties of $A$ (c.f. Proposition $2.3$ of \cite{MR551009}). The other is that for a spectrum $X \in \scr{SH}_{\geq k}$, and a homotopy commutative ring spectrum $E\in \scr{SH}_{\geq 0}$ we have that $X_E\simeq X_{\bb S \pi_0E}$. 
  Given a set $J$ and a finite set $I$, and collections of maps $\{f_i:L_i\ra \bf 1\}_{i \in I}$, $\{g_j:L_j \ra \bf 1\}_{j \in J}$, we introduce a weak version of Moore object 
  \[\bf 1\cal R =C(f_1)\wedge\dots\wedge C(f_r)\wedge \bf 1[\cal J^{-1}]\big.\]
   $\bf 1 \cal R \in \scr{M}_{\geq0}$ and $\tau_0(\bf 1\cal R)\simeq \cal R$, but in general $\bf 1\cal R\wedge \tau_0(\bf 1) \not \in \scr{M}^{\heartsuit}$, and furthermore $\bf 1\cal R$ depends on the choice of the $f_i$'s and $g_j$'s.
   We start the section with reviewing some technical tools for dealing with towers in \ref{sub:towers}. We then treat the $f_i$' and $g_j$' separately in \ref{sub:coning_off_homotopy_elements} and \ref{sub:inverting_homotopy_elements}; in \ref{sub:eta_completion} we make an example on $\eta$-completions in $\scr{SH}(S)$.

  In this section $\scr{M}$ is a presentably symmetric monoidal stable $\infty$-category, endowed with a left-complete accessible multiplicative $t$-structure.

  \subsection{Towers}
    \label{sub:towers}
    Let $K^n$ be the sub simplicial set
      \[\Delta^{\{0,1\}}\coprod_{\Delta ^0} \Delta^{\{1,2\}} \cdots \Delta^{\{n-2,n-1\}} \coprod_{\Delta^0}\Delta^{\{n-1,n\}} \subseteq \Delta^n,\] where the pushouts are taken with respect to the maps
      \[ \Delta^{\{k-1,k\}} \overset{k}{\leftarrow} \Delta^{0} \overset{k}{\rightarrow} \Delta^{\{k,k+1\}}. \] In other words $K^n$ is the sub simplicial set of $\Delta^n$ generated by its non-degenerate $1$-simplices. We also define $K=\cup_n K^n$. A diagram in $\scr{M}$ indexed by ${K^n}^{op}$ (resp. $K^{op}$) is thus a collection of $n$ composable arrows of $\scr{M}$:
      \[\xymatrix{X_n\ar[r] & X_{n-1} \ar[r] & \cdots \ar[r] & X_0}\]
      (resp. a countable collection of composable arrows
      \[\xymatrix{\cdots \ar[r] & X_n\ar[r] & X_{n-1} \ar[r] & \cdots \ar[r] & X_0}).\]
       A composition of $n$ composable arrows $p:{K^n}^{op}\ra \scr{M}$ is an extension of $p$ to ${\Delta^n}^{op}$, and is essentially unique in the following sense.

    \begin{lemma}
      \label{lem:Kn_Deltan}
      \begin{enumerate}
        \item \label{c1} The inclusion $i_n : K^n\subseteq \Delta^n$ is an inner anodyne map;
        \item \label{c2} If $\scr{C}$ is an $\infty$-category, the natural restriction map 
      \[i_n^\ast: \op{Fun}({\Delta^n}^{\op {op}}, \scr C) \ra \op{Fun}({K^n}^{\op {op}}, \scr C)\] is a trivial Kan fibration between $\infty$-categories, and in particular a categorical equivalence.
      \end{enumerate}
      \end{lemma} 
    
    \begin{proof}
      For \eqref{c1} observe that when $n=0,1$ there is nothing to do, while for $n=2$ we have just rewritten Corollary 2.3.2.2 of \cite{HTT}. Assume now that the inclusion $i_n : K^n\subseteq \Delta^n$ is inner anodyne. The inclusion $i_{n+1}$ is the composition
      \[ K^{n+1}= K^{n} \coprod_{\Delta^{\{n\}}}\Delta^{\{n,n+1\}} \subseteq \Delta^n \coprod_{\Delta^{\{n\}}}\Delta^{\{n,n+1\}}  \subseteq \Delta^{n+1},\]
      and we claim that both the above inclusions are inner anodyne. Indeed the central inclusion is $i_n \coprod_{\Delta^{\{n\}}}\Delta^{\{n,n+1\}}$, and this is inner anodyne because the collection of inner anodyne maps is closed under push-outs (being a weakly saturated collection of maps). The rightmost inclusion is identified with the natural inclusion 
      \[ \Delta^{n-1} \ast\Delta^0 \coprod _{\emptyset \ast\Delta^0} \emptyset \ast \Delta^1  \subseteq \Delta^{n-1} \ast \Delta^1,\]
      which is inner anodyne after Lemma 2.1.2.3 of \cite{HTT}.

      For \eqref{c2} one combines Propositions 2.3.2.1 and 1.2.7.3 of \cite{HTT} with \eqref{c1} and with the fact that a map is inner anodyne if and only if the opposite is.
    \end{proof} 
    \begin{lemma}
      \label{lemma:K_Delta}
        let $\bf N^{op}$ be the category associated to the poset $\bb N$ endowed with decreasing order, and $\bf N^{op}_{\leq n}$ be the full subcategory with objects $\{n,\dots,0\}$. Then the obvious isomorphism of simplicial sets $\op N(\bf N^{op}) \simeq \cup_n \op N(\bf N^{op}_{\leq n})$, the inclusions $i_n: K^{n}\subseteq \Delta^n$, and the natural identifications $\op N (\bf N^{op}_{\leq n})\simeq {\Delta^n}^{op}$ altogether induce a trivial Kan fibration of $\infty$-categories
      \[ p: \op{Fun}(\op N (\bf N^{op}),\scr C)\overset{\simeq}{\ra}\lim_n\op{Fun}({K^n}^{\op {op}}, \scr C), \]
      which has an essentially unique section $s$.
      \end{lemma}
      \begin{proof}
      Recall that inclusions are cofibrations in the Joyal model structure on simplicial sets (cf. Theorem 2.2.5.1 of \cite{HTT}), and filtered diagrams of cofibrations are cofibrant. Since Joyal's monoidal model structure is a model for $\rm{Cat}_\infty$, we conclude that $\op N(\bf N)\simeq \colim_n \op N (\bf N_{\leq n})$ as $\infty$-categories, and thus that the natural map
      \[\op{Fun}(\op N (\bf N), \scr C) \simeq \lim_n\op{Fun}(\op N (\bf N_{\leq n}), \scr C)\]
      is a categorical equivalence and a categorical fibration of simplicial sets. Using again  Theorem 2.2.5.1 of \cite{HTT}, which asserts that categorical fibrations and categorical equivalences are respectively the fibrations and weak equivalences of Joyal's structure, we conclude that the above map is a trivial Kan fibration. On the other hand, \ref{lem:Kn_Deltan} yields a chain
      \begin{equation}
        \lim_n\op{Fun}(\op N (\bf N^{op}_{\leq n}), \scr C) = \lim_n\op{Fun}({\Delta^n}^{\op {op}}, \scr C) \overset{\lim i_n^\ast}{\ra} \lim_n\op{Fun}({K^n}^{\op {op}}, \scr C)
      \end{equation}
      of trivial Kan fibrations, since these are stable under inverse limits. Thus $p$ is a trivial Kan fibration, and as such it has an essentially unique section.
    \end{proof}
    
    \subsubsection{}
    \label{subsub:const_of_towers}
    As a consequence, if we have a countable family of composable arrows $f_n: X_n \ra X_{n-1}$ we can essentially uniquely extend this collection to a diagram $\bf f:N(\bf N^{op})\ra\scr{M}$. If $\bf c X_0$ is the constant diagram on $X_0$, we can thus upgrade the datum of the $f_n$'s to a map of towers $f^{\bullet}:\bf f\ra \bf cX_0$ which is a $1$-simplex of $\rm{Fun}(K^{op},\scr{M})$. Finally taking the cofiber of $f^\bullet$ yields a commutative ladder of fiber sequences:
    \[\label{}
      \xymatrix@C=5em{
      X_n \ar[r]^{f^{ n} }\ar[d]_{f_n} &  X_0 \ar[r] \ar@{=}[d] &  C(f^{n}) \ar[d]_{p_n} \\
     X_{n-1} \ar[r]^{f^{ n-1}} &  X_0 \ar[r] &  C(f^{n-1}),\\ 
     }
    \]
    where $f^{ n}$ denotes a composition of the $f_i$ for $i\leq n$.

 \subsection{Quotients} 
    \label{sub:coning_off_homotopy_elements}

    \begin{defin}
    \label{defin:kill_one_homot_elt}
      Let $r\geq 1$ be an integer and for $i\in \{1,\dots,r\}$ let $f_i: L_i \rightarrow \bf 1$ be a map in $\scr{M}$. We denote by $C(f_i)$ the cofiber of $f_i$ and by $M(\underline f)$ the iterated cofiber $C(f_1)\wedge \cdots \wedge C(f_r)$. We call $M(\underline f)$ the Moore object associated with the collection $f_1,\dots,f_r$.     
    \end{defin}

    \begin{rmk}
      \label{rmk:diff_wrt_top_moore_spt}
      If $\M=\scr{SH}$, $r=1$, and $n:\bb S \ra \bb S$ represents the multiplication by $n \in \bb Z$, $C(n)$ is then the usual modulo $n$ Moore spectrum. Note that when $r\geq 2$ the spectrum $M(n_1,n_2)=C(n_1)\otimes C(n_2)$ is not always a Moore spectrum in the classical sense: take $n_1=n_2=2$ for instance. However for every pair of integers $n_1,n_2$, $\bc{M(n_1,n_2)}=\bc{C(\rm{g.c.d.}(n_1,n_2))}$: this follows from Proposition 2.13 of \cite{MR543337}. This observation points out that our definition of Moore object $M$ really depends on choices, and not only on its $\tau_0(M)$. 
    \end{rmk}

    \subsubsection{}
      \label{subsub:right_and_left_mult}
      We introduce some further notation. If $f:L\ra \bf 1$ is a map in $\scr{M}$, we will denote by $l_f(X)$ the left multiplication by $f$ on $X\in \scr{M}$, i.e. the composition 
      \[L\wedge X \overset{f\wedge X}{\rra} \bf 1\wedge X \overset{\simeq}{\ra}X.\]
       Similarly $r_f(X)$ will denote the right multiplication by $f$, given by the composition 
       \[X\wedge L \overset{X\wedge f}{\rra} X\wedge \bf 1 \overset{\simeq}{\ra}X.\]
       Note that left multiplications commutes with any map in $\phi: X \ra Y$ in $\scr{M}$, in the sense that the square
       \[\xymatrix{
       L\wedge X \ar[r]^{l_f(X)} \ar[d]_{L\wedge\phi} & X \ar[d]^{\phi}\\ 
       L\wedge Y \ar[r]^{l_f(Y)}  & Y \\
       }\]
       commutes; right multiplication behaves similarly.

    \begin{defin}
      \label{defin:sdf_object}
      An object $L$ of $\scr{M}$ is called \emph{flat} if the functor $L\wedge-$ respects $\scr{M}_{\geq 0}$ and $\scr{M}_{\leq 0}$. $L$ is called \emph{sdf} (resp. \emph{tif}) if it is strongly dualizable (resp. $\wedge$-invertible) and flat.
    \end{defin}

    \subsubsection{}
      \label{subsub:ex_of_sdf_and_tif_objects}
      In the motivic setting, the main example of a tif object in $\scr{SH}(S)$ is of course $\Sigma^{-\rm{rk}(V)}\rm{Th}(V)\in \scr{SH}(S)$ where $V$ is a virtual vector bundle on the base scheme $S$, and where $\rm{rk}(V)$ is interpreted as a locally constant function on $S$. For instance $\bb G_m$ is a tif object in $\scr{SH}(S)$. 
    \begin{lemma}
      \label{lemma:left_mult_on_heart}
      Let $f:L \ra A$ be a map with flat source and let $X$ be any object of $\scr{M}$. A canonical zig-zag of natural maps induces an equivalence $\tau_i(L\wedge X)\simeq L\wedge \tau_i(X)$. Under this equivalence the map $\tau_i(l_f(X)):\tau_i(L\wedge X) \ra \tau_i(X)$ ans the map $l_f(\tau_i(X)):L\wedge \tau_i(X)\ra\tau_i(X)$ are naturally identified. Similar statements fold for right multiplication.
    \end{lemma}
    \begin{proof}
      Diagram chase using the fiber sequences 
      \[P_{n}(-) \ra (-) \ra P^{n-1}(-)\; \textrm{ and } \; P_{n+1}(-) \ra P_{n}(-) \ra \Sigma^n\tau_n(-)\]
      for $X$ and $L\wedge X$.
    \end{proof}

    \begin{lemma}
    \label{lemma:Moore_is_conn_and_with_corr_pi0}
      For every $i \in \{1,\dots,r\}$ let $f_i:L_i\ra \bf 1$ be a map in $\scr M$, where $L_i$ is a flat object, and let $(f_1,\dots, f_r)$ be the sub-$\tau_0(\bf 1)$-module of $\tau_0(\bf 1)$ obtained as the image of the map
      \[\Sigma_i l_{f_i}(\tau_0(\bf1)): \bigoplus_i L_i\wedge \tau_0(\bf 1) \ra \tau_0(\bf 1)\]
      in $\scr{M}^{\heartsuit}$.
      Then $M(\underline f)\in \scr{M}_{\geq 0}$,  
      \[ \tau_0 M(\underline f)\simeq  \otimes^{\heartsuit}_i\tau_0(C(f_i))\simeq \tau_0(A)/(f_1,\dots,f_r),\]
      and the canonical map $\bf 1 \ra M(\underline f)$ induces on $\tau_0$ the quotient map $ \tau_0 (\bf 1) \ra  \tau_0 (\bf 1)/(f_1,\dots,f_r)$. 
    \end{lemma}
    \begin{proof}
      An easy diagram chase allows to reduce to the following claim: If $\bf A$ is an abelian category, $M,N,P$ are objects of $A$, and we have a diagram $f:N\rightarrow M\leftarrow P:g$, then $(M/N)/\overline {\rm{Im}(g)}\simeq M/\rm{Im}(f+g)$, where $f+g: N\oplus P \ra M$ is the natural map induced by $f$ and $g$. On its turn this claim is easily proved using that colimits commute with each other.
    \end{proof}

    \subsubsection{}
      \label{const:modx_1_x_n_Moore_sp}
      Let $f: L\ra \bf 1$ be a map in $\scr{M}$ and consider the collection of composable maps $r_f(L^{\wedge n-1}):L^{\wedge n} \ra L^{\wedge n-1}$. By applying the argument of \ref{subsub:const_of_towers} we obtain a commutative ladder of fiber sequences
      \begin{equation}
          \label{eqn:mod_x^n_mod_x^n-1_fib_seqs}
          \xymatrix@C=5em{
            L^{\wedge n} \ar[r]^{f^{n}\cdot }\ar[d]_{r_f(L^{\wedge(n-1)})} &  \bf 1 \ar[r] \ar@{=}[d] &  C(f^n) \ar[d]_{p_n} \\
            L ^{\wedge (n-1)} \ar[r]^{f^{n-1}\cdot} &  \bf 1 \ar[r] &  C(f^{n-1}).\\ 
          }
        \end{equation}
        Note that we have equivalences
      \begin{equation}
          \label{eqn:mod_x^n/mod_x^n-1_is_mod_x}
        \begin{split}
          \fib \Big (C(f^n) \overset{p_n}{\ra} C(f^{n-1}) \Big) & \simeq \cofib \Big(L ^{\wedge n} \overset{r_f(L^{\wedge(n-1)})}{\rrra} L^{\wedge(n-1)} \Big )\\
          & \simeq L^{\wedge(n-1)}\wedge C(f)
        \end{split}
      \end{equation}
      yielding a fiber sequence
      \begin{equation}
          \label{eqn:mod_x_mod_x^2_mod_x_cof_seq}
          L^{\wedge(n-1)}\wedge C(f)\ra C(f^n) \overset{p_n}{\ra} C(f^{n-1}).
      \end{equation}

    \begin{defin}
      \label{defin:f_adic_comp}
        Let $f:L\ra \bf 1$ be a map in $\scr{M}$. In view of \ref{const:modx_1_x_n_Moore_sp} we define the \emph{$f$-adic completion} of $X\in \scr{M}$ as the object
      \begin{equation}
      \label{eqn:def_adic_comp}
        X^{\wedge}_f:=\lim_{N(\bf N^{op})} (X \wedge C(f^\bullet)).
      \end{equation}
      We call \emph{$f$-adic completion map} the (essentially unique) map $\chi_f(X): X \ra X^\wedge_f$ induced by $X\simeq \lim (X\wedge \bf 1 )\ra \lim (X\wedge C(f^\bullet))$. An object $X$ is \emph{$f$-complete} if $\chi_f(X)$ is an equivalence.
    \end{defin}
    \begin{rmk}
      Since the operations of taking inverse limits and of smashing with an object $X$ commute with finite limits, the fiber sequences of towers introduced in \eqref{const:modx_1_x_n_Moore_sp} yields the fiber sequence
      \begin{equation}
      \label{eqn:def_seq_of_adic_comp}
        \lim (X\wedge L^{\wedge \bullet}) \ra X \overset{\chi_f(X)}{\ra} X^{\wedge}_f.     
       \end{equation} 
      In particular when $X$ is $k$-connected and $L\in \scr{M}_A(S)_{\geq 1}$ $X$ is $f$-complete.
    \end{rmk}

    \begin{defin}
      If $r\geq 1$ then for every $X\in \scr{M}$ we call $\underline f$-adic completion of $X$ the object 
      \[X^{\wedge}_{\underline f}:=X^{\wedge}_{f_1}{}^{\wedge}_{f_2}{\cdots}^{\wedge}_{f_r} \in \scr{M}.\] The $\underline f$-adic completion map $\chi_{\underline f}(X)$ is defined as a composition of the natural maps 
      \[\chi_{f_r}(X^{\wedge}_{f_1,\dots,f_{r-1}}) \circ \cdots \circ \chi_{f_2}(X^\wedge_{f_1})\circ\chi_{f_1}(X).\]
    \end{defin}
        
    \begin{lemma}
    \label{lemma:dualmaps}
      Let $L$ be a $\wedge$-invertible object of $\scr{M}$ and $f: L\ra \bf 1$ be a map in $\scr{M}$. Then the cofiber $C(f)$ is strongly dualizable in $\scr{M}$ and 
      \[D(C(f))\simeq \rm{fib}(D(f))\simeq \Sigma^{-1} D(L)\wedge C(f),\]
      where $D(-)=\Hom(-,\bf 1)$ and $\Hom(-,-)$ denotes a right adjoint of $-\wedge-$.
    \end{lemma}
    \begin{proof}
      We have that 
      \[D(C(f))\simeq \rm{fib}(D(f))\simeq \rm{fib}(D(L)\wedge f)\simeq D(L)\wedge\rm{fib}(f)\simeq D(L)\wedge \Sigma^{-1}C(f).\] The fact that $\rm{fib}(D(f))\simeq \rm{fib}(D(L)\wedge f)$ follows directly form the definition of dual map via the duality adjunction.
    \end{proof}
    \begin{prop}
      \label{prop:localization_at_mod_x_moore_spectrum}
      Let $L$ be a $\wedge$-invertible object of $f: L \ra \bf 1$ be a map in $\scr{M}$. Then for every object $X$ the natural map $\chi_f(X): X\ra X^{\wedge}_f$ is a $C(f)$-localization of $X$ in $\scr{M}$.
    \end{prop}
    \begin{proof}
      We need to check that $\chi_f(X)$ is a $C(f)$-equivalence and that $X^\wedge_f$ is $C(f)$-local. Let $\phi: M\ra N$ be a map in $\scr{M}$. For every $F\in \scr{M}$ we get an induced map
      \begin{equation}
        \Map\big (C(f)\wedge N,\big ) \ra \Map\big (C(f)\wedge M,F\big ),
      \end{equation}
      and since $C(f)$ has a strong dual $D(C(f))$, we get a natural map 
      \begin{equation}
      \label{eqn:mapdual}
        \Map\big (N, D(C(f))\wedge F\big ) \ra \Map\big (M,D(C(f))\wedge F\big ).
      \end{equation}
      $\phi$ is a $C(f)$-equivalence if and only if the map \eqref{eqn:mapdual} is an equivalence for every $F$. In particular, for every $F$ in $\scr{M}$, $D(C(f))\wedge F$ is $C(f)$-local. Furthermore \ref{lemma:dualmaps} implies that $D(C(f))\simeq \Sigma ^{-1}D(L)\wedge C(f)$. Since local objects are stable under wedging with invertible objects, we conclude that for every $F\in \scr{M}$ the object $C(f)\wedge F$ is $C(f)$-local. 

      Now we show that $X^{\wedge}_f$ is $C(f)$-local: since local objects are closed under inverse limits we reduce to showing that $X\wedge C(f^n)$ is $C(f)$-local. This easily follows by induction: the base case being that $X\wedge C(f)$ is $C(f)$-local, which was observed above. Assume we know that $X\wedge C(f^{n-1})$ is $C(f)$-local. Using the fiber sequence 
      \begin{equation}
        \label{eqn:C_x^n/C_c^n-1_is_Cx}
        L^{\wedge n-1}\wedge C(f) \wedge X \ra C(f^{n}) \wedge X \ra C(f^{n-1})\wedge X
      \end{equation}
      deduced from \eqref{eqn:mod_x_mod_x^2_mod_x_cof_seq} and the 2-out-of-3 property of $C(f)$-local objects in fiber sequences we conclude.

      In order to show that the canonical map $X \ra X^{\wedge}_f$ is a $C(f)$-local equivalence it suffices to show that $C(f)\wedge Y \simeq 0$ in $\scr{M}$, where $Y:=\fib(X \ra X^{\wedge}_f)$. For this note that
      \begin{equation}
        \label{eqn:tower_desc_fiber_of_x_completion}
        Y\simeq \lim (X\wedge L^{\wedge \bullet} )= \lim \Big ( \cdots \ra   X\wedge L^{\wedge n} \overset{r_{f}(L^{\wedge n-1})}{\rra} X\wedge L^{\wedge(n-1)}\ra \cdots \Big )
      \end{equation}
      and that $C(f)\wedge Y\simeq \cofib(r_f(Y): Y\wedge L \ra Y)$. However it is easily checked that the multiplication by $f$ on $Y$ is induced by the multiplication by $f$ on each component of the tower $X\wedge L^{\wedge \bullet}$. Since the inverse limit of a diagram only depend on the pro-object associated with the diagram (see Appendix \ref{sec:pro_spectra}), and since the multiplication by $f$ is clearly an equivalence of pro-objects, we conclude that $r_f(Y)$ is actually an equivalence in $\scr{M}$ and so its cofiber is zero.
    \end{proof}

    \begin{prop}
      \label{prop:localization_at_mod_x_1-x_n_moore_spectrum}
      Let $r\geq 2$ be an integer, for every $i=1, \dots, r$ let $f_i: L_i \ra A$ be a map in $\scr{M}$, and assume that for all $i$, $L_i$ is $\wedge$-invertible. Then for every $X\in \scr{M}$ the natural map $\chi_{\underline f}(X): X\ra {X^{\wedge}_{f_1}} \cdots^{\wedge}_{f_r}$ is a $M(\underline f)$-localization of $X$ in $\scr{M}$.
    \end{prop}
    \begin{proof}
      Let us first set the notation $M(\underline f)=C(f_1)\wedge \cdots \wedge C(f_r)$. The object $M(\underline f)$ has a strong dual 
      \begin{equation}
        D(M(\underline f)\big) \simeq \Sigma^{-r}\wedge_{i=1}^r D(L_i)\wedge M(\underline f).
      \end{equation}
      Hence by running the same argument as in the proof of  \ref{prop:localization_at_mod_x_moore_spectrum} we deduce that $F\wedge M(\underline f)$ is $M(\underline f)$-local for every $F\in \scr{M}$. In order to show that ${X^{\wedge}_{f_1}} \cdots^{\wedge}_{f_n}$ is $M(\underline f)$-local, thanks to the identification
      \[{X^{\wedge}_{f_1}} \cdots^{\wedge}_{f_r} \simeq \lim_{\op N(\bf N^{op})}(X\wedge C(f_1^\bullet)\wedge \cdots \wedge C(f_r^\bullet)),\]
      we only need to prove that each of the objects $X\wedge C(f_1^{n})\wedge \cdots \wedge C(f_r^n)$ is $M(\underline f)$-local. This can be done by induction using iteratively the fiber sequence \eqref{eqn:C_x^n/C_c^n-1_is_Cx} and the fact that $M(\underline f)$-local objects satisfy the $2$-out-of-$3$ property in fiber sequences.

     The natural map 
     \begin{equation}
       X\overset{\chi_{f_1}}{\rra} X^{\wedge}_{f_1} \overset{\chi_{f_2}}{\rra} X^{\wedge}_{f_1}{}^{\wedge}_{f_2}\overset{\chi_{f_3}}{\rra} \cdots \overset{\chi_{f_r}}{\rra} {X^{\wedge}_{f_1}} \cdots^{\wedge}_{f_r} 
     \end{equation}
     is a composition of $M(\underline f)$-equivalences since $\bc{M(\underline f)}\leq \bc{C(f_i)}$ and since $\chi_{f_i}$ is a $C(f_i)$-equivalence for every $i=1,\dots, r$ by \ref{prop:localization_at_mod_x_moore_spectrum}.
    \end{proof}


  \subsection{A remark on $\eta$-completions}
  \label{sub:eta_completion}
      
    Let $\spec K$ be a perfect field of characteristic $\not = 2$, and consider the case where $\scr{M}=\scr{SH}(K)$. In this setting let $\eta \in \pi_0(\bb S)_{-1}(K)$ be the algebraic Hopf map. We have proved above in \ref{prop:localization_at_mod_x_moore_spectrum} that for every spectrum $X$ the $\eta$-completion map $\chi_\eta(X):X\ra X^\wedge_\eta$ is the $M(\eta)$-localization map in $\scr{SH}(K)$. We want to bring the discussion on $\eta$-completions a bit further.
    \begin{lemma}
      \label{lem:Xonehalf_is_eta_complete}
      Assume that the base field $K$ is not formally real. Then for every spectrum $X \in \scr{SH}(K)$ the spectrum $X[\frac{1}{2}]$ is $\eta$-complete.
    \end{lemma}
    \begin{proof}
      It follows from \cite[Ch. 2, Theorem 7.9]{MR770063} that there exists an integer $n$ such that $2^n$ acts as $0$ on the Witt ring of $K$. In particular we deduce that in $\rm{GW}(K)$ the relation $2^n=h\omega$ holds, where $h$ is the rank $2$ hyperbolic space and $\omega$ is some element of $\rm{GW}(K)$. It follows that $2^n\eta=h\omega\eta=0$ in $\op K^{MW}_\ast(K)$. It follows that on $X[\frac{1}{2}]$ the multiplication by $\eta$ is the zero map, which in view of \ref{const:modx_1_x_n_Moore_sp} is enough to conclude. 
    \end{proof}
    \begin{lemma}
      \label{lemma:compact_is_eta_complete}
      Assume that the base field has finite \'etale $2$-cohomological dimension. Then every strongly dualizable object of $\scr{SH}(K)$ is $\eta$-complete.
    \end{lemma}
    \begin{proof}
      If $C$ be a dualizable object of $\scr{SH}(K)$, the operation of smashing with $C$ commutes with inverse limits. In particular $C^\wedge_\eta \simeq C\wedge \bb S^\wedge_\eta$ and thus we reduce to show that $\bb S$ is $\eta$-complete. In view of \ref{prop:fracture_square_for_mod_eta_cohomologies} we just need to show that the spectra $\bb S[\frac{1}{2}]$, $(\bb S^\wedge _2)[\frac{1}{2}]$ and $\bb S^\wedge_2$ are $\eta$-complete. For the two former spectra the previous claim follows from \ref{lem:Xonehalf_is_eta_complete}, while for $\bb S^\wedge_2$  the claim follows from the combination of Proposition 4 and Lemma 21 of \cite{MR2811716}.
    \end{proof}

  \subsection{Inversions} 
   \label{sub:inverting_homotopy_elements}
    
    \subsubsection{}
    \label{ssub:rec_on_triv} Let $J$ be a set and $\cal B=\{B_j\}_{j \in J}$ be a collection of objects of $\scr{M}$. Since $\scr{M}$ is presentable we can choose a set $G$ of $\kappa$-compact objects that generate $\scr{M}$ under $\kappa$-filtered colimits. Note that we can assume that $\kappa$ has been chosen so that $\bf 1$ is a $\kappa$-compact object of $\scr{M}$. We define now $\scr B$ as the smallest full sub-$\infty$-category of $\scr{M}$ which contains the objects of the form
      \[V \wedge B_j\: \textrm{ where } \; V \in G, j\in J,\] and which is closed under small colimits and extensions. Using Proposition 1.4.1.1 of \cite{HA} we deduce that $\scr B$ is a presentable $\infty$-category. We set $\scr B'$ to be the full sub-$\infty$-category of $\scr{M}$ spanned by those objects $X$ such that $\Map_A(C,X)\simeq \ast$ for every $C\in \scr B$. The pair $\{\scr B, \scr B'\}$ forms then an accessible $t$-structure on $\scr{M}$ (Proposition 1.4.4.11 of \cite{HA}). We denote by $\tau^{\cal B}:{}^{\cal B}(-) \ra \op{id}$ the associated co-localization map and by $\lambda^{\cal B}:\op{id}\ra (-)^{\cal B}$ the associated localization map. The notation we used here is inspired by \cite{MR543337}.

    \begin{lemma}
      The functor $(-)^{\cal B}$ is symmetric monoidal and that the natural transformation $\lambda^{\cal B}$ is monoidal.
    \end{lemma}

    \begin{proof}
      This follows essentially from \cite[Proposition 2.2.1.9]{HA}. Indeed according to \cite[Example 2.2.1.7]{HA} we only need to check that $\scr B$ is closed under smashing with arbitrary objects. This follows immediately from the fact that $-\wedge-$ commutes with colimits in both variables, since $\scr{M}$ is presentably symmetric monoidal, and the definition of $\scr B$.
    \end{proof}

    \begin{prop}
      \label{prop:propert_of_trivializ}
      Let $\cal B:=\{B_j\}_{j \in J}$ be a set of strongly dualizable objects of $\scr{M}$. We have that:
      \begin{enumerate}[label=t.\arabic*]
        \item \label{t1} for every pair of objects $X$ and $Y$ of $\scr{M}$, ${}^{\cal B}X\wedge Y^{\cal B}\simeq 0$;
        \item \label{t6} for every $X\in \scr{M}$ we have $X^{\cal B}\simeq \bf 1^{\cal B}\wedge X$ and ${}^{\cal B}X \simeq {}^{\cal B}\bf 1\wedge X$;
        \item \label{t7} The following are equivalent for $X \in \scr{M}$:
        \begin{enumerate}
          \item $X$ is $\cal B$-local;
          \item $X\simeq \bf 1^{\cal B}\wedge X$;
          \item ${}^{\cal B}\bf 1\wedge X\simeq 0$;
          \item $B_j\wedge V\wedge X\simeq 0$ for all $j\in J$ and all $V \in G$
        \end{enumerate} 
        \item \label{t2} $\bc{\vee_j B_j}=\bc {\vee_{G \in G, j \in J} V\wedge B_j}=\bc{{}^{\cal B}\bf1}$;
        \item \label{t3} $\bc{\vee_j B_j}^c=\bc {\vee_{G \in G, j \in J} V\wedge B_j}^c=\bc{\bf1^{\cal B}}$;
        \item \label{t4} The multiplication map $\bf1^{\cal B} \wedge \bf1^{\cal B} \ra \bf1^{\cal B}$ is an equivalences;
        \item \label{t5} an object $X$ is $\vee_j B_j$-acyclic if and only if $X$ is a $\bf1^{\cal B}$-module if and only if $X$ is $\bf 1^{\cal B}$-local.
      \end{enumerate}
    \end{prop}

    \begin{proof}

      We start with \eqref{t1}. For this note that, if $X\in \scr{B'}$ so is $X\wedge Y$ for every $Y\in \scr{M}$. Indeed an object $Z$ is in $\scr{B'}$ if and only if $\Map(V \wedge B_j,Z)\simeq \ast$ for every $j\in J$ and every $V \in G$. On its turn this is equivalent to the condition that $D(B_j)\wedge Z\simeq 0$ for every $j\in J$. The latter condition is clearly stable under smashing over with an arbitrary object of $\scr{M}$. As a consequence ${}^{\cal B}X\wedge Y^{\cal B}$ is both local and co-local (for every $X,Y \in \scr{M}$), and thus equivalent to $0$.

      Now the rest easily follows. \eqref{t4} follows from \eqref{t1} and the fact that $\lambda^{\cal B}: \bf 1\ra \bf 1^{\cal B}$ is a unit for the algebra $\bf 1^{\cal B}$. \eqref{t6} follows from the fact that by design the localization and co-localisation functors are exact. \eqref{t7} is a direct consequence of \eqref{t1}, \eqref{t6} and of the fact the $-\wedge-$ commutes with colimits. \eqref{t2} follows from \eqref{t7}. \eqref{t3} follows from \eqref{t1} and \eqref{t2}. The first implication of \eqref{t5} follows form \eqref{t7} and \eqref{t2}. For the second implication note that $\bf 1^{\cal B}$-modules are obviously $\bf 1^{\cal B}$-local; on the other hand, $\bf 1 ^{\cal B}\wedge X$ is a $\bf 1^{\cal B}$-localization, since the multiplication map of $\bf 1^{\cal B}$ is an equivalence.
      
    \end{proof}

    \begin{defin}
      \label{defin:inverting_one_hom_element} 
     Let $\cal J=\{g_j\}_{j\in J}$ be a collection of maps in $\scr{M}$ of the form $g_j: L_j \ra \bf 1$. Let $B_j:=C(g_j)$ and consider the set of objects $\cal B=\{B_j\}_{j\in J}$. In this case the functor $(-)^{\cal B}$ (resp. the map $\lambda^{\cal B}$) is called $\cal J$-inversion functor (resp. map). In this special case we use the following notational convention: $(-)^{\cal B}=(-)[\cal J^{-1}]$ and $\lambda^{\cal B}=\iota_{\cal J}$.
     We call \emph{$\cal J$-inverted Moore object} the object $\bf 1[\cal J^{-1}]$.
     \end{defin}

    \subsubsection{} The topic of inverting homotopy elements has been treated already in the language of $\infty$-categories, for instance in \cite{MR4124487}. There Hoyois carries over a precise discussion about the possibility of describing $X[\cal J^{-1}]$ in terms of a telescope construction. This is relevant to our discussion, since his discussion directly implies the following.

   \begin{cor}
    \label{cor:inv_elts_is_t_exact}
    Assume that $\scr{M}$ is a presentably symmetric monoidal stable $\infty$-category equipped with an accessible multiplicative $t$-structure. Let $J$ be a set and let $\cal J=\{g_j\}_{j\in J}$ be a collection of maps in $\scr{M}$ of the form $g_j: L_j \ra \bf 1$, where for every $j\in J$, $L_j$ is a tif object of $\scr{M}$. The functor $(-)[\cal J^{-1}]$ is right $t$-exact. If in addition that the $t$-structure on $\scr{M}$ is compatible with filtered colimits, then $(-)[\cal J^{-1}]$ is $t$-exact.
    \end{cor}

    \begin{proof}
      Since the objects $L_i$ are $\wedge$-invertible, we can combine Lemma 3.2, Theorem 3.6, and the discussion before Lemma 3.4 of \cite{MR4124487} to ensure that $X[\cal J^{-1}]$ can be expressed via a filtered colimit, whose terms are of the form 
      \[\wedge_{j\in H}^s D(L_j)^{\wedge a_j}\wedge X,\]
      where $H$ varies among the collection of finite subsets of $J$. In addition, the flatness of the objects $L_j$, and thus of $D(L_j)$, implies that each terms of the diagram is at least as connected as $X$; in particular $(-)[\cal J^{-1}]$ respects $\scr{M}_{\geq 0}$. On the other hand $(-)[\cal J^{-1}]$ respects $\scr{M}_{\leq 0}$ when the $t$-structure is compatible with filtered colimits.
    \end{proof}

    \begin{cor}
    \label{cor:Cf_has_a_comp}
    Let $I$ be a finite set and $J$ be any set. Assume that we have a pair of collections 
    \[\{L_i\}_{i\in I} \; \textrm{ and } \; \{L_j\}_{j \in J}\]
    of strongly dualizable objects of $\scr{M}$ and a pair of collections of maps 
    \[\{f_i: L_i\ra \bf 1\}_{i \in I}\; \textrm{ and } \; \cal J=\{g_j:L_j \ra\bf 1\}_{j\in J}.\]
    Denote by $M$ the object $M(\underline f) \wedge \bf 1[\cal J^{-1}]$. Then $\bc{M}$ has a complement, given explicitly by
    \[\bc M ^{c}=\bc{\vee_{i\in I} \bf1[f_i^{-1}]} \vee \bc{\vee_{j\in J} C(g_j)}.\]
   \end{cor}

    \begin{proof}
    Immediate consequence of formula \eqref{eqn:bc_complements_formulas} together with Proposition
    \ref{prop:propert_of_trivializ}.\eqref{t2}.
    \end{proof}

    \begin{rmk}
    \label{rmk:Moore_lattice}
    In the discussion carried out around Proposition 2.13 of \cite{MR543337}, Bousfield observes that Moore spectra span a subset of $\cal{BA}(\scr{SH})$ isomorphic to the power set of $\spec(\bb Z)$. We believe if would be interesting to investigate a better (than our ad hoc) notion of Moore objects in the motivic category $\scr{SH}(S)$ which involved Thornton's computation of the homogeneous spectrum of $K^{MW}_\ast(K)$ \cite{MR3503978}.
    \end{rmk}

\section{Localization at some homotopy commutative algebras}
  \label{sec:localizations}
  In this section we prove our main results about homology localizations. We dedicate \ref{sub:fracture_squares} to some preliminary results and \ref{sub:tecass} to the statement of our main technical assumption. These are later used along \ref{sub:comparision_of_E_and_pi_0E-loc} and in the proof of Theorem \ref{thm:red_to_moore_spt}, which is the main result of this section. Throughout this section $\scr{M}$ is a presentably symmetric monoidal stable $\infty$-category, endowed with a left-complete accessible multiplicative $t$-structure.

 \subsection{Fracture squares} 
    \label{sub:fracture_squares}
  
    \begin{prop}
      \label{prop:fracture_square_for_mod_eta_cohomologies}
      Let $E,F,G$ be objects of $\scr{M}$ such that $\bc{E}=\bc{F}\vee\bc{G}$ and such that every $G$-local object is $F$-acyclic. Consider the square
      \[
      \xymatrix{
        X_E \ar[r]^{\lambda^E_{F}} \ar[d]^{\lambda ^E_{G}} & X_{F} \ar[d]^{i} \\
        X_{G} \ar[r]^{c} & (X_{F})_{G}.\\
      }
      \] where the maps $\lambda^E_F$ and $\lambda^E_G$ are induced by $\bc{F}\leq \bc{E}\geq \bc{G}$, $c=\lambda_{G}(\lambda_{F}(X))$ is the $G$-localization of the map $\lambda_{F}(X): X \ra X_{F}$, and finally $i=\lambda_{G}(X_{F})$. In this situation the above square is cartesian.

      In particular let $E\in \scr{M}$, $f:L\ra \bf 1$ be a map in $\scr{M}$ with strongly dualizable source $L$, and denote by $E/f:=E\wedge C(f)$. Then for every object $X \in \scr{M}$ the square
      \begin{equation}
        \label{eqn:frc_square_for_mod_eta_cohomologies}
      \xymatrix{
        X_E \ar[r]^{\lambda^E_{E/f}} \ar[d]^{\lambda ^E_{E[f^{-1}]}} & X_{E/f} \ar[d]^{i} \\
        X_{E[f^{-1}]} \ar[r]^{c} & (X_{E/f})_{E[f^{-1}]}.\\
      }
      \end{equation}
      is cartesian.
    \end{prop}

    \begin{proof}
      We start by denoting by $P(X)$ the pull-back of the diagram
      \begin{equation}
        \label{eqn:part_frc_square_for_mod_eta_cohomologies}
      \xymatrix{
         & X_{F} \ar[d]^i \\
        X_{G} \ar[r]^c & (X_{F})_{G}.\\
      }
      \end{equation}
      Since $\bc{F}\leq \bc{E}\geq \bc{G}$ all the objects appearing in the diagram \eqref{eqn:part_frc_square_for_mod_eta_cohomologies} are $E$-local, and thus $P(X)$ is $E$-local as well. We are thus left to prove that the natural map $u: X \ra P(X)$, which is induced by the localization maps $\lambda_{F}(X): X \ra X_{F}$ and $\lambda_{G}(X): X\ra X_{G}$, is an $E$-equivalence. Since $\bc{E}=\bc{F}\vee\bc{G}$ it suffices to show that $u$ is both an $F$-equivalence and a $G$-equivalence. In particular, by the 2-out-of-3 property of local equivalences, we reduce to show that
      \begin{itemize}
        \item the natural map $P(X)\overset{\alpha}{\ra} X_{G}$ is a $G$-equivalence, which is obvious;
        \item the natural map $P(X)\overset{\beta}{\ra} X_{F}$ is a $F$-equivalence, which follows form the assumption that $G$-local objects are $F$-acyclic.
      \end{itemize}
      The second part of the statement follows easily: set $F=E/f$, $G=E[f^{-1}]$, and notice that by Proposition \ref{prop:propert_of_trivializ}\eqref{t3} $\bc{\bf 1}= \bc{C(f)}\vee\bc{\bf 1[f^{-1}]}$, so that $\bc{E}=\bc{E/f}\vee\bc{E[f^{-1}]}$. Finally use the chain of inclusions
      \[\rm{Loc}_{E[f^{-1}]}\subseteq \rm{Loc}_{\bf1[f^{-1}]}= \rm{Ac}_{C(f)}\subseteq \rm(Ac)_{E[f^{-1}]}:\]
      the only non-obvious part is the central equality, which follows from Proposition \ref{prop:propert_of_trivializ}\eqref{t5}.
    \end{proof}

    \begin{cor}
      \label{cor:E/eta_loc_is_E_cofeta_loc}
      Let $E$, $X$ and $f$ be as in the statement of Proposition \ref{prop:fracture_square_for_mod_eta_cohomologies}. Then $X_{E\wedge C(f)}$ is equivalent to $(X_E)_{C(f)}$.
    \end{cor}
    \begin{proof}
      We start by considering the following commutative diagram in $\scr{M}$:
      \begin{equation}
        \label{eqn:E/eta_loc_is_E_cofeta_loc}
      \xymatrix{
       X_E \ar[r]^{\lambda_1}\ar[d]^{\lambda_3} & (X_E)_{C(f)}\ar[d]^{\lambda_4}\\
       X_{E/f} \ar[r]^{\lambda_2} & (X_{E/f})_{C(f)}\\
      }
      \end{equation}
      where $\lambda_1=\lambda_{C(f)}(X_E)$, $\lambda_3=\lambda^E_{E/f}$ is the map introduced above, $\lambda_4=\lambda_{C(f)}(\lambda_3)$, and finally $\lambda_2=\lambda_{C(f)}(X_{E/f})$.
      
      Since $\bc{E/f}\leq\bc{C(f)}$, $\lambda_2$ is actually an equivalence in $\scr{M}$ and we are left to prove that $\lambda_4$ is too. For this we apply the $C(f)$-localization functor to the square \eqref{eqn:frc_square_for_mod_eta_cohomologies} and we use \ref{prop:fracture_square_for_mod_eta_cohomologies} to reduce the proof to checking that 
      \[X_{E[f^{-1}]} \overset{c}{\ra} (X_{E/f})_{E[f^{-1}]} \] is a $C(f)$-equivalence. Now both the source and target of $c$ are $\bf 1[f^{-1}]$-local, being in fact $E[f^{-1}]$-local. Thus, after $C(f)$-localization both source and target of $c$ become zero.
    \end{proof}

    \begin{cor}
      \label{cor:E/x,y-loc_is_E_cofx_cof_y_loc}
      Let $r$ be a positive integer. For every $i=1, \dots, r$ consider strongly dualizable objects $L_i \in \scr{M}$ and maps $f_i: L_i\ra \bf 1$ in $\scr{M}$. Let $M:=C(f_1)\wedge \cdots \wedge C(f_r)$ be the Moore object associated to the maps $f_1, \dots ,f_r$. Then for every pair of objects $E$ and $X$ of $\scr{M}$
      \begin{equation}
        X_{E\wedge M} \simeq (\cdots((X_E)_{C(f_1)})_{C(f_2)}\cdots)_{C(f_r)}\simeq (X_E)_M.
      \end{equation}
    \end{cor}
    
    \begin{proof}
      Apply inductively Corollary \ref{cor:E/eta_loc_is_E_cofeta_loc}.
    \end{proof}

  \subsection{Technical Assumption} 
    \label{sub:tecass}
    We now formulate a technical assumption. This is going to be required for our proves of Theorem \ref{thm:red_to_moore_spt}, \ref{thm:convergence_mod_stuff} and \ref{thm:convergence_inv_stuff} to work.

    \begin{ass}
        \label{sub:assumption_A1}

        $S$ is a base scheme and $E$ is an object of $\scr{M}$ satisfying the following properties:
      \begin{enumerate}
        \item \label{ass1:E_hom_com} $E$ is a homotopy commutative algebra;
        \item \label{ass2:E_con} $E\in \scr{M}_{\geq 0}$;
        \item \label{ass3:phi}there is a finite set of tif objects $\{L_i\}_{i \in I}$ in $\scr{M}$ and maps $f_i: L_i \ra \bf 1$,
        a countable set of tif objects $\{L_j\}_{j \in J}$ and maps $g_j: L_j \ra \bf 1$,
        and a morphism of $\tau_0(\bf 1)$-algebras
        \[\varphi: \tau_0 \Big ( (\tau_0(\bf 1)/\cal I)[\cal J^{-1}] \Big )  \rightarrow \tau_0 E,\]
        where $\cal I$ is the ideal of $\tau_0(\bf 1)$ generated by $\{f_i\}_{i \in I}$ and $\cal J$ is the collection of elements $\{g_j\}_{j \in J}$.
        \item \label{ass4:phi_iso} The map $\varphi$ of point \eqref{ass3:phi} is an isomorphism.
      \end{enumerate}

      In this situation we will denote by $\bf 1\tau_0(E)$ the object $C(f_1)\wedge \cdots\wedge C(f_r)\wedge \bf 1[\cal J^{-1}]$, and refer to it sloppily as of the Moore object associated with $\tau_0(E)$, rather then with the maps $f_i$'s and $g_j$'s.

      Note that the natural map of commutative algebras $(\tau_0(\bf 1)/\cal I)[\cal J^{-1}] \ra \tau_0 \big ( (\tau_0(\bf 1)/\cal I)[\cal J^{-1}] \big )$ is an equivalence as soon as $J=\emptyset$ or when the $t$-structure on $\scr{M}$ is compatible with filtered colimits.
      \end{ass}


  \subsection{Comparison of $E$- and $\tau_0E$-localization} 
    \label{sub:comparision_of_E_and_pi_0E-loc}
    
    \begin{lemma}
        \label{lem:mult_on_ring_sp}
        Let $E\in \scr{M}_{\geq k}$ be a homotopy commutative algebra in $\scr{M}$ and $f: L \ra \bf 1$ be a map in $\scr{M}$ whose source $L$ is an sdf object. Then left multiplication by $f$ on $E$, $l_f(E)$, is an equivalence (resp. homotopic to zero) if and only if left multiplication by $f$ on $\tau_0(E)$, $l_f(\tau_0(E))$, is an equivalence (resp. zero).
      \end{lemma}
      \begin{proof} If $l_f(E)$ is an equivalence (resp. homotopic to zero), the induced map $\tau_0(l_f(E))$ is too, and by Lemma \ref{lemma:left_mult_on_heart} $\tau_0(l_f(E))=l_F(\tau_0(E))$, as morphisms of $\scr{M}^{\heartsuit}$.
      For the other implication we argue as follows. By left-completeness and the fact that $L$ is a compact object, the map $l_f(E): L\wedge E \ra E$ is naturally identified with the limit of $l_f(P^n(E)):L\wedge P^n(E) \ra P^n(E)$, and for checking that it is an equivalence (resp. homotopic to zero) one only needs to check that, for every integer $n$, multiplication by $f$ on $\tau_n(E)$ is an equivalence (resp. homotopic to zero). Now $\tau_n(E)$ is a $\tau_0(E)$-module and thus is endowed with an action map $a:\tau_0(E)\otimes^{\heartsuit}\tau_n(E) \ra \tau_n(E)$. Multiplication by $f$ commutes with such action map by \ref{subsub:right_and_left_mult}, in the sense that one can find a commutative square
        \[\xymatrix{
        L\wedge \tau_n(E) \ar[r]^{l_f(\tau_k(X))} & \tau_n(E)\\
        L\wedge \tau_0(E)\otimes^{\heartsuit}\tau_n(E) \ar[r]\ar[u]_a & \tau_0(E)\otimes^{\heartsuit}\tau_n(E)\ar[u]_a,\\
        }\]
        where the lower horizontal map is $l_f\big (\tau_0(E)\otimes^{\heartsuit}\tau_n(E)\big )$. However, such map is homotopic to $l_f\big (\tau_0(E)\big ) \otimes^{\heartsuit}\tau_n(E)$, and thus is an equivalence (resp. homotopic to zero) by assumption. We conclude noticing that the upper horizontal map of the above square is a retraction of the lower horizontal map via a choice of unit morphism for the $\tau_0(E)$-module $\tau_n(E)$.
      \end{proof}

    \begin{prop}
      \label{prop:fund_ineq_of_loc}
      Let $E$ be an object of $\scr{M}$ satisfying points \eqref{ass1:E_hom_com}-\eqref{ass3:phi} of Assumption \ref{sub:assumption_A1}. Then:
      \[\bc{\tau_0E} \; \leq  \; \bc{E} \; \leq \; \bc{M\tau_0E}.\]
      \end{prop}
    
    \begin{proof}
      Since $E$ is $(-1)$-connected by Assumption \ref{sub:assumption_A1}, the projection to the Postnikov truncation induces a map $p:E\ra \tau_0E$. As it follows from the discussion in \ref{subsub:trunc_of_hom_mod}, $\tau_0E$ has a natural homotopy algebra structure for which $p$ is a ring map. As a consequence $\tau_0(E)$ is a retract of $E\wedge \tau_0(E)$, and it follows that $\bc{\tau_0(E)}=\bc{E\wedge \tau_0(E)}\leq \bc{E}$.

      It remains to show that $\bc{E}= \bc{E\wedge M\tau_0 E}$, which directly implies our claim that $\bc{E}\leq \bc{M\tau_0 E}.$
      For this, recall that
       \[\bc{M\tau_0E}=\bc{C(f_1)\wedge \cdots \wedge C(f_n) \wedge \bf 1[\cal J^{-1}]}.\] 
      From Proposition \ref{cor:Cf_has_a_comp}, the Bousfield class $\bc{M\tau_0E}$ has a complement $\bc{M\tau_0E}^c=\bc{M}$, where 
        \[ M:=\big (\bigvee_{i\in I} \bf 1[f_i^{-1}] \big ) \vee \big(\bigvee_{j\in J} C(g_j)\big).\]
      In particular
      \[ \bc{E} = \bc{E}\wedge \bc{(M\tau_0E \vee M)} = \bc{E\wedge M\tau_0E} \vee \bc{E\wedge M};\]
      since smashing commutes with small sums, it suffices to see that $E\wedge \bf 1[f_i^{-1}]=0=E\wedge C(g_j)$ for every $i\in I$ and every $j\in J$.
      
      Let $f:L\ra\bf 1$ be any of the $f_i$'s and let $l_f(E')$ be left multiplication by $f$ on $E'=E\wedge \bf 1[f_i^{-1}]$. Then $l_f(E')$ is homotopic to the left multiplication induced from $E$, $l_f(E)\wedge \bf 1[f^{-1}]$, but is also homotopic to the left multiplication induces form $\bf 1[\frac{1}{f}]$, $E\wedge l_f(\bf 1[\frac{1}{f}])$. The former is homotopic to zero by Lemma \ref{lem:mult_on_ring_sp}, while the latter is an equivalence by construction. It follows that $E\wedge\bf 1[\frac{1}{f}]\simeq 0$. Similarly let $g:L\ra \bf 1$ be any of the $g_j$'s. Multiplication by $g$ on $E$ is an equivalence by Lemma \ref{lem:mult_on_ring_sp} since it is so on $\tau_0(E)$, and thus $E\wedge C(g)\simeq 0$.
    \end{proof}

    \begin{lemma}
      \label{prop:con_implies_HK^MW-loc}
      Let $X$ be any $k$-connected object of $\scr{M}$ for some integer $k$. Then $X$ is $\tau_0(\bf 1)$-local.
    \end{lemma}
    \begin{proof}
      By the multiplicative properties of the Postnikov tower \ref{ssub:mult_prop_t_str}, the homotopy objects $\tau_p(X)$ are $\tau_0 (\bf 1)$-modules ($P^0(\bf 1)\simeq \tau_0(\bf 1)$ since $\bf 1$ is $(-1)$-connected) and hence $\tau_0 (\bf 1)$-local. By connectivity of $X$ every stage of the Postnikov tower $P^n(X)$ is a finite extension of suspensions of the $\tau_p(X)$'s, and hence it is $\tau_0 (\bf 1)$-local. Finally, $X\simeq \lim P^n(X)$ by left completeness of the $t$-structure.     
    \end{proof}

    \begin{lemma}
    \label{lemma:ker_and_coker_of_multiplications}
       Let $\{\bf A,\otimes,\bf 1\}$ be a symmetric monoidal category where $\bf A$ is also abelian and such that $\otimes$ is right exact. Let $R$ be a commutative monoid in $\bf A$ and denote by $e_R$ and $\mu_R$ its unit and multiplication. Let $L \in \bf A$ and assume we have a map $f: L \ra R$. Then:
       \begin{enumerate}
         \item \label{ck.1}there exist a unique $R$-linear map $\cdot f: R\otimes L \ra R$ that composed with $L\simeq \bf 1\otimes L \overset{\op{e}_R\otimes\op{id}_L}{\ra} R\otimes L$ gives back $f$
         \item \label{ck.2} the $R$-module $C:=\coker (\cdot f)$ has a unique structure of commutative $R$-algebra having the natural projection $p: R\ra C$ as unit;
         \item \label{ck.3} when the multiplicaiton of $R$ is an isomorphism, so is that of $C$;
         \item \label{ck.4} when $L$ is a $\otimes$-invertible object of $\bf A$, $K:=\ker (\cdot f)$ has a natural structure of $C$-module.
         \item \label{ck.5} if $M$ is an $R$-module in $\bf A$ with action map $\alpha: R\otimes M\ra M$ and 
         \[\cdot f_M=\alpha\circ (f\otimes \rm{id}): L\otimes M \ra R\otimes M \ra M\] is the induced multiplication by $f$ on $M$, then $\coker{(\cdot f_M)}$ and $\ker{(\cdot f_M)}$ have a unique structure of $C$-module induced by $\alpha$.
       \end{enumerate}
     \end{lemma}
     \begin{proof}
       Point \eqref{ck.1} follows form the usual free-forget adjunction. Point \eqref{ck.2} and \eqref{ck.3} follows form elementary diagram chases. Regarding \eqref{ck.4} one needs to show that the natural action of $R$ on $K$ factors through $C$. Let us denote by $i$ the monomorphism $K\subseteq R\otimes L$. Consider the map 
       \[\phi:R\otimes L\otimes R\otimes L \ra R\otimes L\; \textrm{ defined by }\; \phi=\mu\otimes \rm{id}_L\circ \cdot f\otimes \rm{id}_R\otimes\rm{id}_L,\]
      so that, with an abuse of notation, $\phi(r_1\otimes l_1\otimes r_2 \otimes l_2)=r_1f(l_1)r_2 \otimes l_2$.
      An easy diagram chase shows that the required factorization exists if and only if the map $\phi$
       composed with $j:=\op{id}\otimes\op{id}\otimes i:R\otimes L\otimes K\ra R\otimes L\otimes R\otimes L$ is zero. What we know, however, is that $\phi\circ t_{(12),(34)}\circ j=0$, where $t_{(12),(34)}$ is the switching of the first and second pair of terms on $R\otimes L\otimes R\otimes L$. We claim that $\phi$ is a multiple of $\phi\otimes t_{(12),(34)}$ by a unit $\epsilon \in \rm{End}_{\bf A}(\bf 1)$. For this note that the permutation $t_{(12),(34)}=t_(13)\circ t_(24)$ of the $R$ terms gives no trouble, since $\mu$ is commutative, so $\phi=\phi\circ t_{(13)}$. Since $L$ is $\otimes$-invertible, $t_{(24)}$ can be identified with left multiplication by a suitable unit $\epsilon \in \op{End}_{\bf A}(\bf 1)$. In conclusion $\phi\circ t_{(12),(34)}=\phi\circ l_\epsilon(R\otimes L\otimes R\otimes L)$, and thus 
       \[\phi\circ j= \phi\circ t_{(12),(34)} \circ l_\epsilon(R\otimes L\otimes R\otimes L)^{-1}\circ j=\phi\circ t_{(12),(34)}\circ j \circ l_\epsilon(R\otimes L\otimes R\otimes L)^{-1}=0, \]
       since multiplication by a unit commutes with every map in $\bf A$ by \ref{subsub:right_and_left_mult}. The proof of point \eqref{ck.5} works similarly, and we omit it.
     \end{proof}

    \begin{lemma}
      \label{lemma:HR/x_vs_H_R/x}
      Let $\cal R$ be an commutative algebra in $\scr{M}^{\heartsuit}$ and let $f: \tau_0(L) \ra \cal R$ be a map in $\scr{M}^{\heartsuit}$ where $L$ is a tif object of $\scr{M}$. Let $C$ denote the cofiber of the induced map $\cdot f:\cal R\wedge L\simeq \cal R\otimes^\heartsuit \tau_0(L) \ra \cal R$ given by \ref{lemma:ker_and_coker_of_multiplications} in $\scr{M}^{\heartsuit}$.
      Then $\tau_0 (C) \simeq \coker{f\cdot}$ and  $\bc{C} \leq \bc{\tau_0(C)}$. 
    \end{lemma}

    \begin{proof}
      From the exact sequence of homotopy objects
      \begin{equation}
        0 \ra \cal K \ra L\wedge \cal R\overset{\cdot f}{\ra} \cal R \ra \coker(\cdot f)\ra 0
      \end{equation}
      we deduce 
      \begin{equation}
        \underline \tau_k (C)=
        \begin{cases}
          \coker(f\cdot) & \text{if $k=0$,}\\
          \cal \ker(f\cdot) & \text{if $k=1$, }\\
          0 & \text{else.}
        \end{cases}
      \end{equation}
      In particular it follows that we have a fiber sequence 
      \begin{equation}
        \Sigma ^1 \tau_1C \ra C \ra \tau_0 C.
      \end{equation}
      relating $C$ with its truncations. Observe that, thanks to Lemma \ref{lemma:ker_and_coker_of_multiplications}, $\tau_0(C)$ is a commutative algebra in $\scr{M}^{\heartsuit}$ and $\Sigma ^1 \tau_1 (C)$ is a module in $\scr{M}^{\heartsuit}$ over $\tau_0(C)$. More generally for every $X\in \scr{M}$, the homotopy object $ \tau_1 (C) \wedge X$ is also a $ \tau_0 (C)$-module in $\scr{M}$. In particular, if $X$ is $\tau_0(C)$-acyclic then it is also $\Sigma^1\tau_1 (C)$-acyclic, and by the above fiber sequence $X$ is $C$-acyclic too.
    \end{proof}

    \begin{cor}
      \label{cor:HK^MW_mod_stuff_vs_H_KMW_mod_stuff}
      Assume that $\scr{M}$ is a presentably symmetric monoidal stable $\infty$-category endowed with an accessible left-complete multiplicative $t$-structure. Assume that $E$ is a homotopy commutative algebra in $\scr{M}$ satisfying assumption \ref{sub:assumption_A1}, and in case $J\not = \emptyset$ we assume in addition that the $t$-structure of $\scr{M}$ is compatible with filtered colimits. Then $\bc{\tau_0 \bf 1\wedge M\tau_0(E)} \leq \bc{\tau_0E}$.
    \end{cor}

    \begin{proof}
      Readily \ref{lemma:HR/x_vs_H_R/x} implies that $\bc{\tau_0(\bf 1)\wedge C(f_1)}\leq \bc{\tau_0(\bf 1)/(f_1)} $. 
      Recall that by construction of Moore object $M\tau_0E=C(f_1)\wedge \cdots \wedge C(f_r)\wedge\bf1[J^{-1}]$, so we can proceed in order by smashing with one $C(f_i)$ at the time. Indeed by the previous case
      \[ \bc{\tau_0(\bf 1)\wedge C(f_1) \wedge  C(f_2)} \leq \bc{\tau_0(\bf 1)/(f_1) \wedge C(f_2)}\]
      and finally, using \ref{lemma:HR/x_vs_H_R/x} with $\cal R=\tau_0(\bf 1)/(f_1)$ and $f=f_2$, we get that
      \[ \bc{\tau_0(\bf 1)/(f_1)) \wedge C(f_2)} \leq \bc{ \big ( \tau_0(\bf 1)/(f_1)\big )/f_2))} =\bc{\tau_0(\bf 1)/(f_1,f_2)}\]
      so we conclude that 
      \[\bc{\tau_0(\bf 1)\wedge C(f_1) \wedge  C(f_2)} \leq  \bc{\tau_0(\bf 1)/(f_1,f_2)}.\] Inductively we arrive to 
      \[\bc{\tau_0(\bf 1) \wedge M(\underline f)}\leq \bc{\tau_0(\bf 1)/(f_1,\dots,f_r)},\]
      where $M(\underline f)=C(f_1)\wedge \cdots \wedge C(f_r)$.
      Finally we observe that 
      \[\bc{\tau_0(\bf 1) \wedge M\tau_0 E} = \bc{\tau_0(\bf 1) \wedge M(\underline f)\wedge \bf 1[\cal J^{-1}]} \leq \bc {\tau_0 (\bf 1/(f_1,\dots,f_r)) \wedge \bf 1[\cal J^{-1}]} = \bc{\tau_0 E },\]
      where the last inequality comes form the fact that $-\wedge \bf 1[\cal J^{-1}]$ is $t$-exact according to Corollary \ref{cor:inv_elts_is_t_exact}.
    \end{proof}

    \begin{thm}
    \label{thm:red_to_moore_spt}
      Let $E$ be a homotopy commutative algebra in $\scr{M}$ satisfying assumption \ref{sub:assumption_A1} in the special case that $J=\emptyset$. Then for every integer $k$, and every $k$-connected $X$ we have that
      \[X_{M\tau_0E}\simeq X_{E}.\]
     
    \end{thm}

    \begin{proof}
      Thanks to Proposition \ref{prop:fund_ineq_of_loc} we know that any $k$-connected $X$ the localization map $X \ra X_{M\tau_0E}$ is an $E$-equivalence so we only have to check that $X_{M\tau_0E}$ is $E$-local. Now consider that 
      \[
      X \ra X_{\tau_0(\bf 1)}
      \] is an equivalence by Proposition \ref{prop:con_implies_HK^MW-loc} so that 
      \[
      X_{M\tau_0E} \ra (X_{\tau_0(\bf 1)})_{M\tau_0E}
      \] is an equivalence too.
      In particular, by combining this with the result of Corollary \ref{cor:E/eta_loc_is_E_cofeta_loc} we deduce that
      \[
      X_{M\tau_0E} \overset{\simeq}{\ra} (X_{\tau_0(\bf 1)})_{M\tau_0E}  \simeq X_{\tau_0(\bf 1)\wedge M\tau_0E}.
      \] 
      Finally we apply Proposition \ref{prop:fund_ineq_of_loc} and \ref{cor:HK^MW_mod_stuff_vs_H_KMW_mod_stuff} to deduce that $X_{\tau_0(\bf 1)\wedge M\tau_0E}$ is $E$-local and this concludes.
    \end{proof}
  

\section{Examples and applications} 
  \label{sec:examples}
  
  We provide some samples of usage in the motivic setting of the results seen so far. In this section $S$ denotes a Noetherian scheme of finite Krull dimension.

  \subsection{Algebraic Cobordisms} 
    \label{ssub:subs_ex_of_loc_mgl}
    Let $\op{MGL}$ (resp. $\op{MSL}$) be the spectrum representing Voevodsky's Algebraic Cobordism (resp. Special Linear 
    Cobordism) over $S$. In \cite{panin_pimenov_rondigs:a_universality_theorem_for_voevodsky_MGL} (resp.  \cite{alg_bor_sp}) the authors construct $\op{MGL}$ (resp. $\op{MSL}$) as a commutative monoid in a monoidal model category presenting $\scr{SH}(S)$. By Theorem 3.8 of \cite{hoyois:from_algebraic_cobordism_to_motivic_cohomology} $\op{MGL}$ satisfies assumption \ref{sub:assumption_A1}; we conclude that the $\op{MGL}$-localization map is canonically identified with the $\eta$-completion map $\chi_\eta(X): X \ra X^\wedge_\eta$ on spectra that are $k$-connected for some $k$. Similarly in \cite{MSL}, when $S$ is the spectrum of a perfect field $ $ of characteristic not $2$, the author verifies that $\op{MSL}$ satisfies assumption \ref{sub:assumption_A1} with $I=J=\emptyset$, and as a consequence $(-)_{\op{MSL}}\simeq \op{id}$ on $k$-connected spectra.

  \subsection{Motivic Cohomologies}
    \label{ssub:ex_of_loc_HZ}
    Let $K$ be a perfect field of charateristic not $2$. Let $\op H\bb Z$ be the spectrum representing Voevodsky's Motivic cohomology with integral coefficients.
    Recall that we have a category of motives $\scr {DM}(K,\bb Z)$ which is related to $\scr {SH}(K)$ by an adjunction
    \begin{equation}
    \label{eqn:add-forget_transfers}
     \bb Z_{tr}:\scr {SH}(K) \rightleftarrows \scr{DM}(K,\bb Z):  u_{tr}.
    \end{equation}
    Since $u_{tr}$ respects algebras, as it follows from Section 4 of \cite{hoyois:from_algebraic_cobordism_to_motivic_cohomology}, $\op H\bb Z=u_{tr}{\bf 1}$ is a commutative algebra in $\scr{SH}(S)$. That $\op H\bb Z$ is $(-1)$-connected follows directly from the representability of motivic cohomology in Proposition 14.16 of $\scr {DM}(K,\bb Z)$ \cite{mazza_voevodsky_weibel;lecture_notes_on_motivic_cohomology} combined with Theorem 3.6 of \cite{mazza_voevodsky_weibel;lecture_notes_on_motivic_cohomology}. Furthermore, if we denote by $e$ the unit of the ring spectrum $\op H\bb Z$ and by $\sigma$ Morel's isomorphism, we can form the following diagram
    \[
    \xymatrix{
    \underline \pi_0(\bb S) \ar[r]^e &  \underline \pi_0(\op H\bb Z)\\
    \cal K^{MW} \ar[u]^{\sigma}_{\simeq} \ar[r] &  \cal K^{MW}/(\eta)=\cal K^M \ar[u]^\lambda, \\
    }
    \]
    since $\eta$ acts trivially on motivic cohomology. In order to check that $\lambda$ is an isomorphism we immediately reduce, using Theorem 2.2 of \cite{MR2275634} to check that the induced map
    \[\lambda_L: \bigoplus_ {n\in \bb N} K^M_n(L) \overset{\simeq}{\ra} \bigoplus_{n \in \bb N} \op H^n(\spec L,\bb Z(n))\]
    on stalks at finitely generated field extensions $L/K$ is an isomorphism. An immediate direct check unveils that $\lambda_L$ is the map of Theorem 5.1 of \cite{mazza_voevodsky_weibel;lecture_notes_on_motivic_cohomology}, and is thus an isomorphism.
  
    We can thus apply Theorem \ref{thm:red_to_moore_spt}, to deduce that the $\op H\bb Z$-localization of a $k$-connected  spectrum $X$ is identified with the $\eta$-completion map $X \ra X^\wedge_\eta$. In a similar fashion let $E=\op H\bb Z/\ell$ be the spectrum representing motivic cohomology with modulo $\ell$ coefficients. The same considerations allow us to conclude that the $\op H\bb Z/\ell$-localization map $\lambda_{\op H\bb Z/\ell}(X): X \ra X_{\op H\bb Z/\ell}$ of a $k$-connected spectrum $X$ is identified with the formal completion map $\chi_{\ell,\eta}(X): X \ra X^\wedge_{\ell,\eta}$.

    The formalism of motives we have just recalled has a quadratic analogue. We have a category of Chow-Witt motives $\widetilde{\scr {DM}}(K,\bb Z)$ with a pair of adjoint functors
    \begin{equation}
    \label{eqn:SH_tDM_adj}
       \widetilde {\bb Z}_{tr}: \scr{SH}(K) \rightleftarrows \widetilde{\scr {DM}}(K):  \widetilde u_{tr}
    \end{equation}
    which are the stabilizations of the functors which respectively add and forget generalized transfers. Quadratic versions of the arguments given above are used in Section 3.3.5 of \cite{MW-mot-comp} to show that $\widetilde{\op H}\bb Z$ is a commutative algebra in $\scr{SH}(K)$. The same argument employed above for $\op H \bb Z$ shows that the ring spectrum $\widetilde {\op H} \bb Z$ is $(-1)$-connected. This time, the unit map $\bb S \ra \widetilde {\op H} \bb Z$ induces an isomorphism on $\tau_0$ as one sees by combining the adjunction \eqref{eqn:SH_tDM_adj} with Theorem 4.2.3 of \cite{MW-mot-comp}. Note that Theorem 4.2.3 of \cite{MW-mot-comp} is exactly a quadratic version of the argument we have given above, and its main ingredient is a quadratic version of the Nesterenko-Suslin-Totaro identification of Theorem 6.19 of \cite{CF}. It thus follows by Theorem \ref{thm:red_to_moore_spt} that for every $k$-connected spectrum $X$, the $\widetilde {\op H} \bb Z$-localization is an equivalence.

  \subsection{Slice completion}
    We recall the following result.
    \begin{thm}[Theorem 3.50 of \cite{MR3898173}]
      Let $K$ be a field of exponential characteristic $p$. Suppose that $X$ is a spectrum having a cell presentation of finite type in $\scr M_A(K)$ where $A=\bb S[\frac{1}{p}]$. Then we  have a canonical commutative square in $\scr {M}_A(K)$
      \[
      \xymatrix{
         X \ar[rr]_{\chi_\eta(X)} \ar[d]^{\sigma(X)} & & X^\wedge_\eta \ar[d]^{\sigma(X)^\wedge_\eta}\\
         X_{sc} \ar[rr]^{\chi_\eta(X_{sc})} & &  (X_{sc})^{\wedge}_\eta,\\
      }
      \]
      where the maps $\chi_\eta(X_{sc})$ and $\sigma(X)^\wedge_\eta$ are equivalences. In particular there is a natural isomorphism $X_{sc} \simeq X^\wedge_\eta$ in $\scr{M}_A(K)$ under which the slice completion map $\sigma(X)$ and the $\eta$-completion map $\chi_{\eta}(X)$ are identified.
    \end{thm}

    By combining \ref{ssub:ex_of_loc_HZ} with the previous result, we deduce that for a cell spectrum $X$ of finite type, $H\bb Z\wedge X=0$ if and only if $X_{sc}=0$. In other words, if $f: X\rightarrow Y \in \scr{M}_A(K)$ is a map between objects that have a cell presentation of finite type, $f$ is an equivalence on slice completions if an only if $f$ induces an equivalence on motivic homology.

  \subsection{Motives of spectra}
    \label{ssub:ex_on_motives_of_spectra}

    Let $K$ be a field of exponential characteristic $p$ and $A=\bb S[\frac{1}{p}]$ and $R=\bb Z[\frac{1}{p}]$. Recall that the adjunction \eqref{eqn:add-forget_transfers} factors as
    \[
     \xymatrix@C=4em{
     \scr{SH}(K)[\frac{1}{p}] \simeq \scr M_{A}(K)\ar@<.5ex>[rr]^{R_{tr}} \ar@<.5ex>[dr]^{\op {HR}\wedge_A-} & &  \scr{DM}(K,R) \ar@<.5ex>[ll]^{u_{tr}} \ar@<.5ex>[dl]^{\Psi}\\
     & \scr M_{\op {HR}}(K) \ar@<.5ex>[ul]^{ U_{\op {HR}}} \ar@<.5ex>[ru]^{\Phi}  & \\
     }
    \]   
    where $R_{tr}(-)\simeq \Phi(\op{HR}\wedge -)$ and $u_{tr}\simeq U_{\op{HR}}(\Psi(-))$.
    
    Similarly, at least when $p\not = 2$, the adjunction \eqref{eqn:SH_tDM_adj} factors as
    \[
     \xymatrix@C=4em{
     \scr{SH}(K)[\frac{1}{p}] \simeq \scr M_{A}(K)\ar@<.5ex>[rr]^{\widetilde R_{tr}} \ar@<.5ex>[dr]^{\op {\widetilde HR} \wedge-} & &  \widetilde{\scr{DM}}(K,R) \ar@<.5ex>[ll]^{\widetilde u_{tr}} \ar@<.5ex>[dl]^{\widetilde\Psi}\\
     & \scr M_{\op {\widetilde HR}}(K) \ar@<.5ex>[ul]^{ \widetilde U_{\op{HR}}} \ar@<.5ex>[ru]^{\widetilde\Phi}  & \\
     }
    \]   
    The previous observations prove the following statement. 
    \begin{cor}
    \label{cor:conserv}
      Let $K$ be a field of exponential characteristic $p$ which is not formally real. Then:
      \begin{enumerate}
        \item the functor $R_{tr}$ is conservative on $k$-connected $ \bb S[\frac{1}{2}]$-modules;
        \item if in addition $cd_2(K)<\infty$ then $R_{tr}$ is conservative on strongly dualizable objects, and in particular on compact objects;
        \item when in addition $p \not = 2$, the functors $(\tilde\Phi,\tilde\Psi)$ are inverse equivalences of categories.
      \end{enumerate}
        
    \end{cor}
    \begin{proof}
      Theorem 1 of \cite{MR2435654} when $p=0$, or Theorem 5.8 of \cite{MR3730515} in general, implies that $(\Phi,\Psi)$ is a pair of adjoint equivalences. Now $X\simeq X[\frac{1}{2}]$ so let us assume that $H\bb R \wedge X=0$. Then by the assumption on $X$, Theorem \ref{thm:red_to_moore_spt}, and Proposition \ref{prop:localization_at_mod_x_moore_spectrum} we have that
      \[0=X_{H\bb Z}=X^\wedge_\eta.\]
      However Lemma \ref{lem:Xonehalf_is_eta_complete} implies that $\chi_{\eta}(X): X \ra X^\wedge_\eta$ is an equivalence in $\scr{M}_A(K)$ and hence $X\simeq 0$. Let us now assume that $K$ has finite $2$-cohomological dimension. We achieve the second point by running the same argument, but using \ref{lemma:compact_is_eta_complete} instead of \ref{lem:Xonehalf_is_eta_complete}. We deduce that, if $X$ is a dualizable object with $H\bb Z \wedge X=0$, then $X\simeq 0$ in $\scr M_A(K)$. Since $p$ is inverted in the coefficients, strongly dualizable objects and compact objects of $\scr{SH}(K)[\frac{1}{p}]$ are the same. The third point works similarly to the first, but using only Theorem \ref{thm:red_to_moore_spt}, combined with the fact that that $(\tilde\Phi,\tilde\Psi)$ are mutually inverse equivalences, which follows from Lemma 5.3 of \cite{Bac-Fas}.
    \end{proof}
    
    \begin{rmk}
      \label{rmk:our_cons_res_vs_Bachmann}

      In \cite[Theorem 1]{Bacmann_Conservativity} Bachmann proves that over a perfect field $K$ of exponential characteristic $p\not = 2$ and with $cd_2(K)<\infty$, the functor $\bf L \bb Z_{tr}$ is conservative on effective and $k$-connected spectra where $p$ acts invertibly. 
      In the case where $2$ in inverted in the coefficients, Corollary \ref{cor:conserv} is only apparently more general. A direct inspection of Bachmann's argument shows that, upon inverting two, he does not need to assume that $cd_2(K)<\infty$.

      However, if we don not wish to invert the prime $2$ in the coefficients, Bachmann's argument needs an extra non-trivial input, namely a description of the slice filtration on homotopy groups of spectra, coming from \cite{MR2918174}. In this case our approach in genuinely different and not as powerful, but in any case it recovers a very non-trivial portion of the stated result of Bachmann.

      Nevertheless we have not been able to use Bachmann's results to recover our results on $\op H\bb Z$-localizations. It would probably be interesting to employ his techniques, particularly those using the real \'etale topology, for the study of homology localizations. 
    \end{rmk}


  \subsection{$K$-theories} 
  \label{sub:ex_K-loc}
  Let $S$ be an essentially smooth scheme over a field $K$ of exponential characteristic $p$ (resp. let $S$ be the spectrum of a field of characteristic $p\not =2$).
  Recall that we have spectra $\op{KGL}$ (resp. $\op{KQ}$) in $\scr{SH}(S)[\frac{1}{p}]$ representing algebraic $K$-theory (resp. algebraic Hermitian $K$-theory). The tensor product on bundles induces natural homotopy commutative algebras structures on $\op{KGL}$ and on $\op{KQ}$. This can be found respectively in Theorem 2.2.1 of \cite{MR2597741} and in Theorem 1.5 of \cite{MR3882540}: in both cases the multiplicative structure is constructed over $\spec(\bb Z)$ and $\spec(\bb Z[\frac{1}{2}])$ respectively, and then pulled back over more general bases. Thanks to the multiplicative properties the slice tower, the effective cover $f_0(\op{KGL})$ (resp. the very effective cover $\widetilde f_0(\op{KQ})$ have a structure of homotopy commutative algebras as well. 

  We first deal with $\op{KGL}$. Consider the fiber sequence
  	\[f_1 \op{KGL} \ra f_0\op{KGL} \ra s_0 \op{KGL}.\]
  On one hand we have an isomorphism of commutative ring spectra $s_0 \op{KGL} \overset{\simeq}{\ra} \op H\bb Z$, see \cite[Section 11]{MR2365658} and Theorem 8.5 of \cite{hoyois:from_algebraic_cobordism_to_motivic_cohomology}. On the other hand $f_1\op{KGL}\in \scr{SH}(K)[\frac{1}{p}]_{\geq 1}$ thanks to Lemma 8.11 of \cite{hoyois:from_algebraic_cobordism_to_motivic_cohomology}: this can be applied since $\op{KGL}$ is the spectrum representing the cohomology theory associated with the Landweber exact formal group law $X+Y-\beta XY$ on $\bb Z[\beta,\beta^{-1}]$ by \cite{MR2496504}. In conclusion $f_0\op{KGL}$ satisfies the assumption \ref{sub:assumption_A1} with $\tau_0f_0\op{KGL} \simeq \cal K^M_\ast$.

  Now we deal with $\op K=\op{KQ}$. Consider the fiber sequence \[\widetilde f_1 \op{KQ} \ra \widetilde f_0\op{KQ} \ra \widetilde s_0 \op{KQ}.\]
  This time we have an isomorphism of commutative ring spectra $\widetilde s_0 \op{KQ} \simeq \widetilde{\op{H}} \bb Z$ by Theorem 5.2 of Chapter 7 of \cite{mwmotives}, while $\widetilde f_1 \op{KQ}\in \scr{SH}(K)_{\geq 1}$ by construction.
  As a consequence $\widetilde f_0 \op{KQ}$ satisfies the assumption \ref{sub:assumption_A1} with $\tau_0 \widetilde f_0 \op{KQ}\simeq\cal K^{MW}_\ast$. 


\section{The $E$-based Adams-Novikov spectral sequence} 
  \label{sec:the_E_based_MANSS}
  In this section we briefly recall the construction of the Adams-Novikov spectral sequence based on a homotopy commutative algebra in $\scr{M}$. This section works in any symmetric monoidal stable $\infty$-category $\scr{M}$.

  \subsection{Construction of the spectral sequence} 
    \label{sub:construction_of_the_spectral_sequence}
    Let $E$ be a homotopy commutative algebra in $\scr{M}$. We start by considering fiber sequence
    \begin{equation}
      \label{eqn:Ebar-S-E}
      \overline E \overset{\bar e}{\ra} \bf 1 \overset{e}{\ra}E
    \end{equation}
     where $\overline E:=\fib(e:\bf 1\ra E)$. We set the notation $\overline {E}^1=\overline E$ and $\overline E^0=\bf 1$. By induction, assuming we have already defined $\overline {E}^n$, we obtain a new fiber sequence by applying $-\wedge \overline E^n$ to the fiber sequence \eqref{eqn:Ebar-S-E}: we get the fiber sequence
    \begin{equation}
      \label{eqn:fund_fib_seq_E^n+1_E^n_EE^n}
    \xymatrix{
      \overline E \wedge \overline E^n \ar[r]^{\bar e\wedge \op{id}} &  \bf 1\wedge \overline E^n\ar[r]^{e \wedge \op{id}} &  E\wedge \overline E^{n}\\
    }.
    \end{equation}
    We set $\overline E^{n+1}:=\overline E \wedge \overline E^n$ and as well $W_n:=E\wedge \overline E^{n}$. Furthermore we name the maps $\bar e^{n+1}:=\bar e\wedge \op{id}_{\overline E^n}$ and $e^{n+1}:=e\wedge \op{id}_{\overline E^n}$. In this way we have produced a tower $\{\overline E^n\}_{n \in \bb N}$ over $\bf 1$ fitting in the following diagram
    \begin{equation}
      \label{eqn:Adams_tower_over_S}
      \xymatrix{
        \bf 1=\overline E^0 \ar[d]_e & \overline E^1 \ar[d]_{e^2} \ar[l]_{\bar e^1} & \overline E^2 \ar[d]_{e^3} \ar[l]_{\bar e^2} & \cdots \ar[l]_{\bar e^3} \\
        E = W_0 \ar@{-->}[ur] & E\wedge \overline E^1 = W_1 \ar@{-->}[ur] & E\wedge\overline E^2=W_2 \ar@{-->}[ur] & \cdots\\
      }
    \end{equation}
    where each dashed arrow is pictured to remind that the triangle it bounds is a fiber sequence.
    Given any object $X$ we can smash every part of the previous construction with $X$ and get a tower $\{X\wedge\overline E^n\}_{n\in \bb N}$ over $X$ and actually a whole diagram similar to \eqref{eqn:Adams_tower_over_S}.

      \subsubsection{}
      We could use an exact couple coming from \eqref{eqn:Adams_tower_over_S} to construct an Adams spectral sequence, but for our purposes it is more helpful to consider the tower under $\bf 1$ of cofibers induced by \eqref{eqn:Adams_tower_over_S}.
      With this aim in mind we proceed. The techniques of \ref{sub:towers} allow us to upgrade \eqref{eqn:Adams_tower_over_S} to a diagram $N(\bf N^{op}) \ra \scr{M}$, and thus to a fiber sequence of towers. 
      \begin{equation}
      \label{eqn:fib_seq_of_Ad_towers}\xymatrix{\overline E^\bullet \ar[r]^{e^{\bullet}} & \bf 1 \ar[r] & \overline E_{\bullet -1}}
      \end{equation}
      We visualize it as a commutative ladder of fiber sequences in $\scr{M}$:
    \begin{equation}
      \label{eqn:fund_comp_of_fib_seq_E^n_to_En}
      \xymatrix{
      \overline E^{n+1} \ar[r]^(.55){e^{n+1}} \ar[d]_{\bar e^{n+1}} & \bf 1 \ar[r] \ar@{=}[d] & \overline E_{n} \ar[d]^{f_n}\\
      \overline E^{n}  \ar[r]^(.55){e^{n}} & \bf 1 \ar[r]  & \overline E_{n-1},\\
      }
    \end{equation}
    where $e^n$ denotes an $n$-fold composition $\bar e^1\circ\cdots \circ \bar e^n$. Note that implicitly we have $\overline E_{-1} \simeq 0$ and $\overline E_0\simeq E$. Moreover we have equivalences $W_n=\cofib(\bar e^{n+1})\simeq \fib({f_n})$ in $\scr{M}$, and in particular fiber sequences
      \begin{equation}
        \label{eqn:fund_fib_seq_Wn_En_En-1}
        W_n \overset{l_n}{\ra} \overline E_n \overset{f_n}{\ra} \overline E_{n-1} \overset{\partial_n}{\ra} \Sigma^1 W_n.
      \end{equation}
     We thus get a new diagram
    \begin{equation}
      \label{eqn:Adams_tower_under_S}
      \xymatrix{
        \cdots  & W_3 \ar[d]^{l_3} &  W_2 \ar[d]^{l_2} & W_1 \ar[d]^{l_1} & W_0 \ar@{=}[d]^{l_0} & \\
        \cdots \ar[r] & \overline E_3 \ar[r]_{f_3} \ar@{-->}[ul] & \overline E_2 \ar[r]_{f_2} \ar@{-->}[ul] & \overline E_1 \ar[r]_{f_1} \ar@{-->}[ul]  & \overline E_0 \ar[r]_{f_0} \ar@{-->}[ul]& 0 \ar@{-->}[ul].\\
      }
    \end{equation}
    where again the dashed arrows are pictured to remind us that the triangles they bound are fiber sequences. Note that now the maps $f_n$ form a tower under $\bf 1$.

    As we did above, given any object $X$ we can build similar diagrams by applying $X\wedge-$ to \eqref{eqn:Adams_tower_under_S}. We obtain the following:
    \begin{equation}
      \label{eqn:Adams_tower_under_X}
      \xymatrix{
        \cdots  & X \wedge W_3 \ar[d]^{l_3} & X\wedge W_2 \ar[d]^{l_2} & X \wedge W_1 \ar[d]^{l_1} & X\wedge W_0 \ar[d]_{\simeq}^{l_0} & \\
        \cdots \ar[r] & X\wedge\overline E_3 \ar[r]_{f_3} \ar@{-->}[ul] & X \wedge\overline E_2 \ar[r]_{f_2} \ar@{-->}[ul]_{\partial_3} &X \wedge\overline E_1 \ar[r]_{f_1} \ar@{-->}[ul]_{\partial_2}  & X\wedge\overline E_0 \ar[r]_{f_0} \ar@{-->}[ul]_{\partial_1}& 0 \ar@{-->}[ul]_{\partial_0=0}.\\
      }
    \end{equation}
    Here we abuse the notation and keep denoting the maps involved in \eqref{eqn:Adams_tower_under_X} with the same names used above in \eqref{eqn:Adams_tower_under_S}.

    \begin{defin}
    \label{defin:std_adams_tower}
      The tower under $X$
      \begin{equation}
        \label{eqn:std_adams_tower}
        \cdots \overset{f_{n+1}}{\ra} X\wedge \overline E_n \overset{f_n}{\ra} \cdots\overset{f_2}{\ra} X \wedge\overline E_1 \overset{f_1}{\ra} X \wedge\overline E_0 \ra \ast
      \end{equation}
      is called the \emph{Standard $E$-Adams tower}. The \emph{$E$-nilpotent completion} of $X$ is the object 
      \[ X^\wedge_E:=\lim_{N(\bf N^{op})} (X\wedge \overline E_n) \in \scr{M}.\]
      The natural map $\alpha_E(X): X \ra X^\wedge_E$ is called the $E$-nilpotent completion map of $X$.
    \end{defin}

    \subsubsection{}
    For every object $Y$ we can apply the functor $[Y,-]$ to \eqref{eqn:Adams_tower_under_X} and get an exact couple 
    \begin{equation}
      \label{eqn:ANSS_second_ex_couple}
      \xymatrix{
        [\Sigma^{\bullet}Y,X\wedge\overline E_\bullet] \ar[rr]^{j} & & [\Sigma^{\bullet}Y,X\wedge\overline E_\bullet] \ar[dl]^k\\
         & [\Sigma^{\bullet}Y, X\wedge W_\bullet]. \ar[ul]^{i} & \\
        }
    \end{equation}
    Here the map 
    \[j: [\Sigma^{p}Y,X\wedge\overline E_n] \ra [\Sigma^{p}Y,X\wedge\overline E_{n-1}]\]
    is the natural map induced by $f_n$ and has bi-degree $(0,-1)$; the map
    \[k: [\Sigma^{p}Y,X\wedge\overline E_n]  \ra [\Sigma^{p-1}Y, X\wedge W_{n+1}]\]
    is the natural map induced by the dashed arrow $\partial_{n+1}$ and has bi-degree $(-1,1)$. Finally the map
    \[i: [\Sigma^{p}Y, X\wedge W_n] \ra [\Sigma^{p-1}Y,X\wedge\overline E_{n}] \]
    is the map induced by $l_n$ and has bi-degree $(0,0)$.

    \subsubsection{}
    The spectral sequence obtained form the exact couple \eqref{eqn:ANSS_second_ex_couple} is the \emph{$E$-based Adams-Novikov spectral sequence}. Note that this is an example of the general procedure described in IX.4 of \cite{MR0365573} for associating the so called \emph{Homotopy Spectral Sequence} to a tower of fibrations over a given space. In our specific example the tower we used is $\{X\wedge\overline E_n, f_n\}$. 

\section{Nilpotent Resolutions} 
  \label{sec:nilpotent_resolutions}
  In this section we introduce $E$-nilpotent (resp. strongly $\cal R$-nilpotent) resolutions associated with a homotopy commutative algebra $E$ in $\scr{M}$ (resp. a commutative algebra $\cal R$ in $\scr{M}^{\heartsuit}$). This takes place in \ref{sub:E-nilpotent-resol} (resp. \ref{sub:R_nilpotent_resolutions}). We use these constructions to describe a universal property for the Adams tower after passing to pro-objects, obtaining thus a more ductile construction of $E$-nilpotent completions. As an application we obtain an explicit description (see Theorem \ref{thm:convergence_mod_stuff}) for the $E$-nilpotent completion of a $k$-connected object. Throughout this section $\scr{M}$ is a presentably symmetric monoidal stable $\infty$-category, endowed with a left-complete accessible multiplicative $t$-structure.
  
  \subsection{$E$-Nilpotent resolutions} 
    \label{sub:E-nilpotent-resol}
    \begin{defin}
    \label{defin:E-nilpotence}
    Let $E$ be a homotopy commutative algebra in $\scr{M}$. We define the $\infty$-category of \emph{$E$-nilpotent} objects as the smallest full sub-$\infty$-category $\mathrm{Nilp}(E) \subseteq \scr{M}$ satisfying the following properties:
    \begin{enumerate}
      \item $E \in \mathrm{Nilp}(E)$;
      \item given any $X\in \scr{M}$ and any $F \in \mathrm{Nilp}(E)$      then $X\wedge  F \in \mathrm{Nilp}(E)$;
      \item $\mathrm{Nilp}(E)$ has the 2-out-of-3 property on fiber sequences, i.e. given a fiber sequence
        \[X\ra Y\ra Z\]
        in $\scr{M}$ where any two of the three objects $X,Y,Z$ are in $\mathrm{Nilp}(E)$, then the third is in $\mathrm{Nilp}(E)$ as well;
      \item $\mathrm{Nilp}(E)$ is closed under retracts;
    \end{enumerate}
     \end{defin}

    \begin{rmk}
      \label{rmk:mods_over_E-nilp_ringsp_are_E-nilp}
      If $R$ is a homotopy algebra and $M$ is a homotopy $R$-module then the action map $R\wedge  M \ra M$ in split by the unit. So if $R$ is in $\mathrm{Nilp}(E)$ then $M$ is $E$-nilpotent too. 
    \end{rmk}

    \begin{lemma}
      \label{lemma:E-nilp_implies_E-loc}
        If $E$ is a homotopy commutative algebra and $X$ is any $E$-nilpotent object, then $X$ is $E$-local  
    \end{lemma}

    \begin{proof}
      The proof goes exactly as in \cite[Lemma 3.8]{MR551009}. We filter $\rm{Nilp}(E)$ by inductively constructed subcategories $C_i$. $C_0$ is defined as the full sub-$\infty$-category of $\scr{M}$ whose objects are equivalent to $E\wedge  X$ for some $X \in \scr{M}$. If $i\geq 1$ we set $C_i$ to be the full subcategory of $\scr{M}$ of those objects that are equivalent to a retract of an object in $C_{i-1}$ or an extension of objects in $C_{i-1}$. It is formal to check that the union of the $C_i$'s coincides with $\rm{Nilp}(E)$.
      Indeed, thanks to \ref{rmk:mods_over_E-nilp_ringsp_are_E-nilp} we have that $C_0\subset \rm{Nilp}(E)$, and since $E$-nilpotent objects are closed under retractions and extensions, we get by induction that each of the $C_n$ is contained in $\rm{Nilp}(E)$. Now the $C_n$'s form an increasing sequence of subcategories of $\rm{Nilp}(E)$ and we need to check that their union, which we denote by $C$, is the whole $\rm{Nilp}(E)$. 
      However this is clear: by construction $C$ satisfies all the four axioms of \ref{defin:E-nilpotence} so we must have $C\supseteq \rm{Nilp}(E)$, and so $\rm{Nilp}(E) = C$. For proving the $E$-locality: $E$-modules are $E$-local, so $C_0\subseteq \rm{Loc}(E)$; since $E$-local objects are closed under extensions and retractions $C_i\subseteq \rm{Loc}(E)$, and hence $\rm{Nilp}(E) = \cup_i C_i\subseteq \rm{Loc}(E)$.
     \end{proof}
    
    \begin{defin}
      Let $E$ be a homotopy commutative algebra. An object $X$ is called $E$-pre-nilpotent if $X_E$ is $E$-nilpotent.
    \end{defin}
    \begin{prop}
      \label{prop:pre-nilp_sphere}
      Let $E$ be a homotopy commutative algebra. Then the following are equivalent:
      \begin{enumerate}[label=P.\arabic*]
        \item \label{prenilp_a} $\bf 1$ is $E$-pre-nilpotent, i.e. $\bf 1_E$ is $E$-nilpotent;
        \item \label{prenilp_b} for every object $X$, $\bf 1_E\wedge  X$ is $E$-nilpotent;
        \item \label{prenilp_c} every object $X$ is $E$-pre-nilpotent, i.e. $X_E$ is $E$-nilpotent for every $X$;
        \item \label{prenilp_d} $\rm{Nilp}(E)=\rm{Loc}(E)$.
      \end{enumerate}
      Moreover the following are equivalent:
      \begin{enumerate}[label=S.\arabic*]
        \item \label{prenilp_e} for every $X\in \scr{M}$, the map $\lambda_E(\bf 1)\wedge  \op{id}_X: X \ra \bf 1_E\wedge  X$ is an $E$-localization of $X$;
        \item \label{prenilp_f} the multiplication map of the $E$-local sphere $\bf 1_E \wedge  \bf 1_E \ra \bf 1_E$ is an equivalence and the natural inequality $\bc E \leq \bc {\bf 1_E}$ is an equality.
        \end{enumerate} 
      In addition the statement \ref{prenilp_b} implies \ref{prenilp_e}.
      Furthermore if $E$ has a multiplication map $E\wedge  E\ra E$ which is an equivalence, then the unit $e:\bf 1 \ra E$ is a localization map $\lambda_E(\bf 1)$, and condition \ref{prenilp_a} holds for $E$.
    \end{prop}
    \begin{proof}
       We start by observing that a localization map $\lambda_E(X)$ factors as
      \begin{equation}
      \label{eqn:E-loc_vs_smash_with_loc_sphere}
      \xymatrix{
      X \ar[r]^{\lambda_E(X)} \ar[d]_{\lambda_E(A)\wedge  \op{id}_X} & X_E\\
      \bf 1_E\wedge  X. \ar[ur]_{\tilde\lambda(X)} & \\
      }
      \end{equation}
      Since all the maps in the diagram are $E$-equivalences, $\tilde \lambda(X)$ is an equivalence if and only if $\bf 1_E\wedge  X$ is $E$-local. 

      $E$-nilpotent objects are closed under smashing with arbitrary objects, so that $\bf 1_E$ is $E$-nilpotent if and only if for every $X\in \scr{M}$, $\bf 1_E\wedge  X$ is $E$-nilpotent too (i.e. \ref{prenilp_a} is equivalent to \ref{prenilp_b}). In view of \eqref{eqn:E-loc_vs_smash_with_loc_sphere}, if $\bf 1_E\wedge  X$ is $E$-nilpotent then it is $E$-local and hence $\tilde\lambda(X)$ is an equivalence (i.e. \ref{prenilp_b} implies \ref{prenilp_c}). Clearly \ref{prenilp_c} implies \ref{prenilp_a} and \ref{prenilp_c} is equivalent to \ref{prenilp_d}. 

      Using \eqref{eqn:E-loc_vs_smash_with_loc_sphere} we immediately deduce that \ref{prenilp_b}$\Rightarrow$\ref{prenilp_e}.

      By applying \ref{prenilp_e} to $X=\bf 1_E$, in view of \eqref{eqn:E-loc_vs_smash_with_loc_sphere}, we deduce that the multiplication map $\tilde\lambda(\bf 1_E)$ of $\bf 1_E$ is an equivalence. On the other hand, by smashing the fiber sequence
      \[{}_E\bf 1 \ra \bf 1 \ra \bf 1_E\]
      with an object $X$, we deduce that if $X$ is $\bf 1_E$-acyclic, then it is also $E$-acyclic. This means that $\bc E\leq \bc{\bf 1_E}$. The reverse equality is immediate from \ref{prenilp_e}, so \ref{prenilp_e} implies \ref{prenilp_f}. Assume now \ref{prenilp_f}. Since the multiplication of $\bf 1_E$ is an equivalence, for every $X\in \scr{M}$ the map $\lambda_{E}(A)\wedge  \op{id}_X: X \ra \bf 1_E\wedge  X$ is an $\bf 1_E$-localization of $X$; however $\bc{\bf 1_E}=\bc E$ so that \ref{prenilp_f} implies \ref{prenilp_e}.

      If $E$ is a homotopy commutative algebra with the property that the multiplication $E\wedge  E \ra E$ is an equivalence, then the unit $e: \bf 1 \ra  E$ is an $E$-equivalence; since $E$ is $E$-nilpotent, and thus $E$-local, we conclude.
    \end{proof}

    \begin{defin}
    \label{defin:smashing_localization}
      We say that an object $E$ induces a smashing localization if the map $\tilde \lambda_E(X): \bf 1_E \wedge  X \ra X_E $ of \eqref{eqn:E-loc_vs_smash_with_loc_sphere} is an equivalence in $\scr{M}$.
    \end{defin}

    \begin{ex}
    \label{ex:nilp_loc}
      The objects $\bf 1^{\cal B}$ appearing in Proposition \ref{prop:propert_of_trivializ} give smashing localizations of $\scr{M}$.

      Let $S$ be the spectrum of a field, $\scr{M}=\scr{SH}(S)$, and $E=\op H\bb Q$ be the spectrum representing Voevodsky's motivic cohomology with rational coefficients. Combining Proposition 14.1.6 with Corollary 16.1.7 of \cite{cisdeg_triang_mixed}, we deduce that the multiplication map of $\op H \bb Q$ is an equivalence. In particular the localization at $\op H\bb Q$ is smashing. 
    \end{ex}

    \begin{lemma}
    \label{lemma:smashing_implies_Loc_is_Mod}
      Let $E$ be an object inducing a smashing localization, and let $\bf 1_E$ be an $E$-localization of $\bf 1$. Then the forgetful functor $U_{\bf 1_E}: \scr{Mod}_{\bf 1_E}(S) \ra \scr{M}$ factors through an equivalence $U_{\bf 1_E}:\scr{Mod}_{\bf 1_E}(S)\ra \rm{Loc}_E$.
    \end{lemma}
    \begin{proof}
      Since $\bf 1_E$-modules are $\bf 1_E$-local and $\rm{Loc}_{\bf 1_E}=\rm{Loc}_E$ by \ref{prop:pre-nilp_sphere}\eqref{prenilp_f} we have the desired factorization of $U$ through the inclusion $\rm{Loc}_{E}\subseteq \scr{M}$. Note that every $E$-local object is in the underlying object of a $\bf 1_E$-module, since by definition $\lambda_E(X): X \ra \bf 1_E\wedge X$ is an equivalence. In order to show that $U$ is fully-faithful, we just need to see that for every $\bf 1_E$-module $X$, the natural "action map" $\alpha: F_{\bf 1_E}U_{\bf 1_E}(X) \ra X$ is an equivalence. 
      Indeed we have a commutative diagram of spaces
      \[\xymatrix{
      \Map_{\bf 1_E}(X,Y)\ar[rr] \ar[dr] & & \Map(U_{\bf 1_E}(X),U_{\bf 1_E}(Y))\\
        & \Map_{\bf 1_E}(F_{\bf 1_E} U_{\bf 1_E} (X),Y)\ar[ur]^{\simeq} & \\ 
      }\]
      where $\Map_{\bf 1_E}(-,-)$ denotes the mapping space in $\bf 1_E$-modules.
      Checking that $\alpha$ is an equivalence can be done after forgetting to $\scr{M}$, where $\alpha$ has a right inverse induced by the unit $e\wedge \op{id}: U_{\bf 1_E}(X) \ra \bf 1_{E}\wedge U_{\bf 1_E}(X)$. However $e\wedge \op{id}$ is an equivalence, since $E$ is smashing and $U_{\bf 1_E}(S)$ is $E$-local, so $\alpha$ is an equivalence too.
    \end{proof}
    \subsubsection{}
      \label{subs:Relation between localization and nilpotent completion}
      Let $E$ be a homotopy commutative algebra. We wish to point out how Definition \ref{defin:std_adams_tower} and Lemma \ref{lemma:E-nilp_implies_E-loc} imply that, for every $X\in \scr{M}$, the $E$-nilpotent completion $X^{\wedge}_E$ is $E$-local. Indeed
      \[X^\wedge_E=\lim \big(  \cdots \overset{f_{n+1}}{\ra} X\wedge  \overline E_n \overset{f_n}{\ra} \cdots \ra  X\wedge  \overline E_0 \ra 0 \big), \]
      and each of the maps in the tower sits in a fiber sequence
        \[
           E \wedge  \overline E^n \wedge  X \ra \overline E_{n}\wedge  X \overset{f_n}{\ra} \overline E_{n-1}\wedge  X
        \]
      that we have deduced from \eqref{eqn:fund_fib_seq_Wn_En_En-1}. As a consequence, by induction, each of the terms in the tower is $E$-nilpotent, hence $E$-local, and thus $X^\wedge_E$ is $E$-local too. In particular the natural map $\alpha_E(X) : X \ra X^{\wedge}_E$ factors as 
      \begin{equation}
      \label{eqn:alphabetagamma}
        \xymatrix{
        X \ar[rr]^{\alpha_E(X)} \ar[dr]_{\lambda_E(X)} & & X^{\wedge}_E \\
        & X_E \ar[ur]_{\beta_E(X)} & \\
        }.
      \end{equation}
      It follows that $\alpha_E(X)$ is an $E$-equivalence if and only if the induced map $\beta_E(X)$ is an equivalence.

      \subsubsection{}
      We wish to point out another fact. On one hand, if $X\ra Y$ is an $E$-equivalence, then it induces an equivalence of the Standard $E$-Adams Towers (\ref{defin:std_adams_tower}) associated to $X$ and $Y$, so that the natural map induced on homotopy inverse limits $X^\wedge_E \ra Y^\wedge_E$ is an equivalence.
      On the other hand the composition of $\alpha_E(X)$ with the projection to the $0$-th term of the tower
      \[  X \ra X^\wedge_E \ra X\wedge \overline E_0=X\wedge  E\] is identified with $\op{id}_X\wedge  e$, where $e: \bf 1 \ra  E$ is the unit of the algebra $E$. Thus, after smashing with $E$, 
      the map $\alpha_E(X)\wedge  E: X\wedge  E \ra X^\wedge_E\wedge  E$ has a splitting which is functorial in $X$. So if $f: X\ra Y$ is a map inducing an equivalence on $E$-nilpotent completions $X^\wedge_E \ra Y^\wedge_E$, then $f$ is an $E$-equivalence. We conclude that $\alpha_E(X)$ is an $E$-equivalence if and only if the induced map 
      \[\alpha_E(X)^\wedge_E: X^\wedge_E\ra (X^\wedge_E)^\wedge_E\] is an equivalence.  
        \begin{defin}
      \label{defin:E-nilp_resol}
      For an object $X\in \scr{M}$, an \emph{$E$-nilpotent resolution of $X$} is a tower of objects under $X$
      \[X \ra \cdots  \ra X_n \ra X_{n-1} \ra \cdots \ra X_0 \]
      satisfying the following two properties:
      \begin{enumerate}
        \item $X_n \in \mathrm{Nilp}(E)$ for every $n \in \bb N$;
        \item for any $Y\in \mathrm{Nilp}(E)$ the map of pro-objects $\{X\}\ra \{X_\bullet\}$ defined by the tower induces an equivalence 
        \[\Map_{\rm{Pro(\scr{M})}}(\{X_\bullet\},\{Y\}) \ra \Map_{\rm{Pro(\scr{M})}}(\{X\},\{Y\}),\]
        where $X$ and $Y$ are considered as constant pro-objects.
      \end{enumerate}
    \end{defin}
    \subsubsection{}
    Recall that for pro-objects $X_\bullet,Y_\bullet$ in an $\infty$-category $\scr C$
    \[\Map_{\rm{Pro}(\scr M)}(X_\bullet,Y_\bullet)\simeq \rm{lim}_m\colim_n\Map_{\scr{M}}(X_n,Y_m),\]
    (c.f. Lemma \ref{lemma:map_of_probjects}).
    In our situation, by applying homotopy groups (of the geometric realization) to the previous formula, we get that a tower $X \ra X_\bullet$ of nilpotent objects under a given $X \in \scr{M}$ is an $E$-nilpotent resolution if and only if for every $E$-nilpotent object $Y \in \scr{M}$ the induced map
    \[\pi_i\colim_n \Map_{\scr{M}}(X_n,Y) \ra \pi_i\Map_{\scr{M}}(X,Y)\]
    is an isomorphism. Since taking homotopy groups commutes with filtered colimits, and since $E$-nilpotent objects are closed under shifts, this is equivalent to asking that for every $E$-nilpotent object $Y\in \scr{M}$ the natural map
    \[\colim_n [X_n,Y]  \ra [X,Y] \]
    is an isomorphism. The definition we have given is thus compatible with the classical definition for the stable homotopy category.

    \begin{prop}
      \label{prop:E-nilp_resol_are_unique_and_com_the_nilp_comp}
      Let $X$ be any object of $\scr{M}$. Then:
      \begin{enumerate}
        \item \label{ENR1} the standard Adams tower $ \overline{E}_\bullet \wedge  X$ is an $E$-nilpotent resolution of $X$;
        \item \label{ENR2} the pro-object under $X$ associated with an $E$-nilpotent resolution of $X$ is unique up to a contractible space of choices;
        \item \label{ENR3} if $X \ra X_\bullet$ is a $E$-nilpotent resolution of $X$, there is a natural equivalence $\rm{lim}_{N(\bf N^{op})} X_\bullet \simeq X^{\wedge}_{E}$ in $\scr{M}_{X/}$.
      \end{enumerate}
    \end{prop}

    \begin{proof}
     We start with \eqref{ENR1}. As we already observed in Remark \ref{subs:Relation between localization and nilpotent completion}, the terms $\overline E_n\wedge  X$ of the tower are $E$-nilpotent. Let $Y$ be any object of $\scr{M}$. By smashing the fiber sequence of towers \eqref{eqn:fib_seq_of_Ad_towers} with $X$ we get a commutative ladder of long exact sequences
     \begin{equation}
      \xymatrix{
       \cdots \ar[r] & [\overline E_n\wedge  X,Y] \ar[r] & [X,Y] \ar[r] & [\overline E^{n+1}\wedge  X,Y] \ar[r] & \cdots\\
       \cdots \ar[r] & [\overline E_{n-1}\wedge  X,Y] \ar[r]\ar[u]^{f_n^\ast} & [X,Y] \ar[r] \ar@{=}[u] & [\overline E^{n}\wedge  X,Y] \ar[r]\ar[u]^{\bar e_{n+1}^\ast} & \cdots .\\
      }
     \end{equation}
      We deduce that we only need to show that $\varinjlim _n [\overline E^n\wedge  X,Y]=0$ for every $E$-nilpotent object $Y$. We will proceed by induction on the family of subcategories $C_i$ that we used in the proof of \ref{lemma:E-nilp_implies_E-loc}. Assume thus that $Y\in C_0$, i.e. that $Y\simeq E\wedge  Z$ for some $Z\in \scr{M}$: we will show that the transition maps in the colimit vanish, hence the colimit vanishes too. For this we look at the fiber sequence \eqref{eqn:fund_fib_seq_E^n+1_E^n_EE^n}: it gives a long exact sequence 
     \begin{equation}
       \cdots \ra [E\wedge  \overline E^n\wedge  X,E\wedge  Z] \overset{(e\wedge  \rm{id}_X)^\ast}{\rra} [\overline E^n\wedge  X, E\wedge  Z] \overset{(\overline e^{n+1}\wedge  \rm{id}_X)^\ast}{\rra}[\overline E^{n+1} \wedge  X, E\wedge  Z] \ra \cdots
     \end{equation}
     where the maps $(e\wedge  \rm{id_X})^\ast$ are surjective since $E\wedge  Z$ is an $E$-module. The transition maps in the direct limit are thus $0$. Now observe that the property that $\varinjlim _n [\overline E^n\wedge X,Y]=0$ is stable in the $Y$ variable under retracts and extensions: this implies that if every $Y \in C_{i-1}$ satisfies this property then also every $Y \in C_i$ does as well. Since the union of the $C_i$'s exhausts $\rm{Nilp}(E)$ the first point is done.

    \eqref{ENR2}. Existence is \eqref{ENR1}, and for uniqueness we argue as follows. Let $X_\bullet \overset{p}{\leftarrow} X \overset{q}{\rightarrow} Y_\bullet$ be $E$-nilpotent resolutions. Then we have a chain of equivalences 
      \[
      \xymatrix{
      \Map_{\rm{Pro}(\scr{M})}(\{X_\bullet\},\{Y_\bullet\}) \ar[r]^(.45){\simeq} &  \lim_{\op N(\bf N^{op})} \Map_{\rm{Pro}(\scr{M})}(\{X_\bullet\},\{Y_n\})\ar[d]^{\circ p}  & \\
      & \lim_{\op N(\bf N^{op})} \Map_{\rm{Pro}(\scr{M})}(\{X\},\{Y_n\}) &  \Map_{\rm{Pro}(\scr{M})}(\{X\},\{Y_\bullet\})\ar[l]_(.45){\simeq},\\
      }\]
      where the first and last the obvious maps, while the second map is induced by the the projection $p$. This is enough to conclude that the diagram of pro-objects $\{X_\bullet\} \overset{p}{\leftarrow} \{X\} \overset{q}{\rightarrow} \{Y_\bullet\}$ can be filled essentially in a unique way to a $2$-simplex
      \[\xymatrix{
      & \{X\}\ar[dl]_p  \ar[dr]^q& \\
      \{X_\bullet\}\ar[rr]^f & & \{Y_\bullet\},
      }\]
      as the next lemma shows. The same argument with $X$ and $Y$ interchanged implies that any choice for $f$ must be an equivalence.
      Point \eqref{ENR3} follows by combining \eqref{ENR1} and \eqref{ENR2} with the observation that the operation of taking inverse limits factors through pro-objects.

    \end{proof}
    \begin{lemma}
    \label{lemma:com_tri_inf_cat}
      Let $X \overset{p}{\leftarrow} Z\overset{q}{\rightarrow} Y$ be a diagram in an $\infty$-category $\scr C$, and assume that composition with $p$ induces an equivalence $\Map(X,Y) \ra \Map(Z,Y)$. Then there exist a unique $2$-simplex of $\scr C$ extending the horn $(p,q)$ up to a contractible space of choices.
    \end{lemma}
    \begin{proof}
      Consider the over-category $\scr C_{p/}$: the restriction along the source and target of $p$ induce respectively functors
      \[\xymatrix{
      \scr C_{Z/} & \scr C_{p/}\ar[l]_{\sigma}\ar[r]^{\tau}_{\simeq}  & \scr C_{X/},\\
      }\]
      where $\sigma$ is a left fibration (Proposition 2.1.2.1 of \cite{HTT}), and $\tau$ is a trivial Kan fibration (Proposition 2.1.2.5 of \cite{HTT}). Moreover the three categories are themselves total spaces of left fibrations over $\scr C$, and both $\sigma$ and $\tau$ commute with the projection to $\scr C$. When we take fibers over $Y \in \scr C$ we get 
      \[
      \xymatrix{
      \Map^L(Z,Y) & (\scr C_{p/})_{Y} \ar[l]_(0.4){\sigma_{Y}} \ar[r]^(0.4){\tau_{Y}}_(0.4){\simeq} &\Map^L(X,Y),\\
      }
      \]
      where $\Map^L(-,-)$ denotes the left mapping space. Here $\sigma_Y$ is a left fibration over a Kan complex, and thus it is a Kan fibration by Lemma 2.1.3.3. of \cite{HTT}. Similarly $\tau_{Y}$ is again trivial Kan fibration.  Our assumption then implies that $\sigma_{Y}$ is a trivial Kan fibration. By definition of left mapping space, the zero simplexes of $\Map^L(Z,Y)$ are exactly the arrows $Z \ra Y$ of $\scr C$, while the zero simplexes of $(\scr C_{p/})_{Y}$ are commutative triangles
      \[\xymatrix{
       & Z \ar[dr]_{q}\ar[dl]^p & \\
        X \ar[rr] & & Y.}\\\]
        Our argument then says that there exist only a contractible spaces of arrows $f: X \ra Y$ making the triangle commutative.
    \end{proof}

  \subsection{Strongly $\cal R$-nilpotent resolutions} 
    \label{sub:R_nilpotent_resolutions}
    \begin{lemma}
      \label{lemma:props_of_solid_rings}
      Let $\cal R$ be a commutative algebra in $\scr{M}^{\heartsuit}$ with the property that the multiplication map of $\mu_{\cal R} : \cal R\otimes^{\heartsuit}\cal R \ra \cal R$ of $\cal R$ is an isomorphism. Then:
      \begin{enumerate}
        \item \label{PSR2} For every $\cal R$-module $\cal M$ in $\scr{M}^{\heartsuit}$, the action map $\cal R\otimes^{\heartsuit} \cal M\ra \cal M$ is an isomorphism. In particular an $\tau_0(\bf 1)$-module $\cal M$ has at most one  $\cal R$-module structure.
        \item \label{PSR3} Every map $\phi: \cal M\ra \cal N$ in $\scr{M}^{\heartsuit}$ where $\cal M$ and $\cal N$ are $\cal R$-modules is $\cal R$-linear. In particular the category of $\cal R$-modules is a full subcategory of $\scr{M}^{\heartsuit}$.
      \end{enumerate}
    \end{lemma}
    \begin{proof}
      For \eqref{PSR2} observe that, given a $\cal R$-module $\cal M$ in $\scr{M}^{\heartsuit}$, we have a co-equalizer diagram defining the monoidal product on the category of $\cal R$-modules
      \[
      \xymatrix{
      \cal R \otimes^{\heartsuit} \cal R \otimes^{\heartsuit} \cal M \ar@<.4ex>[r]^{1\otimes a} \ar@<-.4ex>[r]_{\mu_{\cal R}\otimes 1} & \cal R\otimes^{\heartsuit} \cal M \ar[r]^q & \cal R\otimes^{\heartsuit}_{\cal R}\cal M. \\
      }
      \]
      Moreover the action map $a: \cal R \otimes^{\heartsuit} \cal M \ra \cal M$ induces an isomorphism $\bar a: \cal R \otimes^{\heartsuit}_{\cal R} \cal M\ra \cal M$. The map $\cal R\otimes^{\heartsuit} \cal M \ra \cal R\otimes^{\heartsuit}\cal R\otimes^{\heartsuit} \cal M$ defined by $r\otimes m \mapsto r\otimes 1 \otimes m$ is an inverse of both $\mu_{\cal R}$ and $a$, so that $1\otimes a=\mu_{\cal R}\otimes 1$ and $q$ is as isomorphism. In particular every $\cal R$-module $\cal M$ is isomorphic to the free $\cal R$-module on the object $\cal M$. Point \eqref{PSR3} follows by combining the free-forget adjunction and point \eqref{PSR2}.
    \end{proof}

    \begin{lemma}
      \label{lem:AssA1_implies_solidity_of_pi_0}
      Assume $E$ is a homotopy commutative algebra satisfying assumption \ref{sub:assumption_A1}. Then $\cal R=\tau_0(E)\in \scr{M}^{\heartsuit}$ and its multiplication map $\cal R\otimes^\heartsuit\cal R \rightarrow \cal R$ is an isomorphism. 
    \end{lemma}

    \begin{proof}
       Surely $\tau_0(E)$ is a commutative algebra in $\scr{M}^{\heartsuit}$ by the multiplicative properties of the Postnikov tower \ref{subsub:trunc_of_hom_mod}. Let us start from the special case $J=\emptyset$, so that $\cal R=\tau_0 E\simeq \tau_0(\bf 1)/(f_1,\dots,f_r)$ for some $f_i : L_i \ra \bf 1$. In this case the map
       \[(f_1,\dots,f_r):\bigoplus_i L_i\ra \bf 1\] induces 
       \[(f_1,\dots,f_r)^{\heartsuit}:\bigoplus_i L_i\wedge  \tau_0 (\bf 1) \ra \tau_0(\bf 1)\] in $\scr{M}^{\heartsuit}$. Since $\tau_0 E\simeq \coker ((f_1,\dots,f_r)^{\heartsuit})$, Lemma \ref{lemma:ker_and_coker_of_multiplications} implies that $\cal R$ is solid.
      
       Let us now consider the case of $J\not = \emptyset$. Denote by $\cal C:=\tau_0(\bf 1)/\rm{Im}(f_1,\dots,f_r)^{\heartsuit}$ so that we already know that $m: \cal C\otimes^{\heartsuit}\cal C \ra \cal C$ is an isomorphism by the above argument, and that $\tau_0(E)\simeq \tau_0(\cal C[\cal J^{-1}])$ where $\cal J=\{g_j\}_{j \in J} $. 
       Now by \ref{cor:inv_elts_is_t_exact} the functor $[\cal J^{-1}]$ is right $t$-exact, so we have a natural equivalence
       \[\tau_0(\cal C[\cal J^{-1}])\otimes ^{\heartsuit} \tau_0(\cal C[\cal J^{-1}]) \simeq \tau_0((\cal C\otimes^{\heartsuit} \cal C)[\cal J^{-1}])\]
       under which the multiplication of $\tau_0(\cal C[\cal J^{-1}])$ is identified with the map $\tau_0(m[\cal J^{-1}])$, which is an equivalence by the previous part.
    \end{proof}

    \subsubsection{}
    	From the rest of the section we fix a commutative algebra $\cal R$ in $\scr{M}^{\heartsuit}$ with the property that its multiplication is an isomorphism. 

    \begin{defin}
      \label{defin:S-nilpotent_hom_module_and_spectra}
        We say that an object $\cal M$ of $\scr{M}^\heartsuit$ is \emph{strongly $\cal R$-nilpotent} if it has a finite filtration $\cal M = \cal M_0 \supseteq \cal M_1 \supseteq \cdots \supseteq \cal M_r$ in $\scr{M}^\heartsuit$ such that $\cal M_i/\cal M_{i+1}$ has an $\cal R$-module structure for every $i$. An object $X$ is called \emph{strongly $\cal R$-nilpotent} if for each $k\in \bb Z$ the homotopy object $\tau_k(X)$ is $\cal R$-nilpotent and there exist integer $s$ and $t$ such that $X\simeq P_s(X)\simeq P^t(X)$. We denote the full subcategory of strongly $\cal R$-nilpotent objects by $\mathrm{SNilp}(\cal R)$.
    \end{defin}

    \subsubsection{} 
    	\label{lemma:S-nilp_implies_E-nilp}
    	Observe that if $\cal M \in \scr{M}^{\heartsuit}$ is strongly $\cal R$-nilpotent then it is strongly $\cal R$-nilpotent as an object of $\scr{M}$; conversely any strongly $\cal R$-nilpotent object of $\scr{M}$ which is concentrated in degree $0$ for the $t$-structure is a strongly $\cal R$-nilpotent object of $\scr{M}^{\heartsuit}$. If an object $X$ in $\scr{M}$ is strongly $\cal R$-nilpotent then it is $\cal R$-nilpotent in the sense of Definition \ref{defin:E-nilpotence}. Indeed since $X$ is bounded in the $t$-structure of $\scr{M}$, $X$ is an iterated extension of finitely many of its homotopy objects $\Sigma^i\tau_i(X)$ which are, on their turn, iterated extensions of $\cal R$-modules. In particular if $E$ is a $(-1)$-connected homotopy commutative algebra in $\scr{M}$ and $\tau_0(E)\simeq \cal R$, then any strongly $\cal R$-nilpotent object $X$ is also $E$-nilpotent.

    \begin{defin}
      \label{defin:S-nilp_resol}
      Let $X$ be an object in $\scr{M}$. A \emph{strongly $\cal R$-nilpotent resolution of $X$} is a tower of objects under $X$
      \[X \ra \cdots  \ra X_n \ra X_{n-1} \ra \cdots \ra X_0 \]
      satisfying the following two properties:
      \begin{enumerate}
        \item $X_n \in \mathrm{SNilp}(\cal R)$ for every $n \in \bb N$;
        \item for any $Y\in \mathrm{SNilp}(\cal R)$ the induced map of pro-objects $\{X\}\ra \{X_\bullet\}$ defined by the tower induces an equivalence 
        \[\Map_{\rm{Pro}}(\{X_\bullet\},\{Y\}) \ra \Map_{\rm{Pro}}(\{X\},\{Y\}),\]
        where $X$ and $Y$ are considered as constant pro-objects.
      \end{enumerate}
    \end{defin}

     We now prepare for an analogue of Proposition \ref{prop:E-nilp_resol_are_unique_and_com_the_nilp_comp} for strongly $\cal R$-nilpotent resolutions.
    
    \begin{lemma}
      \label{lemma:S-nilp_in_sequences}
      If 
      \[\cal M\ra \cal N\ra \cal O\ra \cal P\ra \cal Q\]
      is an exact sequence in $\scr{M}^{\heartsuit}$ where $\cal M,\cal N,\cal P,\cal Q$ are strongly $\cal R$-nilpotent, then $\cal O$ is as well.
    \end{lemma}

    \begin{proof}
      By breaking up the exact sequence in shorter pieces, the statement follows by combining \ref{lemma:S-nilp_in_sequences-extensions} and  \ref{lemma:S-nilp_in_sequences-ker_and_coker}.      
    \end{proof}

    \begin{lemma}
      \label{lemma:S-nilp_in_sequences-extensions}
      If $0 \ra \cal M \ra \cal N \ra \cal O \ra 0$ is a short exact sequence in $\scr{M}^{\heartsuit}$ and both $\cal M$ and $\cal O$ are strongly $\cal R$-nilpotent, then $\cal N$ is too.
    \end{lemma}
    \begin{proof}
      A suitable filtration on $\cal N$ can be obtained by combining the filtration on $\cal M$ and the pre-image in $\cal N$ of the filtration on $\cal O$.
    \end{proof}

    \begin{lemma}
      \label{lemma:S-nilp_in_sequences-ker_and_coker}
      If $\phi: \cal C^0\ra \cal C^1$ is a map in $\scr{M}^{\heartsuit}$ and both $\cal C^0$ and $\cal C^1$ are strongly $\cal R$-nilpotent, then both $\cal H^0=\rm{Ker} \phi$ and $\cal H^1=\rm{Coker}\phi$ are strongly $\cal R$-nilpotent too.
    \end{lemma}
    \begin{proof}
      Up to increasing the length of the filtrations $\cal C^0_i$ and $\cal C^1_i$ we can assume that $\phi$ respects the filtrations. As a consequence $(\cal C^\bullet,\phi)$ is a filtered complex in $\scr{M}^{\heartsuit}$.
      With respect to this filtration we use the spectral sequence for filtered complexes. By assumption the terms in the $E_2$-page are $\cal R$-modules and the differentials are $\cal R$-linear by  \ref{lemma:props_of_solid_rings}. Hence we gain a finite filtration on the cohomology of $\cal C^\bullet$ with associated graded $\cal R$-module quotients.
    \end{proof}
    
    \begin{lemma}
    \label{lemma:mul_by_f_in_a_B_nilp_A-mod}
      Let $f:L\ra \tau_0(A)$ be a morphism in $\scr{M}$, where $L$ is a tif object. Denote by $\cdot f$ the induced map $\cal R\wedge L \ra \cal R$ in $\scr{M}^{\heartsuit}$ given by \ref{lemma:ker_and_coker_of_multiplications}, and let $\cal S:=\coker(\cdot f)$. Then $\cal S$ is a commutative $\cal R$-algebra in $\scr{M}^{\heartsuit}$ and its multiplication map $\mu_{\cal S}: \cal S\otimes^{\heartsuit}\cal S\ra \cal S$ is an isomorphism. Furthermore for every strongly $\cal R$-nilpotent object $\cal M$ of $\scr{M}^{\heartsuit}$ both kernel and cokernel of the multiplication by $f$ on $\cal M$ are strongly $\cal S$-nilpotent.
    \end{lemma}

    \begin{proof}
      First of all since $L\in \scr{M}_{\geq 0}$, $f$ induces a map $\tau_0(L) \ra \tau_0(A) \ra \cal R$, and $\tau_0(L)$ is a strongly dualizable object of $\scr{M}^{\heartsuit}$ so Lemma \ref{lemma:ker_and_coker_of_multiplications} fully applies. We deduce that $f$ extends uniquely to the map $\cdot f:\cal R\wedge L \ra \cal R$, the commutative algebra structure of $\cal R$ in $\scr{M}^{\heartsuit}$ descends uniquely to a commutative algebra structure on $\cal S$, such that the projection $\cal R\ra \cal S$ is an algebra map. Furthermore the multiplication induced on $\cal S$ is an isomorphism.

      Let $0=\cal M_n \subseteq \cal M_{n-1} \subseteq \cdots \subseteq \cal M_0=\cal M$ be a filtration of $\cal M$ by objects $\cal M_i$ of $\scr{M}^{\heartsuit}$ whose associated graded pieces are $\cal R$-modules, and denote by $r_f(\cal M): \cal M\wedge L \ra \cal M$ the induced right multiplication by $f$ on $\cal M$. If $n=1$ then $\cal M$ is an $\cal R$-module, $r_f(\cal M)$ is automatically $\cal R$-linear by Lemma \ref{lemma:props_of_solid_rings}, and thus kernel and co-kernel of $r_f(\cal M)$ are $\cal S$-modules by \ref{lemma:ker_and_coker_of_multiplications}. If $n>1$ one can proceed by induction. Indeed, since the multiplication map $r_f(\cal M): \cal M\wedge L \ra \cal M$ commutes with any map in $\scr{M}$ by \ref{subsub:right_and_left_mult}, it respects the filtration. As a consequence, for every integer $k\geq 1$, $r_f(\cal M)$ induces a map of short exact sequences
      \begin{equation}
        \xymatrix{
        0 \ar[r] & \cal M_{k-1} \ar[r] \ar[d]^{r_f(\cal M)} & \cal M_{k} \ar[r] \ar[d]^{r_f(\cal M)} & \cal M_k/\cal M_{k-1} \ar[r] \ar[d]^{r_f(\cal M)} & 0\\
        0 \ar[r] & \cal M_{k-1} \ar[r] & \cal M_k \ar[r] & \cal M_k/\cal M_{k-1} \ar[r] & 0.\\
        }
      \end{equation}
      Thus by combining the snake lemma together with \ref{lemma:S-nilp_in_sequences} and the inductive assumption we conclude.
    \end{proof}    
    

  \subsection{Relation between localizations and nilpotent completions}

  	\begin{notat}
      For the rest of the section we fix a homotopy commutative algebra $E$ in $\in \scr{M}_{\geq 0}$, and we assume that the induced multiplication on $\tau_0(E)$ is an isomorphism. We also denote by $\cal R$ the homotopy object $\tau_0 E$. 
    \end{notat}
    \begin{prop}
      \label{prop:PnEn_tower_is_an_S-nilp_resol}
      Let $X \in \scr{M}_{\geq k}$ be a for some integer $k$. Then the tower $X \ra P^\bullet(\overline{E}_\bullet\wedge X)$ is a strongly $\cal R$-nilpotent resolution of $X$.
    \end{prop}

    \begin{proof}
      We first need to check that $P^n(\overline{E}_n\wedge X) \in \mathrm{SNilp}(\cal R)$ for every $n\in \bb Z$. By the connectivity of $X$, each of the $P^n(\overline{E}_n\wedge X)$ in bounded in the $t$-structure, so we only need to check that $\tau_k(\overline{E}_n\wedge X)$ is strongly $\cal R$-nilpotent for every pair of integers $k,n$. Recall that for every $n$ in $\bb N$ we have a fiber sequence of the form
      \begin{equation}
        \label{eqn:the_main_fib_seq_in_ANSS}
        X \wedge \overline{E}^n\wedge E \ra X \wedge \overline{E}_n \overset{f_n}{\ra} X \wedge \overline{E}_{n-1}   
      \end{equation}
      obtained from \eqref{eqn:fund_fib_seq_Wn_En_En-1} by smashing with $X$. By construction, for $n=0$, such a fiber sequence is the cone sequence of the identity of $X\wedge E$. In particular for every $k\in \bb Z$, $\tau_k( X\wedge \overline E_0)$ is a $\tau_0E$-module and thus it is strongly $\cal R$-nilpotent ($\cal R = \tau_0 E$). We can now proceed by induction on $n$. Note that $\tau_k(X\wedge \overline E^n\wedge E)$ is an $\cal R$-module too and thus is strongly $\cal R$-nilpotent. This observation, once combined with the inductive assumption, allows to apply Lemma \ref{lemma:S-nilp_in_sequences} to the long exact sequence of homotopy objects associated to the fiber sequence \eqref{eqn:the_main_fib_seq_in_ANSS}. 

      As a second step we need to prove that for any $Y\in \mathrm{SNilp}(\cal R)$ the projection $p:X\ra P^{\leq \bullet}(\overline E_{\bullet}\wedge X)$ induces an equivalence
      \[\Map_{\rm{Pro}(\scr{M})}(\{P^{\leq \bullet}(\overline E_{\bullet}\wedge X)\},\{Y\}) \ra \Map_{\rm{Pro}(\scr{M})}(\{X\},\{Y\}).\]
      However the factorization through the Adams tower $X\ra \overline E_{\bullet}\wedge X$ induces a commutative diagram
      \begin{equation}
        \xymatrix{
         \colim_n \Map(P^n(\overline E_n\wedge X),Y) \ar[d]^{\pi_n} & \Map_{\rm{Pro}(\scr{M})}(P^{\leq \bullet}(\overline E_{\bullet}\wedge X),Y)  \ar[d]^{\pi} \ar[r]\ar[l]_(0.47){\simeq} & \Map(X,Y)  \\
         \colim_n \Map(\overline E_n\wedge X,Y) &  \Map_{\rm{Pro}(\scr{M})}(\overline E_{\bullet}\wedge X,Y)\ar[ur]\ar[l]_(0.47){\simeq} , &\\
        } 
      \end{equation}
      where the map $\pi$ is induced by the projection maps to the Postnikov truncations. The facts that $Y$ is $E$-nilpotent (\ref{lemma:S-nilp_implies_E-nilp}) and that the standard Adams tower $X\ra\overline E_\bullet \wedge X$ is a $E$-nilpotent resolution (\ref{prop:E-nilp_resol_are_unique_and_com_the_nilp_comp}), imply that the diagonal map is an equivalence. We are left to show that $\pi$ in an equivalence too. This follows form the fact that $Y$ is bounded in the homotopy $t$-structure on $\scr{M}$, so $\pi_n$ is an equivalence for large $n$.      
    \end{proof}

    \begin{prop}
      \label{prop:S-nilp_resol_are_unique_and_com_the_nilp_comp}
      For every object $X \in \scr{M}_{\geq k}$, where $k\in \bb Z$, the following holds.
      \begin{enumerate}
        \item \label{SNR1} The Postnikov truncation of the standard Adams tower $X \ra P^\bullet(\overline{E_\bullet}\wedge X)$ is a strongly $\cal R$-nilpotent resolution of $X$ in $\scr{M}$.
        \item \label{SNR2} The pro-object under $X$ associated to a strongly $\cal R$-nilpotent resolution of $X$ is unique up to a contractible space of choices.
        \item \label{SNR3} If $X \ra X_\bullet$ is a strongly $\cal R$-nilpotent resolution of $X$, there is a natural equivalence $\lim_{N(\bf N^{op})} X_\bullet \simeq X^{\wedge}_{E}$ in $\scr{M}_{X/}$.
      \end{enumerate}
    \end{prop}
    \begin{proof}
      Point \eqref{SNR1} is \ref{prop:PnEn_tower_is_an_S-nilp_resol}. Point \eqref{SNR2} works similarly to second point of \ref{prop:E-nilp_resol_are_unique_and_com_the_nilp_comp}. For point \eqref{SNR3} we combine the equivalence 
      $\{X_\bullet\}\simeq \{P^\bullet(\overline{E_\bullet}\wedge X)\}$ of pro-objects under $X$ coming from the previous points with the $\tau$-equivalence (c.f. \ref{defin:pro-we}) $\{P^\bullet(\overline{E_\bullet}\wedge X)\} \leftarrow \{\overline{E_\bullet}\wedge X\}$ of pro-objects under $X$ induced by the projection to the tower of truncations. The inverse limit of a composition of these two maps induces an equivalence under $X$ as required in the statement.
    \end{proof}
    
    \begin{lemma}
    \label{lemma:PnMoore_tower_is_an_S-nilp_resol}
      Consider a collection of tif objects $L_1,\dots, L_r$ of $\scr{M}$. For every $i=1,\dots,r$ let $f_i: L_i \ra\tau_0(\bf 1)$ be a map in $\scr{M}$, and let $\cal R$ be the commutative algebra in $\scr{M}^{\heartsuit}$ defined by $\cal R=\tau_0(\bf 1)/(f_1,\dots,f_r)$. Then for every integer $k$ and every $X \in \scr{M}_{\geq k}$, the tower $ X\ra P^\bullet(C(f_1^\bullet)\wedge \cdots \wedge C(f_r^\bullet)\wedge X)$ is a strongly $\cal R$-nilpotent resolution of $X$. 
    \end{lemma}

    \begin{proof}
      We need to check that for every pair of integers $k,n$ the homotopy objects $\tau_k(C(f_1^n)\wedge \cdots \wedge C(f_r^n)\wedge X)$ are $\cal R$-nilpotent. We accomplish this by induction on $n$, the base case being that of $n=1$. Assume thus that $n=1$: we proceed by induction on $r$. When $r=0$, $\cal R=\tau_0(\bf 1)$, and thus the homotopy objects $\tau_k(X)$ are $\tau_0(\bf1)$-nilpotent. When $r\geq 1$, we set $\cal R_a:=\tau_0(\bf 1)/(f_1,\dots,f_{a})$.
      We can apply Lemma \ref{lemma:mul_by_f_in_a_B_nilp_A-mod} to the homotopy objects $\cal M := \tau_k (C(f_1)\wedge \cdots \wedge C(f_a)\wedge X)$ and $\cal M^{'}:=\tau_{k-1}(C(f_1)\wedge \cdots \wedge C(f_a)\wedge X)$, which are $\cal R_a$-nilpotent by the inductive assumption. It follows that the external homotopy objects of the exact sequences (for varying $k$)
        \begin{equation}
          \label{eqn:}
          0 \ra \coker (r_{f_{a+1}}(\cal M)) \ra \tau_k(C(f_1)\wedge \cdots \wedge C(f_{a+1})\wedge X)) \ra \ker (r_{f_{a+1}}(\cal M')) \ra 0 
        \end{equation}
       are strongly $\cal R_a/(f_{a+1})$-nilpotent. In particular they are $\cal R_{a+1}$-nilpotent, since $\cal R_a/(f_{a+1})\simeq \cal R_{a+1}$, and by Lemma \ref{lemma:S-nilp_in_sequences-extensions} we conclude that the central homotopy objects are $\cal R_{a+1}$-nilpotent too. Given now any $r$-tuple of positive integers $(n_1,\dots,n_r)$ we can show that $\tau_k(C(f_1^{n_1})\wedge \cdots \wedge C(f_r^{n_r})\wedge X)$ is $\cal R$-nilpotent by induction, using the fiber sequences
      \begin{equation}
        L_i^{\wedge n_i-1}\wedge C(f_i) \ra C(f_i^{n_i}) \ra C(f_i^{n_i-1}) 
      \end{equation}
      for $i=1,\dots,r$ and Lemma \ref{lemma:S-nilp_in_sequences} to reduce to the above base case.

      Let us assume now that $Y$ is a strongly $\cal R$-nilpotent object of $\scr{M}$. Let us denote by $C_n$ the object $C(f_1^n)\wedge \cdots \wedge C(f_r^n) \in \scr{M}$. We thus need to check that the map 
      \[ \colim_{N(\bf N^{op})} \Map \big (P^n(C_n\wedge X), Y\big )\ra \Map(X,Y),\]
      induced by the projection to the tower, is an equivalence. Since $\tau_k(Y)$ is strongly $\cal R$-nilpotent then each of the $f_i$'s acts on such a homotopy module nilpotently. It follows that there is an integer $N$ for which each of the $f^N_i$'s acts by $0$ on each the homotopy object of $Y$. Moreover, since $Y$ is bounded in the $t$-structure, up to enlarging $N$ we may actually assume that each of the $f^k_i$'s acts by $0$ on $Y$ whenever $k\geq N$. From this we get that the projection to the tower induces an equivalence
      \[\xymatrix{\Map \big (P^k(C_k\wedge X),Y\big )\ar[r]^(.6){\simeq} & \Map(X,Y),\\}\]
      for $k$ large. 
    \end{proof}

    \begin{thm}
      \label{thm:convergence_mod_stuff}
      Let $E$ be a homotopy commutative algebra of $\scr{M}$ satisfying assumption \ref{sub:assumption_A1} in the special case that $J=\emptyset$. Then for every $k$-connected object of $\scr{M}$ the natural map $\alpha_E(X): X \ra X^{\wedge}_E$ is an $E$-equivalence. In particular the map $\beta_E(X): X_E\ra X^{\wedge}_E$ of \eqref{eqn:alphabetagamma} is an equivalence in $\scr{M}$. 
    \end{thm}

    \begin{proof}
      Let $C_n:=C(f_1^n)\wedge  \cdots\wedge  C(f_r^n)$. Proposition \ref{prop:PnEn_tower_is_an_S-nilp_resol} and Lemma \ref{lemma:PnMoore_tower_is_an_S-nilp_resol} imply that both the tower $P^\bullet(\overline E_\bullet\wedge  X)$ and the tower $P^\bullet(C_\bullet\wedge  X)$ are strongly $\cal R$-nilpotent resolutions of $X$. Thanks to Proposition \ref{prop:S-nilp_resol_are_unique_and_com_the_nilp_comp} there is an equivalence $\varphi$ of their associated pro-objects under $X$. As a consequence we have an induced equivalence under $X$
      \[\psi: \lim_{N(\bf N^{op})} P^\bullet(\overline E_\bullet\wedge  X) \overset{\simeq}{\ra}  \lim_{N(\bf N^{op})} P^\bullet(C_\bullet\wedge  X)\]
      between their limits. Since these limits are naturally identified with $X^{\wedge}_E$ and $X^{\wedge}_{f_1,\dots,f_r}$ respectively, we get a commutative square of pro-objects under $X$
      \begin{equation}
        \label{eqn:E-loc_equals_E-nilp:first_sq}
        \xymatrix{
        \{X^{\wedge}_E\} \ar[r]^{\psi}_{\simeq} \ar[d] & \{X^\wedge_{f_1,\dots,f_r}\} \ar[d] \\
        \{P^{\leq \bullet}(\overline E_\bullet\wedge  X)\}  \ar[r]^{\phi} & \{P^\bullet(C_\bullet\wedge  X)\}, \\
        }
      \end{equation} 
      where the vertical maps are the natural projections. 
      
      After smashing the towers in \eqref{eqn:E-loc_equals_E-nilp:first_sq} with $E$ we obtain a new commutative square 
        \begin{equation}
          \label{eqn:E-loc_equals_E-nilp:first_sq_E}
          \xymatrix{
           \{E\wedge  X^{\wedge}_E\} \ar[r] \ar[d] &  \{E\wedge  X^\wedge_{f_1,\dots,f_r}\} \ar[d]\\
           \{E\wedge  P^\bullet(\overline E_\bullet\wedge  X)\}  \ar[r] & \{E\wedge  P^\bullet(C_\bullet\wedge  X)\}.\\
          }
        \end{equation}
        of pro-objects. Here the lower horizontal map remains an equivalence of pro-objects. We claim that the right vertical map of \eqref{eqn:E-loc_equals_E-nilp:first_sq_E} is a $\tau$-equivalence of pro-objects. 

        For showing this claim, consider that the vertical map on the right hand side of \eqref{eqn:E-loc_equals_E-nilp:first_sq}
        \begin{equation}
          \label{eqn:E-loc_eq_E-nilp:first_sq_rhs}
          X^\wedge_{f_1,\dots,f_r} \ra  \{P^\bullet(C_\bullet\wedge  X)\}
        \end{equation}
        factors as the composition of two maps: the projection 
        \begin{equation}
          \label{eqn:E-loc_eq_E-nilp:first_sq_proj}
          \{X^{\wedge}_{f_1,\dots,f_r}\}  \ra  \{C_\bullet\wedge  X\},
        \end{equation}
        and the projection to the Postnikov tower 
        \begin{equation}
          \label{eqn:E-loc_eq_E-nilp:first_sq_pt}
          \{C_\bullet\wedge  X \}\ra \{P^\bullet(C_\bullet\wedge  X)\}.
        \end{equation}
        The map \eqref{eqn:E-loc_eq_E-nilp:first_sq_pt} is a $\tau$-equivalence by \ref{subs:proj_to_post_is_pwe} and stays a $\tau$-equivalence after smashing with $E$ by \ref{lemma:E_wedge_preserves_pro_equiv}. 
        Concerning \eqref{eqn:E-loc_eq_E-nilp:first_sq_proj}, we observe that by construction we have a commutative diagram of pro-objects,
        \begin{equation}
          \label{eqn:E-loc_eq_E-nilp:first_sq_proj_dec}
        \xymatrix{
            X^{\wedge}_{f_1,\dots,f_r}  \ar[r] &  \{C_n\wedge  X\}  \\
            & X. \ar[ul]_{\chi} \ar[u] \\ 
           }
        \end{equation}
        Since $\chi$ is an $E$-equivalence by \ref{prop:fund_ineq_of_loc}, in order to finish the proof of the above claim we are left to show that the vertical map of \eqref{eqn:E-loc_eq_E-nilp:first_sq_proj_dec} induces, after smashing with $E$, a $\tau$-equivalence of pro-objects in $\scr{M}$.

        To accomplish this task we consider the tower $F^{(r)}_\bullet=\rm{fib}(X \ra C_\bullet \wedge X)$, and show that the pro-object $\{E\wedge  F^{(r)}_\bullet\}$ is $\tau$-equivalent to $0$.
        For this, remember that in the tower $C_\bullet\wedge  X$ the transition maps
        \[\psi_{n+1}: C_{n+1}=C(f_1^{n+1})\wedge  \cdots\wedge  C(f_r^{n+1}) \ra C(f_1^{n})\wedge  \cdots\wedge  C(f_r^{n})=C_n \] 
        are defined as $\psi_n=p_{n+1}(f_r)\circ \cdots \circ p_{n+1}(f_1)$: here for every integer $n$ the maps $p_{n}(f_i)$ are induced by the maps $p_n$ defined in \ref{const:modx_1_x_n_Moore_sp} and displayed in \eqref{eqn:mod_x^n_mod_x^n-1_fib_seqs}. 
        If $r=1$ we have $F^{(1)}_\bullet=\bf {r_{f_1}}(L^{\wedge n-1})$ and its transition maps are $r_{f_1}(L^{\wedge n-1})$, as in \eqref{eqn:mod_x^n_mod_x^n-1_fib_seqs}. Since $r_{(f_1)}$ acts trivially on the towers of homotopy objects $\tau_k(E\wedge  F^{(1)}_\bullet)$, we deduce that $\{E\wedge  F^{(1)}_\bullet\}$ is $\tau$-equivalent to $0$. If $r>1$ one can argue by induction. Indeed, using the octahedral axiom, we can find a fiber sequence of towers
        \[ F^{(s-1)}_\bullet \ra F^{(s)}_\bullet \ra G^{(s)}_\bullet  \]
        having the following properties: the pro-object $\{E\wedge  F^{(s-1)}_\bullet\}$ is $\tau$-equivalent to $0$ by the inductive assumption, and $\{E\wedge  G^{(s)}_\bullet\}$ is $\tau$-equivalent to $0$ by the case $r=1$ treated above.
        Hence, thanks to \ref{cor:2-3prozero}, we conclude that \eqref{eqn:E-loc_eq_E-nilp:first_sq_proj} becomes a $\tau$-equivalence after smashing with $E$.
        
        Let us consider now the commutative square of towers:
         \begin{equation}
          \label{eqn:E-loc_equals_E-nilp:second_sq}
          \xymatrix{
           X  \ar[r] \ar[d] &  X^{\wedge}_E \ar[d] & \\
           \overline E_\bullet\wedge  X \ar[r] &  P^\bullet(\overline E_\bullet\wedge  X)  \\
          }
        \end{equation}
        and note that, in order to conclude, we only need to show that the map $E\wedge  \alpha_E(X) :E\wedge  X\ra E\wedge  X^{\wedge}_E$ is an equivalence. This map is the upper horizontal arrow of the diagram 
        \begin{equation}
          \label{eqn:E-loc_equals_E-nilp:second_sq_E}
          \xymatrix{
           E\wedge  X  \ar[r] \ar[d] & E\wedge  X^{\wedge}_E  \ar[d] & \\
           E\wedge  (\overline E_\bullet\wedge  X)\ar[r] &  E\wedge  P^\bullet(\overline E_\bullet\wedge  X) \\
          }
        \end{equation}
        which is obtained from \eqref{eqn:E-loc_equals_E-nilp:second_sq} by smashing with $E$. In \eqref{eqn:E-loc_equals_E-nilp:second_sq_E} the right vertical map induces a $\tau$-equivalence on pro-objects as it follows from the previous part of the argument. The lower horizontal map induces as well a $\tau$-equivalence of the associated pro-objects by \ref{lemma:E_wedge_preserves_pro_equiv}. Finally the left vertical map of \eqref{eqn:E-loc_equals_E-nilp:second_sq_E} induces a $\tau$-equivalence too, as we now explain. We have a fiber sequence of towers
          \begin{equation}
          \label{eqn:bo}
            \overline E^\bullet \wedge  X \ra X \ra \overline E_{\bullet-1}\wedge  X
          \end{equation}
        which we obtain from diagram \eqref{eqn:fib_seq_of_Ad_towers} upon smashing with $X$. The left vertical map in the square \eqref{eqn:E-loc_equals_E-nilp:second_sq} is the map induced by the right hand side maps of \eqref{eqn:bo}. We claim that, after smashing \eqref{eqn:bo} with $E$, the tower $ \overline E^\bullet \wedge  X$ on the left hand side of \eqref{eqn:bo} becomes $\tau$-equivalent to zero. Indeed by the very inductive definition of $\overline E^n$ we have fiber sequences deduced from \eqref{eqn:fund_fib_seq_E^n+1_E^n_EE^n}:
        \[\xymatrix{ \overline E^{n+1}\wedge  X \ar[r]^{\bar e\wedge \rm{id}}  & \overline E^n\wedge  X \ar[r]^(.45){e_E\wedge\rm{id}} & E\wedge  \overline E^n\wedge  X.}\] 
        Here $e_E: \bf 1 \ra E$ is the unit of the algebra $E$, while $\bar e\wedge  \op{id}$ (see \eqref{eqn:fund_fib_seq_E^n+1_E^n_EE^n}) appears as the transition map in the tower $E\wedge  \overline E^\bullet \wedge  X$.
        After smashing with $E$ we thus have an induced long exact sequences of homotopy objects
        \[ \cdots \ra \tau_k(E\wedge  \overline E^{n+1}\wedge  X) \ra \tau_k(E\wedge  \overline E^n\wedge  X) \ra \tau_k(E\wedge  E\wedge  \overline E^n\wedge  X) \ra \cdots\]
        and the maps $\tau_k(E\wedge  \overline E^n\wedge  X) \ra \tau_k(E\wedge  E\wedge  \overline E^n\wedge  X)$ are split by the multiplication of $E$. In this case the previous map in the long exact sequence, which is the same as the transition map in the tower $\tau_k(E\wedge  \overline E^\bullet\wedge  X)$ is zero, and hence the associated pro-object is equivalent to $0$ for every $k$. Using Corollary \ref{cor:2-3prozero}, we deduce that the upper horizontal map of \eqref{eqn:E-loc_equals_E-nilp:second_sq_E} is a $\tau$-equivalence. Both source and target of this map are constant towers so the map is actually an equivalence in $\scr{M}$, and this concludes the proof.
    \end{proof}

    \begin{rmk}
      \label{rmk:HKOcomparison}
      Observe that the combination of Theorem \ref{thm:convergence_mod_stuff}, Theorem \ref{thm:red_to_moore_spt} and Proposition \ref{prop:localization_at_mod_x_1-x_n_moore_spectrum} recovers and generalizes the results of \cite{MR2811716} using different techniques. Similar results have appeared \cite{topologmodels}.
      \end{rmk}

    \begin{lemma}
      \label{lem:Pn_inverted_Moore_tow_is_S_nilp_resol}
      Let $E$ be a homotopy commutative algebra in $\scr{M}$ satisfying Assumption \ref{sub:assumption_A1} in the special case that $I=\emptyset$. Then for every $k$-connected object $X$, the tower $\{P^n(\bf 1[\cal J^{-1}]\wedge  X)\}_n$ is a strongly $\tau_0E$-nilpotent resolution of $X$.
    \end{lemma}

    \begin{proof}
      Since the unit $\bf 1 \ra E$ induces an equivalence $\tau_0(\bf 1[\cal J^{-1}]) \overset{\simeq}{\ra} \cal R$ we immediately conclude that the homotopy objects $\tau_k(X\wedge  \bf 1[\cal J^{-1}])$ are all $\cal R$-modules and hence they are strongly $\cal R$-nilpotent. The mapping property of strongly $\cal R$-nilpotent resolutions is an immediate consequence of the universal property of the Postnikov truncations.
    \end{proof}

    \begin{thm}
      \label{thm:convergence_inv_stuff}
      Let $E$ be an homotopy commutative algebra in $\scr{M}$ satisfying assumption \ref{sub:assumption_A1} in the special case that $I=\emptyset$. Then for every $k$-connected object $X$, the natural map $\lambda_E(X): X \ra X^{\wedge}_E$ is a $E$-equivalence. In particular the map $\beta_E(X): X_E\ra X^{\wedge}_E$ of \eqref{eqn:alphabetagamma} is an equivalence in $\scr{M}$. 
    \end{thm}

    \begin{proof}
      The proof proceeds along the same lines as the proof of \ref{thm:convergence_mod_stuff}. More precisely we start by observing that both the towers $P^n(\overline E_n\wedge  X)$ and $P^n(X\wedge  \bf 1[\cal J^{-1}])$ are both strongly $\cal R$-nilpotent resolutions of $X$ by \ref{prop:PnEn_tower_is_an_S-nilp_resol} and \ref{lem:Pn_inverted_Moore_tow_is_S_nilp_resol} respectively. We deduce, as in the proof of \ref{thm:convergence_mod_stuff}, that there is an equivalence $\psi: X^\wedge_E \ra X\wedge  \bf 1[\cal J^{-1}]$ under $X$ making the following square of pro-objects under $X$
      \begin{equation}
          \label{eqn:Inv-E-loc_equals_E-nilp:first_sq}
          \xymatrix{
          \{ X^{\wedge}_E\} \ar[r]^(0.45)u \ar[d] & \{ X\wedge  \bf 1[\cal J^{-1}]\} \ar[d]\\
          \{ P^\bullet(\overline E_\bullet\wedge  X) \} \ar[r]^(0.45){\simeq}_(0.45)\varphi & \{ P^\bullet(X\wedge  \bf 1[\cal J^{-1}]) \}\\
          }
        \end{equation}
      commutative. After smashing \eqref{eqn:Inv-E-loc_equals_E-nilp:first_sq} with $E$ the previous diagram, the lower horizontal map remains an equivalence of pro-objects. The vertical map on the right hand side of \eqref{eqn:Inv-E-loc_equals_E-nilp:first_sq} is a $\tau$-equivalence by \ref{subs:proj_to_post_is_pwe}, and stays a $\tau$-equivalence after smashing with $E$ by \ref{lemma:E_wedge_preserves_pro_equiv}. These two observations show that, after smashing with $E$, also the left vertical map of \eqref{eqn:Inv-E-loc_equals_E-nilp:first_sq} is a $\tau$-equivalence. The remaining part of the proof follows step by step the proof of Theorem \ref{thm:convergence_mod_stuff}.
    \end{proof}

\appendix

  \section{Pro-Objects} 
  \label{sec:pro_spectra}
  \stepcounter{subsection}
  We gather here some well known statements on pro-objects in an $\infty$-category $\scr C$ that we have freely used in the previous sections. We will also recall some properties of pro-objects in a symmetric monoidal stable $\infty$-category $\scr{M}$ endowed with a left-complete multiplicative $t$-structure.

  \subsubsection{} Recall that the $\infty$-category of pro-objects of $\scr C$, denoted by $\rm{Pro}(\scr C)$ is defined as $\rm{Ind}(\scr C^{op})^{op}$ (c.f Section 5.3.5 of \cite{HTT}). Recall that every diagram $p:I \ra \scr C$ indexed on a co-filtered simplicial set $I$ gives a pro-object of $\scr C$: in $\scr C^{op}$ we have a filtered diagram $I^{op} \ra \scr C^{op}$, and the colimit of its composition with the Yoneda embedding 
  \[\xymatrix{I^{op}\ar[r]^p & \scr C^{op} \ar[r]^j & \scr P(\scr C^{op}),}\]  is the desired pro-object.
  Starting with Section \ref{sec:moore_spectra} we have considered towers in $\scr C$, i.e. diagrams indexed on $\bf N^{op}$ or $K^{op}=\cup_n {K^n}^{op}$: these are clearly co-filtered simplicial sets, and thus every diagram in $\scr C$ indexed on them gives rise to a pro-object of $\scr C$. Along the text we have used the symbol $\{X_\bullet\}$ to denote the pro-object associated with a diagram $X_\bullet$.

  \begin{lemma}
    \label{lemma:map_of_probjects}
     Let $X_\bullet, Y_\bullet:I \ra \scr C$ be diagrams in $\scr C$. Then for the associated pro-objects,
     \[\Map_{\rm{Pro}(\scr C)}(\{X_\bullet\},\{Y_\bullet\})\simeq \lim_n\colim_k\Map_{\scr C}(X_k, Y_n).\]
   \end{lemma} 
  \begin{proof}
    Follows from the fact $\scr P(\scr C^{op})$ admits small colimits (Corollary 5.1.2.3 of \cite{HTT}) and constant ind-objects are compact (Proposition 5.3.5.5 of \cite{HTT}).
  \end{proof}

  \begin{defin}
  \label{defin:pro-we}
  	Assume that $\scr M$ is a stable $\infty$-category endowed with a left-complete multiplicative $t$-structure.
    A map of pro-objects $f:\{X_\bullet\} \ra \{Y_\bullet\}$ is a \emph{$\tau$-equivalence} if, for every integer $p$, the induced map $ \{\tau_p (X_\bullet)\} \ra  \{\tau_p (Y_\bullet)\}$ is an equivalence in $\rm{Pro}(\scr{C}^{\heartsuit})$. A pro-object $\{X_\bullet\}$ is \emph{$\tau$-equivalent to $0$} if, for every integer $p$, the homotopy objects $\{\tau_p(X_\bullet)\}$ are equivalent to $0$ as objects of $\rm{Pro}(\scr{C}^{\heartsuit})$.
  \end{defin}
  
  \subsubsection{}
  \label{subs:proj_to_post_is_pwe}
  Of course every equivalence of pro-objects is a $\tau$-equivalence. Moreover, given any pro-object $\{X_\bullet\}$ of $\scr{M}$, the projection to the Postnikov tower 
  \[ \pi_k:X_k \ra P^k(X_k)\]
  induces a $\tau$-equivalence $\{\pi_\bullet\}:\{X_\bullet\} \ra \{ P^\bullet (X_\bullet)\}.$ Note that in general the projection map $\{\pi_\bullet\}$ does not need to be a pro-equivalence. 
  If the $t$-structure on $\scr M$ is left complete then the map $\{X_\bullet\} \ra \{P^\bullet(X_\bullet)\}$ induces an equivalence between the respective inverse limits, as it is easily checked using the Milnor sequence.

  \subsubsection{}
  Since $\scr{M}^{\heartsuit}$ is an abelian category, then $\rm{Pro}(\scr{M}^{\heartsuit})$ is an abelian category by Proposition 4.5 of \cite[Appendix]{MR0245577}. Moreover the full subcategory $\rm{Tow}(\scr{M}^{\heartsuit})\subseteq \rm{Pro}(\scr{M}^{\heartsuit})$ spanned by pro-objects indexed on $\op N(\bf N^{op})$ is closed under finite limits and colimits: this can be proved following the same strategy of proof of Proposition 2.7 in \cite{MR1428551}. In particular $\rm{Tow}(\scr{M}^{\heartsuit})$ is an abelian category. It follows that a map $\{f_\bullet\} :\cal \{M_\bullet\} \ra \cal \{N_\bullet\}$ of $\rm{Tow}(\scr{M}^{\heartsuit})$ is an equivalence if and only if both $\ker(\{f_\bullet\})$ and $\coker(\{f_\bullet\})$ are pro-objects equivalent to $0$. In particular we conclude the following.
  \begin{cor}
  \label{cor:2-3prozero}
    Let $X_\bullet, Y_\bullet, Z_\bullet$ be towers in $\scr{M}$. Assume we have a fiber sequence of towers
    \[\xymatrix{X_\bullet \ar[r]^{f_\bullet} & Y_\bullet \ar[r]^{g_\bullet} &  Z_\bullet}\]
    Then the induced map of pro-objects $\{g_\bullet\}:\{Y_\bullet\} \ra \{Z_\bullet\}$ is a $\tau$-equivalence if and only if the pro-object $\{X_\bullet\}$ is $\tau$-equivalent to $0$.
  \end{cor}

  \begin{lemma}
  \label{lemma:E_wedge_preserves_pro_equiv}
    Assume $\{ X_\bullet\}$ is an object of $\rm{Tow}(\scr{M})$. Then for every $k$-connected object $E\in \scr{M}$, the projection to the Postnikov tower induces a $\tau$-equivalence 
    \[\{E\wedge X_\bullet\}\ra \{ E\wedge P^\bullet (X_\bullet)\}.\]
  \end{lemma}
  \begin{proof}
      Consider the fundamental fiber sequence
      \[ E\wedge P_{n+1}(X_{n}) \ra E\wedge X_{n}\ra E\wedge P^{n}(X_{n}) \]
      induced from \eqref{eqn:postfibseq}.
      If $E$ is $k$-connected then $E\wedge P_{n+1}(X_{n})$ is $(k+n+1)$-connected, and in particular $\{\tau_p X_\bullet\}$ is equivalent to zero. The statement then readily follows from \ref{cor:2-3prozero}. 
   \end{proof}

\section{Categorical Recollection} 
    \label{sec:cat_rec}
  
  We quickly gather some known general properties of the $\infty$-category $\scr{SH}(S)$, where $S$ is a Noetherian scheme of finite Krull dimension. Using the formalism of \cite{HA} we deduce analogous properties for categories of modules over commutative algebras in $\scr{SH}(S)$. Our arguments are by no mean original; they are rather a reader's guide to the navigation of the relevant statements of \cite{HA}. 

  \subsection{} In \cite{jard:motsym} the author introduces the motivic model structure on symmetric $T$-spectra $\spt_T^\Sigma(S)$, which is the base for our constructions. Section 4, and in particular Theorem 4.15 and 4.31 of \cite{jard:motsym} contain a proof of the following statement.

  \begin{thm}
    \label{thmB:mod_str_on_sym_spt}
    The stable motivic model structure on $\spt_T^\Sigma(S)$ is simplicial stable proper cofibrantly generated and combinatorial. The smash product of motivic symmetric spectra can be completed to the datum of a monoidal structure on $\spt_T^\Sigma(S)$ making the stable motivic model structure a symmetric monoidal model structure. 
  \end{thm}

  \begin{defin}
  \label{defin:sym_mon_infty_cat}
    A symmetric monoidal $\infty$-category is an infinity category $\scr C^{\otimes}$ with a coCartesian fibration $p: \scr C^\otimes \ra \op N (\op{Fin_\ast)}$ such that the functors $\rho^i_!: \scr C^\otimes_{\langle n \rangle} \ra \scr C^\otimes_{\langle 1 \rangle}$ induce and equivalence of $\infty$-categories $\prod_i \rho^i_!: \scr C^\otimes_{\langle n \rangle} \simeq (\scr C^\otimes_{\langle 1 \rangle})^n$ (Definition 2.0.0.7 of \cite{HA}). The $\infty$-category $\scr C^\otimes _{\langle 1 \rangle}$ is called the underlying $\infty$-category associated with $\scr C^\otimes$, and is usually abusively called a symmetric monoidal $\infty$-category, hiding the reference to the map $p$.
  \end{defin}
  \subsubsection{} 
  \label{ssubB:operad_vs_mon_cat}
  As in ordinary category theory coloured operads generalize the notion of symmetric monoidal category, similarly in higher category theory $\infty$-operads generalize the notion of symmetric monoidal $\infty$-categories (see Proposition 2.1.2.12 of \cite{HA}). We do not recall the definition of an $\infty$-operad which can be found in 2.1.1.10 of \cite{HA}. We however record that, as a consequence of Proposition 2.1.2.12 of \cite{HA}, a coCartesian fibration $p: \scr C^\otimes \ra \op N (\op{Fin_\ast)}$ between $\infty$-categories is a $\infty$-operad if and only if it is a symmetric monoidal $\infty$-category.

  \subsubsection{$\scr{SH}(S)$ is symmetric monoidal, presentable and stable} 
  \label{ssubB:SH_is_mon_pres_stab}
  
  In the case of interest for us, the work of Lurie on $\infty$-categories gives us a streamlined way of producing a symmetric monoidal $\infty$-category $\scr{SH}(S)$ associated with a base scheme $S$. One possible construction is to set
  \[\scr{SH}(S):=\op N_{\Delta}(\spt_T^\Sigma(S)^o),\]
  where $\op N_\Delta$ denoted the simplicial nerve construction, and $\spt_T^\Sigma(S)^o\subseteq \spt_T^\Sigma(S)$ denotes the full subcategory spanned by cofibrant-fibrant objects.

  With this definition, using Proposition 1.2.3.9 of \cite{HTT}, a direct check allows to conclude that we have equivalences of homotopy categories 
  \[h \op N_\Delta(\spt_T^\Sigma(S)^o)\simeq \op{Ho} (\spt_T^\Sigma(S)^o)\simeq \op{Ho} (\spt_T^\Sigma(S))\simeq \cal{SH}(S).\]

  By design, $\scr{SH}(S)$ is the fiber over $\langle 1 \rangle$ of the operadic nerve construction 
  \[p: \op{N}^\otimes(\spt_T^\Sigma(S)^o)) \ra \op N (\op{Fin_\ast)}\]
  c.f. Notation 2.1.1.22 and Definition 2.1.1.23 of \cite{HA}), and the map $p$ is in fact a coCartesian fibration (Proposition 4.1.7.15 of \cite{HA}) and an $\infty$-operad \cite[Proposition 2.1.1.27]{HA}. In view of \ref{ssubB:operad_vs_mon_cat} we conclude that $\scr{SH}(S)$ is the underlying $\infty$-category of a symmetric monoidal $\infty$-category. The monoidal product, which is induced by the smash product of spectra is denoted by $-\wedge-$ and the unit is denoted $\bb S$.

  Recall that $\spt_T^\Sigma(S)$ is combinatorial simplicial model category. In this situation to A.3.7.6 of \cite{HTT} implies that $\scr{SH}(S)$ is a presentable $\infty$-category. In particular $\scr{SH}(S)$ admits all small limits and colimits (Proposition 4.2.4.8 of \cite{HTT}). 

  In addition $\scr{SH}(S)$ is a stable $\infty$-category in the sense of Definition 1.1.1.9 of \cite{HA}. Indeed in view of Corollary 1.4.2.27 of \cite{HA}, we only need to check that $\scr{SH}(S)$ has finite colimits and the suspension functor $\Sigma: \scr {SH}(S) \rightarrow \scr{SH}(S)$ is an equivalence. Now the first condition is satisfied since $\scr{SH}(S)$ is presentable. The second condition is implied by the fact that the suspension functor $\Sigma$ is induced from the Quillen equivalence $\Sigma_{\bb S^1} : \spt_T^\Sigma(S) \rightleftarrows \spt_T^\Sigma(S): \Omega_{\bb S^1}$. Indeed from Theorem 4.2.4.1 of \cite{HTT} it follows that $\op N_\Delta$ carries homotopy (co)limit diagrams in $\spt_T^\Sigma(S)^o$ to (co)limit diagrams in $\scr{SH}(S)$. 

  A more informative construction of $\scr{SH}(S)$ was given in Theorem 2.26 of \cite{MR3281141}, where the author proves that $\scr{SH}(S)$ is the stable presentable symmetric monoidal $\infty$-category initial with respect to the property that $\bb P^1\wedge-$ be an invertible endofunctor. 


  \subsubsection{$\scr{SH}(S)$ is presentably symmetric monoidal} 
  \label{ssubB:sh_is_presentable_monoidal}
  
  The properties of the monoidal structure on $\scr{SH}(S)$ are actually a bit stronger than what already stated, making it a presentably symmetric monoidal $\infty$-category according to Definition 3.4.4.1 of \cite{HA}. The essential point is that $\scr{SH}(S)$ is the underlying $\infty$-category of a combinatorial simplicial symmetric monoidal model category. Indeed: 
  \begin{enumerate}
    \item $p$ is a coCartesian fibration of $\infty$-operads, as recorded in \ref{ssubB:SH_is_mon_pres_stab};
    \item $p$ is compatible with small colimits according to Definition 3.1.1.18 of \cite{HA}. Indeed $\spt_T^\Sigma(S)$ is co-complete and the smash product of spectra commutes with all colimits, being a left adjoint.
    \item For each object $\bc n \in N (\op{Fin}_\ast) $, the fiber $\scr{SH}(S)^{\otimes}_{<n>}$ is presentable. Indeed 
    \[\scr{SH}(S)^{\otimes}_{<n>}\simeq \prod_{i=1}^n \scr{SH}(S)^{\otimes}_{<1>}\simeq \prod_{i=1}^n\scr{SH}(S),\] whose presentability of was addressed in \ref{ssubB:SH_is_mon_pres_stab}.
  \end{enumerate}

  \begin{defin}[\cite{HA} 2.1.3.1]
  \label{definB:com_alg}
    A commutative algebra of $\scr{SH}(S)$ is a map of $\infty$-operads $s: \op N(\op{Fin}_\ast) \rightarrow \scr{SH}(S)^{\otimes}$ which is a section of the natural projection $p: \scr{SH}(S)^{\otimes} \rightarrow \op N(\op{Fin}_\ast)$. The object $s(\bc 1): \Delta \rightarrow \scr {SH}(S)$ is called the commutative algebra underlying $s$. Along this section we will often refer to $s (\langle 1 \rangle)$ as commutative algebra in $\scr{SH}(S)$, living the section $s$ implicit.
  \end{defin}

  \subsubsection{Modules over a commutative algebra} 
  \label{ssubB:infty_mods}

  Let us denote $\scr{SH}(S)$ simply by $\scr C$ in order to make references to \cite{HA} more directly traceable. Let $A$ be a commutative algebra in $\scr C$. Lurie constructs in Section 3.3.3 of \cite{HA} a fibration of $\infty$-operads $p_A: \rm{Mod}_A(\scr C)^\otimes \rightarrow \op N(\op{Fin}_\ast)$. What in \ref{ssub:modules_over_ring_spectra} we denoted by $\scr{Mod}_A(S)$ is the $\infty$-category $\rm{Mod}_A(\scr C):=\rm{Mod}_A(\scr C)^\otimes_{\langle 1 \rangle}$. We are omitting the relevant $\infty$-operad from the notation, since we are only dealing with $\op N(\op{Fin}_\ast)$. The construction of $\rm{Mod}_A(\scr C)^\otimes$ also produces a forgetful functor $\phi: \rm{Mod}_A(\scr C)^\otimes \rightarrow \scr C^\otimes$. In \ref{ssub:modules_over_ring_spectra} we have denoted by $U_A=\phi_{\langle 1\rangle}$ the functor induced by $\phi$ on the underlying $\infty$-categories.
  
  \subsubsection{$U_A$ commutes with small limits and colimits}
  \label{ssubB:inf_mods_are_bicomplete}
  We can use Theorem 3.4.4.2 or more precisely Corollary 3.4.4.6 of \cite{HA} to deduce that $p_A: \rm{Mod}_A(\scr C)^\otimes  \rightarrow \op N(\op{Fin}_\ast)$ is actually a presentably-symmetric-monoidal $\infty$-category and that the forgetful functor $\phi: \rm{Mod}_A(\scr C)^\otimes \rightarrow \scr C^\otimes$ detects and commutes with all small colimits. Similarly using Corollary 3.4.3.6 of \cite{HA} we deduce that $\phi$ detects and commutes with all small limits. 
  In particular the underlying functor $U_A=\phi_{\langle 1 \rangle}: \rm{Mod}_A(\scr C) \rightarrow \scr C$ detects and commutes with all small limits and colimits.

  \subsubsection{$U_A$ is conservative}
  This is not spelled out explicitly in \cite{HA}, but it can be easily deduced from the combination of some of the main statements. In first place we use that forgetting the commutativity of $A$ induces a canonical equivalence between $A$-modules and left modules over the associative algebra $A$
  \[\rm{Mod}_A(\scr C) \overset{\simeq}{\ra} \rm{LMod}_A(\scr C),\]
  as proven in Corollary 4.5.1.6 of \cite{HA}. In second place we use that $\rm{Mod}_A(\scr C)$ is naturally identified as the fiber over $A$ of a cartesian fibration $\vartheta: \rm{Mod}(\scr C) \rightarrow \rm{Alg}(\scr C)$. More precisely we have a diagram of $\infty$-categories:
  \begin{equation}
    \xymatrix{
    \rm{LMod}(\scr C) \ar[r]^{U} \ar[d]^{\vartheta} & \scr C\\
    \rm{Alg}(\scr C). & \\
    }
  \end{equation}
  The behavior of $\vartheta$ on objects is mapping $(R,M) \mapsto R$, where $R$ is an associative algebra and $M$ a left $R$-module. The functor $\eta$ instead operates on objects as $(R,M)\mapsto M$.
  The functor $\vartheta$ is a cartesian fibration of $\infty$-categories and an arrow of $\rm{LMod}(\scr C)$ is $\vartheta$-cartesian if and only if its image in $\scr C$ is an equivalence (see. 4.2.3.2 of \cite{HA}). However, an $A$-module map $a: \Delta^1 \rightarrow \rm{LMod}_A(\scr C)$ is a map of $\rm{LMod}(\scr C)$ covering $\op{id}_A: \Delta^1 \rightarrow \rm{Alg}(\scr C)$, and thus $a$ is $\theta$-cartesian if and only if $a$ is an equivalence in $\rm{LMod}(\scr C)$ (see for instance Proposition 2.4.4.3 of \cite{HTT}) if and only if $a$ is an equivalence in $\rm{LMod}_A(\scr C)$. Hence $\phi$ is conservative.

  \subsubsection{The left adjoint of $U_A$}
  We keep using the identification mentioned above $\rm{Mod}_A(\scr C) \overset{\simeq}{\ra} \rm{LMod}_A(\scr C)$. The functor $U_A: \rm{LMod}_A(\scr C) \rightarrow \scr C$ has a left adjoint $F_A$, that on objects acts by mapping a spectrum $X$ to the free $A$-module generated by $X$. The composition $U_A\circ F_A\simeq A\wedge -$. Moreover if $\lambda: F_A(X)\overset{\simeq}{\ra} M$ is a free $A$-module, $\lambda$ induces an equivalence of spaces
  \[ \Map_A(M,N) \ra \Map_{\bb S}(X,U_A(N)).\]
  All these claims can be found in Corollary 4.2.4.6 and 4.2.4.8 of \cite{HA}.

  The functor $\phi: \rm{Mod}_A(\scr C)^\otimes \rightarrow \scr C^\otimes$ has a symmetric monoidal left adjoint $\psi: \scr{C}^\otimes \rightarrow \rm{Mod}_A(\scr C)^\otimes$. This follows from Theorem 4.5.3.1 and Remark 4.5.3.2 of \cite{HA}, whose assumptions are automatically satisfied, since $p_A$ is a presentably-symmetric-monoidal $\infty$-category. In particular $F_A$ is equivalent to $\psi_{\langle 1\rangle}$.

      
\printbibliography

\end{document}